\documentclass{ucthesis}

\usepackage{latexsym, amsmath, amssymb, amsfonts, amscd}
\usepackage{amsthm}
\usepackage{t1enc}
\usepackage[mathscr]{eucal}
\usepackage{indentfirst}
\usepackage{graphicx, pb-diagram}
\usepackage{fancyhdr}
\usepackage{fancybox}
\usepackage{enumerate}
\usepackage[all]{xy}

\theoremstyle{plain}
\newtheorem{thm}{Theorem}[section]

\newtheorem{prop}[thm]{Proposition}
\newtheorem{lemma}[thm]{Lemma}
\newtheorem{cor}[thm]{Corollary}
\newtheorem{con}[thm]{Conjecture}

\theoremstyle{definition}
\newtheorem{defi}[thm]{Definition}
\newtheorem{pdef}[thm]{Proposition-Definition}

\theoremstyle{remark}
\newtheorem{remark}[thm]{Remark}
\newtheorem{ep}[thm]{Example}

\newcommand{\lra}{\longrightarrow}

\newcommand{\btd}{\bigtriangledown}

\newcommand{\os}{\overset}

\newcommand{\Z}{\ensuremath{\mathbb Z}}

\newcommand{\R}{\ensuremath{\mathbb R}}

\newcommand{\T}{\ensuremath{\mathbb{T}}}

\newcommand{\g}{\ensuremath{\mathfrak{g}}}

\newcommand{\di}{d}
\newcommand{\ci}{\ensuremath{C^{\infty}}}

\newcommand{\G}{\Sigma}

\newcommand{\call}{{\mathcal L}}
\newcommand{\cC}{\mathcal{C}}            
\newcommand{\cE}{\mathcal{E}}
\newcommand{\cF}{\mathcal{F}}
\newcommand{\cX}{\mathcal{X}}
\newcommand{\cY}{\mathcal{Y}}
\newcommand{\cZ}{\mathcal{Z}}
\newcommand{\cG}{\mathcal{G}}
\newcommand{\cH}{\mathcal{H}}
\newcommand{\cV}{\mathcal{V}}

\DeclareMathOperator{\pr}{pr}

\newcommand{\half}{\frac{1}{2}}

\newcommand{\br}{[\cdot,\cdot ]}

\newcommand{\ta}{\tilde{a}}
\newcommand{\tb}{\tilde{b}}

\newcommand{\tA}{\tilde{A}}

\newcommand{\dds}{\frac{\partial}{\partial s}}

\newcommand{\gpoid}{$G_1 \underset{\mathbf{t}}{ \overset{\mathbf{s}}
{\rightrightarrows}} G_0$ }
\newcommand{\hpoid}{$H_1 \underset{\mathbf{t}}{ \overset{\mathbf{s}}
{\rightrightarrows}} H_0$ }

\newcommand{\gammapoid}{$\Gamma \underset{\mathbf{t}_1}{ \overset{\mathbf{s}_1}
{\rightrightarrows}} P$ }
\newcommand{\hgammapoid}{$\Gamma^h \underset{\mathbf{t}_1}{ \overset{\mathbf{s}_1}
{\rightrightarrows}} P$ }
\newcommand{\bt}{\mathbf{t}}                  
\newcommand{\bs}{\mathbf{s}}                  
\newcommand{\bbt}{\bar{\mathbf{t}}}           
\newcommand{\bbs}{\bar{\mathbf{s}}}           
\newcommand{\bm}{\bar{m}}                     

\begin{document}

\title{Integrating Lie algebroids via stacks and  applications to Jacobi manifolds}
\author{Chenchang Zhu}
\degreesemester{Spring} \degreeyear{2004} \degree{Doctor of
Philosophy} \chair{Professor Alan Weinstein}
\othermembers{Professor Allen Knutson \\
  Professor Hitoshi
Murayama }
\numberofmembers{3} \prevdegrees{BS (Peking University)
1999} \field{Mathematics} \campus{Berkeley}




\maketitle
\approvalpage
\copyrightpage
\pagestyle{fancyplain}
\begin{abstract}

Lie algebroids can not always be integrated into Lie groupoids. We
introduce a new object---``Weinstein groupoid'', which is a
differentiable stack with groupoid-like axioms. With it, we have
solved the integration problem of Lie algebroids. It turns out
that every Weinstein groupoid has a Lie algebroid, and every Lie
algebroid can be integrated into a Weinstein groupoid.

Furthermore, we apply this general result to Jacobi manifolds and
construct contact groupoids for Jacobi manifolds. There are
further applications in prequantization and integrability of
Poisson bivectors.

\abstractsignature
\end{abstract}

\pagestyle{plain}
\begin{frontmatter}
\pagenumbering{roman}
\begin{dedication}
\null\vfil
{\large
\begin{center}
To my dear grandma and grandpa.
\end{center}}
\vfil\null
\end{dedication}

\tableofcontents
\lhead[\fancyplain{}{}]{\fancyplain{}{}}
\rhead[\fancyplain{}{}]{\fancyplain{}{}}
\begin{acknowledgements}

I want to express my deepest thanks to my advisor Alan Weinstein.
Under his gentle mathematical guidance, I was able to enter this
particular research area, and I was able to appreciate the chance
to do math every day!! His support is felt at all time, and I
admire him not only as a mathematician but also as a great person!

Allow me to go back to my far-away young years and express my thanks
to those teachers who have selflessly helped me. Ms Ruan, she was
such a brave and wise teacher. It was she who asked me to work out
five mathematics problems every day, and gave them to her
boyfriend (then, husband now) to correct. That initiated my interest
toward math. Also, Ms Di, my math teacher at that time, her trust in my
humble ability gave me a lot of
encouragement. Mr. Liu, he was really a wonderful math teacher,
and to me, also a poet in his deep mind. It was he and Mr. Qian
who guided me through the math Olympiads. The support of both of them
was much further  beyond mathematics. Without them, I wouldn't have achieved
that much in math competitions.  Prof. Liu and Prof. Qian in
Peking University, when they call me ``the buds of the mother
land'', I feel so much deep care from previous generation of
the wonderful Chinese mathematicians. Prof. Qi in Wuhan
University, he is such a great person, who loves math, music,
paintings, and almost everything! I was so lucky to have the
chance to meet all of them, and learned from all of them. Without
them, my life today would be different. One day, I wish I could
spread my love towards my students like them.

Thanks also go to my collaborators Henriques Burzstyn, Marius Crainic,
Hsian-Hua Tseng and Marco Zambon. I have learned  many
things from them. By talking to them,  I am able to do math so happily!

Apart from the above people, for the content in this thesis,
I would also like to thank Kai Behrend, David Farris, Tom Graber,
 Andr\'e Henriques, Eli Lebow, Joel Kamnitzer,
 David Metzler, Ieke Moerdijk, Janez Mr{\v{c}}un, Ping Xu   for very
 helpful discussions and suggestions.

I would also like to thank my dear lovely friends at Berkeley, whose
names are so many.  I hope they would forgive me for not listing them
here.

Finally, I devote my endless thanks to my dear parents, my
grandparents and  my aunts. Without them, there would be no me...

\end{acknowledgements}

\end{frontmatter}

\pagestyle{fancyplain}
\pagenumbering{arabic}


\chapter{Introduction}
Integrating Lie algebroids is a
long-standing problem: unlike  (finite dimensional) Lie
algebras\footnote{Non-integrability already appears in the case of
infinite dimensional Lie algebras \cite{dl}. In this paper, Lie
algebroids are assumed to be finite dimensional. }  which always have their
associated Lie groups, Lie algebroids do not always have their
associated Lie groupoids \cite{am1} \cite{am2}. So the
Lie algebroid version of Lie's third theorem poses the question
indicated by the following chart:

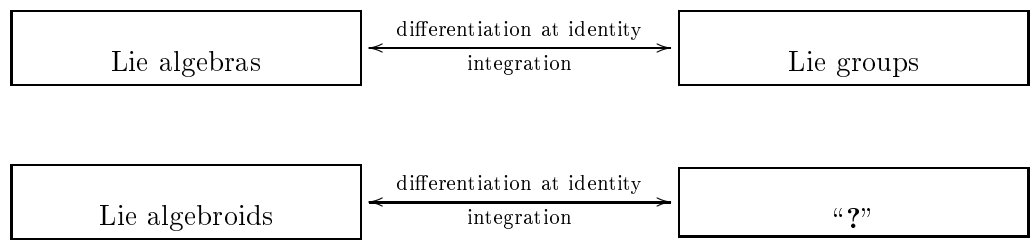
\[ \xymatrix{
& \fbox{\parbox{.3\linewidth}{\center{Lie algebras}}}
  \ar[rrrr]^{\text{differentiation at identity}} &  & & &
  \fbox{\parbox{.3\linewidth}{\center{Lie groups}}}
  \ar[llll]^{\text{integration}}   \\
 & \fbox{\parbox{.3\linewidth}{\center{Lie
   algebroids}}} \ar[rrrr]^{\text{differentiation at identity}} &  & & &
   \fbox{\parbox{.3\linewidth}{\center{``{\bf ?}''}}}
  \ar[llll]^{\text{integration}} }
\]

There have been several approaches to the object ``{\bf ?}''.
The first important approach, due to Pradines \cite{pradines}, constructs a
local Lie groupoid. Another
important approach, due to Crainic and Fernandes \cite{cf} constructs a universal topological
groupoid. Some special
cases of this construction were also observed through the Poisson
sigma model in
\cite{cafe}. In \cite{picard}, Weinstein further conjectured
that this topological groupoid must have some smooth structure.
But normal differential structures such as manifolds or orbifolds
can not serve this purpose. In this thesis, differentiable stacks (see
\cite{bx} \cite{metzler} \cite{pronk} and references therein) are
used to study this ``smooth structure''. We introduce an object which
we call Weinstein groupoid, and prove a version of
Lie's third theorem for Lie algebroids. The approach described in this paper
enriches the structure of the topological
groupoids constructed in \cite{cf}, and provides a global alternative
to the local groupoids in \cite{pradines}.

\begin{defi} [Weinstein groupoid]\label{wgpd}
A Weinstein groupoid over a manifold $M$ consists of the following
data:
\begin{enumerate}
\item an \'etale differentiable stack $\cG$ (see Definition \ref{stack});

\item (source and target) maps $\bar{\bs}$,
$\bar{\bt}$: $\cG \to M$ which are surjective submersions between
differentiable stacks;

\item (multiplication) a map $\bar{m}$: $\cG\times_{\bbs, \bbt} \cG \to
\cG$, satisfying the following properties:
\begin{itemize}
  \item $\bbt \circ \bar{m}=\bbt\circ pr_1$, $\bbs \circ \bar{m}=\bbs\circ
  pr_2$, where $pr_i: \cG \times_{\bbs, \bbt} \cG \to \cG$ is the
  $i$-th projection $\cG\times_{\bbs, \bbt} \cG \to
\cG$;
  \item associativity up
to a 2-morphism, i.e. there is a unique 2-morphism $\alpha$
between maps $\bar{m}\circ (\bar{m} \times id)$ and
$\bar{m}\circ(id\times \bar{m})$;
\end{itemize}

\item  (identity section) an injective immersion  $\bar{e}$: $M\to \cG$
such that, up to 2-morphisms, the following identities
\[
\bar{m}\circ ((\bar{e}\circ \bbt)\times id)=id, \,\,\bar{m}\circ
(id\times (\bar{e}\circ\bbs) )=id,\] hold (In particular,by
combining with the surjectivity of $\bbs$ and $\bbt$, one has
$\bbs \circ \bar{e}= id$, $\bbt \circ \bar{e}= id$ on $M$);

\item (inverse) an isomorphism of differentiable stacks
$\bar{i}$: $\cG \to \cG$ such that, up to 2-morphisms, the following
identities
\[ \bar{m}\circ (\bar{i}\times id \circ \Delta)=\bar{e}\circ\bbs, \;\;
\bar{m}\circ (id\times\bar{i}\circ \Delta)=\bar{e}\circ \bbt,\]
hold, where $\Delta$ is the diagonal map: $\cG\to \cG\times\cG$.
\end{enumerate}
Moreover, restricting to the identity section, the above
2-morphisms between maps are the $id$ 2-morphisms. Namely, for
example, the 2-morphism $\alpha$ induces the $id$ 2-morphism
between the following two maps:\[ \bar{m}\circ ((\bar{m} \circ
(\bar{e}\times\bar{e}\circ \delta))\times \bar{e} \circ
\delta)=\bar{m}\circ(\bar{e}\times(\bar{m}\circ(\bar{e}\times\bar{e}\circ\delta))\circ\delta),
\]where $\delta$ is the diagonal map: $M\to M\times M$.
\end{defi}

\noindent{\em General Remark}: the terminology involving stacks in
the above definition, as well as in the following theorems, will
be explained in detail in Chapter \ref{sect: stack}. For now,  to
get a general idea of these statements, one can take stacks simply
to be manifolds.

Our main result is the following theorem:

\begin{thm}[Lie's third theorem] \label{lieIII}
To each Weinstein groupoid one can associate a Lie algebroid. For
every Lie algebroid $A$, there are naturally two Weinstein
groupoids $\cG(A)$ and $\cH(A)$ with Lie algebroid $A$.
\end{thm}

Chapter \ref{w} is devoted to the proof of this theorem.

These two Weinstein groupoids  $\cG(A)$ and $\cH(A)$ come from the
monodromy and holonomy groupoids of some path space
respectively. Moreover we have the following conjectures.

\begin{con} $\cH(A)$ is the unique source-simply connected Weinstein
  groupoid whose Lie algebroid is $A$.
\end{con}

\begin{con} [Lie's second theorem] \label{lieII}
For any morphism of Lie algebroids $\phi: A \to B$, there is a
unique morphism $\Phi$ from the Weinstein groupoid
$\cH(A)$ to any Weinstein groupoid $\cG$ integrating $B$, such
that $d\Phi = \phi$.
\end{con}

We can apply our result to the classical integrability problem,
which studies when a Lie algebroid can be integrated into a Lie
groupoid.

\begin{thm}\label{integ}
A Lie algebroid $A$ is integrable in the classical sense iff
$\cH(A)$ is representable, i.e. it is an ordinary manifold. In
this case $\cH(A)$ is the source-simply connected Lie groupoid of
$A$ (it is also called the Weinstein groupoid of $A$ in
\cite{cf}).
\end{thm}

We can also relate our work to previous work on the integration of
Lie algebroids via the following two theorems:

\begin{thm}\label{local}
Given a Weinstein groupoid $\cG$, there is an\footnote{It is canonical up to isomorphism near the
identity section.} associated local Lie groupoid $G_{loc}$ which
has the same Lie algebroid as $\cG$.
\end{thm}

\begin{thm}\label{cf}
As topological spaces, the orbit spaces of $\cH(A)$ and $\cG(A)$
are both isomorphic to the universal topological groupoid of $A$
constructed in \cite{cf}.
\end{thm}

The three theorems above are also proved in Chapter \ref{w}

In Chapter \ref{app-ja} and Chapter \ref{fur-app}, we apply the theory
above to integration of {\em Jacobi manifolds} as
introduced by Kirillov \cite{kirillov} and Lichnerowicz
\cite{jacobi1}. Just as each Poisson manifold $P$ has an
associated Lie algebroid $T^*P$, a Jacobi manifold $M$ also has an
associated Lie algebroid $T^*M \oplus_M \R$, the direct sum of
$T^*M$ and the trivial $\R$ bundle over $M$ \cite{ks}. Integrating
the Lie algebroid of a Poisson manifold gives the {\em symplectic
groupoid} of the Poisson manifold. In 1993, Kerbrat and
Souici-Benhammadi noticed that the base manifold of a contact
groupoid is a Jacobi manifold and that the contact groupoid
integrates the Lie algebroid associated to the Jacobi manifold
\cite{ks}. In 1997, Dazord \cite{dazord} proved that locally,
every Jacobi manifold can be integrated into a contact groupoid.
Very recently, Iglesias-Ponte and Marrero \cite{igma} have
generalized contact groupoids to Jacobi groupoids and found that
the infinitesimal invariants of Jacobi groupoids are generalized
Lie bialgebroids.

We begin our study of a Jacobi manifold $M$ (integrable or not) by
constructing the Weinstein groupoid  of the associated Lie
algebroid $T^*M \oplus_M \R$. In this way, we recover the contact
groupoid of $M$ and can see clearly when  $M$ is integrable. In
Chapter \ref{app-ja}, we prove the following theorem.

\begin{thm} \label{main} Let $M$ be a Jacobi manifold, $M\times \R$ its
Poissonization (see Section \ref{homogeneous}). Then
\begin{enumerate}
\item[i)] there is an isomorphism between Weinstein groupoids
\[ \cG(T^*(M\times\R))  \cong \cG(T^*M\oplus \R)\times \R,
\quad \cH(T^*(M\times\R))\cong \cH(T^*M\oplus \R) \times \R \]
\item[ii)] $M$ is integrable as a Jacobi manifold iff $M\times
\R$ is integrable as a Poisson manifold.
\item[iii)] when $M$ is integrable, $\cH(T^*M\oplus \R)$ is the
source-simply connected contact groupoid (see Section \ref{jc}) of
$M$.
\end{enumerate}
\end{thm}

As applications of contact groupoids, we view a Poisson manifold as a
Jacobi manifold and consider its contact groupoid.  The contact
groupoid of a Poisson manifold is closely related to
the integrability of the Poisson bivector and
prequantization of its symplectic groupoid.

A Poisson bivector is a Lie algebroid 2-cocycle of the Lie
algebroid $T^*M$. It is called integrable iff it comes from a Lie
groupoid 2-cocycle. The relation between Lie algebroid cocycles
(cohomologies) and Lie groupoid cocycles (cohomologies) is
explained in \cite{m-von-est} \cite{xw1}. We have the following
result on the integrability of Poisson bivectors:

\begin{thm}\label{main2}
The Poisson bivector $\Lambda$ of a Poisson manifold $M$ can be
integrated into a Lie groupoid 2-cocycle
if and only if $M$ is integrable as a Poisson manifold and the
symplectic form of the source-simply connected symplectic groupoid
is exact.
\end{thm}
Equivalent conditions more convenient for computations will be given in Chapter
\ref{fur-app}.

Prequantizations of symplectic
groupoids were introduced by Weinstein and Xu in \cite{xw1}, as the first step in quantizing symplectic groupoids for the
purpose of quantizing Poisson manifolds. Using
the contact groupoids constructed above in Theorem \ref{main}, we
are able to construct the prequantizations of symplectic
groupoids. In Chapter \ref{fur-app}, we have the following result:

\begin{thm}\label{main3}
If $(\Gamma_s(M), \Omega)$ is a symplectic groupoid with
$\Omega\in H^2(\Gamma_s(M), \Z)$, then $M$ can be integrated into
a contact groupoid $(\Gamma_c(M), \theta, 1)$. Furthermore, if we
quotient out by a $\Z$ action, $\Gamma_c(M)/\Z$ is a
prequantization of $\Gamma_s(M)$ with connection 1-form
$\bar{\theta}$ induced by $\theta$. Moreover, $(\Gamma_c(M)/\Z,
\bar{\theta}, 1)$ is also a contact groupoid of $M$.
\end{thm}

\chapter{Differentiable stacks}\label{sect: stack}
The
notion of stack has  been extensively studied in
algebraic geometry for the past few decades (see for example
\cite{beffgk} \cite{dm} \cite{lmb} \cite{v1}). However
stacks can also be defined over other categories, such as the
category of topological spaces and category of smooth manifolds
(see for example \cite{SGA4} \cite{bx} \cite{metzler} \cite{pronk}
\cite{v2}). In this section we collect certain facts about stacks in the
differentiable category that we will use later. Many of them
already appeared in the literature \cite{bx} \cite{metzler}.

\section{Stacks over the category of differentiable (or Banach)
manifolds}
\subsection{The definitions}
In general, a stack over a category is a category
fibred in groupoids satisfying some sheaf-like conditions
\cite{artin} \cite{dm}. In particular, here, we suppose that the
base category $\cC$ is either the category of differentiable manifolds
or the category of
Banach manifolds. Banach manifolds are
possibly infinite dimensional and have Banach spaces as their
local charts \cite{lang}. We endow $\cC$ with a Grothendieck
topology \cite{SGA4} by declaring $\{ f_i: U_i \to S\} $ to be a
{\em covering family} if each $f_i$ is an open embedding and
$\cup_i f_i(U_i)=S$. It is easy to check that this forms a basis
for a Grothendieck topology on $\cC$.  See the above citation or
\cite{metzler} Section 2 for the detailed definition of a
Grothendieck topology.

\begin{defi} [categories fibred in groupoids]
A  category fibred in groupoids $\cX \to \cC$ is a category $\cX$
over a base category $\cC$, together with a functor $\pi: \cX \to
\cC$, such that the following two axioms are satisfied:
\begin{enumerate}
\item[i)](pullback) for every morphism $V\to U$ in $\cC$, and every object
$x$ in $\cX$ over $U$ (i.e. $\pi(x)=U$), there exists an object
$y$ over $V$ and a morphism $y\to x$ lifting $V\to U$;
\item[ii)] for every composition of morphisms $W\to V\to U$ in
$\cC$ and morphisms $z\to x$ lying over $W\to U$ and $y\to x$
lying over $V\to U$, there exists a unique morphism $z\to y$ such
that the triangle of morphisms between $x, y, z$ commute.
\end{enumerate}
\end{defi}
\begin{remark}
Here the object $y$ over $V$ exists (if $x$ exists) and is unique
up to a unique morphism by ii). We call $y$ a
pullback\footnote{``A pullback'' is used here since $y$ is not
really unique, but nevertheless, by abuse of notation, we will
still denote it by $f^*x$ or $x|V$, where $f $ is the morphism
$V\to U$.}  of $x$ through $f$. Let $\cX_U$ be the category whose
objects are all the objects in $\cX$ lying over $U$ and whose
morphisms are all
the morphisms lying over $id_U$. By i), $\cX_U$ is not empty if $\cX$
is not empty. By
ii), using $U\overset{id}{\to} U \overset{id}{\to} U$, any
morphism between two objects $x$ and $x'$ in $\cX_U$ is
invertible, i.e. it is an isomorphism. Therefore such a ``fibre''
$\cX_U$ of $\cX$ over $\cC$ is a groupoid (set-theoretically).
\end{remark}

\begin{defi}[stacks]
We call $\cX$ a  stack over $\cC$ if,
\begin{enumerate}
\item[i)] $\cX\to \cC$ is a category fibred in groupoids;
\item[ii)] for any $S$ in $\cC$ and any two objects $x, y$ in
$\cX_S$, the contravariant functor $\mathrm{Isom} (x, y)$, defined
by
\[ \begin{split}
\mathrm{Isom}(x, y)(U)= &\{ (f,\phi)| f:U\to S \; \text{is a
morphism in} \; \cC, \\ & \phi: f^* x \to f^* y \; \text{is a
morphism in} \; \cX_U \}
\end{split}
\]
is a sheaf, where $f^*x$ can be any possible pullback of $x$ via
$f$.
\item[iii)] for any $S$ in $\cC$, and every covering family
$\{U_i\}$ of $S$, every family $\{x_i\}$ of objects $x_i \in
\cX_{U_i}$ and every family of morphisms $\{ \phi_{ij}\}$,
$\phi_{ij}$: $x_j|U_{ij} \to x_i|U_{ij}$, satisfying the cocycle
condition $\phi_{kj} \circ \phi_{ji} =\phi_{ki}$ (which holds in
the fibre $\cX_{U_{ijk}}$), there exists a global object $x$
over $S$, together with isomorphisms $\phi_i: x|U_i \to x_i$ such
that we have $\phi_{ij}\circ \phi_j =\phi_i $ over $U_{ij}$.
\end{enumerate}
\end{defi}
\begin{remark}
Roughly, i) says that pullbacks exist and are unique up to a
unique morphism; iii) says that the elements satisfying the gluing
conditions can glue together; ii) tells us that the element glued is
unique up to a unique isomorphism.
\end{remark}

\subsection{Representability}
\begin{ep}\label{mfd}
Given a (Banach) manifold  $M$, one can view it as a stack over
$\cC$. Let $\underline{M}$ be the category where
\[Obj(\underline{M})=\{ (S, u): S\in \cC, u\in Hom(S, M)\}, \]and
a morphism $(S, u) \to (T, v)$ of objects is a morphism $f: S\to
T$ such that $u=v\circ f$. This category encodes all the
information of $M$ and no more than this in the sense that the
morphisms between stacks $\underline{M}$ and $\underline{M}'$ all
come from the ordinary morphisms between $M$ and $M'$, i.e. $\cC$
is a full subcategory of the category of stacks.
In this way, the notion of stacks generalizes the notion of
manifolds. A stack isomorphic to $\underline{M}$ for some
$M\in\cC$ is called {\bf representable}.
\end{ep}

\section{Differentiable (Banach) stacks}
From now on, by ``manifolds'' we mean finite dimensional smooth
manifolds unless we put in front the word ``Banach''. However, all
the theory can be easily extended to Banach manifolds (also see Remark \ref{banach}). {\bf Morphisms} between
stacks are functors between categories over $\cC$ viewing the
stacks as categories over $\cC$, and {\bf 2-morphisms} between two
stack morphisms are natural transformations between functors.




\begin{defi}[monomorphisms and epimorphisms]
A morphism of stacks $f: \cX\to\cY$ is called a monomorphism if
for any two objects $x, x'$ in $\cX$ over $S\in \cC$ and any arrow
$\eta: f(x)\to f(x')$ there is a unique arrow $x\to x'$ as the
preimage of $\eta$ under $f$.

A morphism of stacks $f: \cX\to\cY$ is called an epimorphism if
for any objects $y$ in $\cY$ over $S\in \cC$ there is a covering
$S_i$ of $S$ and objects $x_i$ in $\cX$ over $S_i$ such that
$f(x_i)\cong y|_{S_i}$, for all $i$.
\end{defi}

\begin{defi}[fibre product \cite{bx}] \label{def-fp} Given two morphisms of
stacks $\phi:\cX\to \cZ$ and $\varphi\cY \to \cZ$, one can form
the fibre product $\cX\times_{\phi, \cZ, \varphi} \cY$ in the
following way: the objects over $S\in \cC$ are $(x, \eta, y)$
where $x\in \cX_S$ and $y\in \cY_S$ and $\eta$ is an arrow from
$\phi(x)$ to $\varphi(y)$ in $\cZ_S$; the morphisms over $S\to S'$
are arrows from $(x, \eta, y)$ to $(x', \eta', y')$ consist of
compatible morphisms $x\to x'$ and $y\to y'$ and $\eta \to \eta'$.
\end{defi}

\begin{defi}[representable surjective submersions \cite{bx}]
A morphism $f:\cX\to\cY$ of stacks is a  representable
submersion if for every (Banach) manifold $M$ and every morphism
$M\to\cY$ the fibred product $\cX\times_{\cY}M$ is representable
and the induced morphism $\cX\times_{\cY}M\to M$ is a submersion.
$f$ is a  representable surjective submersion if it is also
an epimorphism.
\end{defi}

\begin{defi}[differentiable (Banach) stacks \cite{bx}] \label{stack}
A  differentiable (Banach) stack $\cX$ is a stack over the
category $\cC$ of differentiable (Banach) manifolds with a
representable surjective submersion $\pi: X\to\cX$ from a (Banach)
Hausdorff manifold $X$. $X$ together with the structure morphism
$\pi: X\to\cX$ is called an atlas for $\cX$.
\end{defi}

\begin{ep}
A Hausdorff (Banach) manifold is a differentiable (Banach) stack
by definition.
\end{ep}

\begin{ep} \label{bg}
Let $G$ be a Lie group. The set of  principal $G$-bundles forms a
stack $BG$ in the following way. The objects of $BG$ are
\[ Obj (BG)=\{\pi:  P\to M | P \;\text{ is a principal $G$-bundle over
  $M$}. \} \]
A morphism between two objects $(P, M)$ and $(P', M')$ is a
morphism $M\to M'$ and a $G$-equivariant morphism $P\to P'$
covering $M\to M'$. Moreover $BG$ is a differentiable stack. The
atlas is simply a point $pt$. The map
\[\begin{split}
\pi: &(f: M \to pt) \mapsto (M\times_{f,pt, pr} G),
\\ &(a: (f_1: M_1\to pt)\to (f_2: M_2 \to pt)) \mapsto ((x_1, g_1)\mapsto (a(x_2),
g_2)),
\end{split}
\] 
(where $pr$ is the projection from $G$ to the point $pt$) is a
representable surjective submersion.
\end{ep}

\section{Morphisms and 2-morphisms}

We have the following
two easy properties of representable surjective submersions:

\begin{lemma} [composition]\label{comp}
The composition of two representable (surjective) submersions is
still a representable (surjective) submersion.
\end{lemma}
\begin{proof}
For any manifold $U$ with map $U\to\cZ$, consider the following
diagram
\[
\begin{CD}
\cX\times_\cY \cY \times_\cZ U
@>\tilde{f}>>\cY \times_\cZ U @>\tilde{g}>> U \\
@VVV @VVV @VVV \\
\cX @> f>> \cY @>g>> \cZ .
\end{CD}
\]
Since $f$ and $g$ are representable submersions, $\cY \times_\cZ
U$ is a manifold so that  $\cX\times_\cZ U=\cX\times_\cY\cY
\times_\cZ U$ is  also a manifold. Since $\tilde{g}$ and
$\tilde{f}$ are submersions, $\tilde{g}\circ \tilde{f}$ is also a
submersion. The composition of two epimorphisms is still an
epimorphism.
\end{proof}

\begin{lemma}[base change] \label{BC}
In the following diagram
\begin{equation}\label{bc}
\begin{CD}
\cX \times_{\cY}\cZ @>g>> \cZ \\
@VVV @VVV \\
\cX @>f>>\cY,
\end{CD}
\end{equation}
where $\cX$ and $\cY$ are differentiable stacks (but not
necessarily $\cZ$), if $f$ is a representable (surjective)
submersion, then so is $g$.
\end{lemma}
\begin{proof}
For any manifold $U$ mapping to $\cZ$, we have the following
diagram
\[
\begin{CD}
\cX \times_{\cY}\cZ\times_{\cZ} U @>h>> U \\
@VVV @VVV \\
\cX \times_{\cY}\cZ @>g>> \cZ
\end{CD}
\]
Composing the above diagram with \eqref{bc}, one can see that
$\cX{\times}_{\cY}\cZ\times_{\cZ} U=\cX \times_{\cY} U$ is a
manifold and $h$ is a submersion because $f$ is a representable
submersion. Therefore $g$ is a representable submersion. Moreover,
the base change of an epimorphism is clearly still an epimorphism.
\end{proof}
\begin{remark} In general, we call the procedure of obtaining $g$
from $f$ base change of $\cX \to \cY$ by $\cZ \to \cY$ and we call
the result map $g$ the base change of  $\cX \to \cY$ by $\cZ \to
\cY$.
\end{remark}

\begin{defi}[smooth morphisms of differentiable stacks]
A morphism $f:\cX\to\cY$ of differentiable stacks is smooth if for
any atlas $g:X\to\cX$ the composition $X\to\cX\to\cY$ satisfies
the following: for any atlas $Y\to\cY$ the induced morphism
$X\times_{\cY}Y\to Y$ is a smooth morphism of manifolds.
\end{defi}

In the rest of the article, morphisms between differentiable
stacks are referred to as smooth morphisms without special
explanations.

\begin{defi}[embeddings]
A morphism $f:\cX\to\cY$ of stacks is an embedding if for every
submersion $M\to\cY$ from a manifold $M$ the product
$\cX\times_{\cY}M$ is a manifold and the induced morphism
$\cX\times_{\cY}M\to M$ is an embedding of manifolds.
\end{defi}

\begin{defi}[immersions, \'etale maps and closed immersions \cite{metzler}]
A morphism $f:\cX\to\cY$ of stacks is an immersion (resp an
\'etale map, or closed immersion) if for every representable
submersion $M\to\cY$ from a manifold $M$ the product
$\cX\times_{\cY}M$ is a manifold and the induced morphism
$\cX\times_{\cY}M\to M$ is an immersion (resp. an \'etale map, or
closed immersion) of manifolds.
\end{defi}

\begin{defi}[\'etale differentiable stacks]
A differentiable stack $\cX$ is called  \'etale if there is a
presentation $\pi: X\to\cX$ with $\pi$ being \'etale.
\end{defi}

\begin{lemma}\label{sm}
A morphism $\cX\to\cY$ is smooth if and only if there exist an
atlas $X\to\cX$ of $\cX$ and an atlas $Y\to\cY$ of $\cY$ such that
the induced morphism $X\times_{\cY}Y\to Y$ is a smooth morphism of
manifolds.
\end{lemma}

\begin{proof}
One implication is obvious. Suppose that $X\times_{\cY}Y\to Y$ is
smooth. Let $T\to\cX$ be another atlas. Then using base change of
$X\times_\cY Y\to \cX$,  $T\times_{\cX}X\times_{\cY}Y$ is a
manifold and $T\times_{\cX}X\times_{\cY}Y \to X\times_\cY Y$ is a
submersion, hence a smooth map. The map
$T\times_{\cX}X\times_{\cY}Y\to Y$ factors as
$T\times_{\cX}X\times_{\cY}Y\to X\times_{\cY}Y\to Y$. Hence
$T\times_{\cX}X\times_{\cY}Y\to Y$ is smooth. It also factors as
$T\times_{\cX}X\times_{\cY}Y\to T\times_{\cY}Y\to Y$. Similarly,
the map $T\times_{\cX}X\times_{\cY}Y\to T\times_{\cY}Y$ is a
submersion, hence $T\times_{\cY}Y\to Y$ is smooth.

Now assume that $U\to\cY$ is an atlas of $\cY$. The induced map
$T\times_{\cY}Y\times_{\cY}U\to Y\times_{\cY}U$ is smooth because
it is the base-change of a smooth map $T\times_{\cY}Y\to Y$ by a
submersion $Y\times_{\cY}U\to Y$. One can find a collection of
locally closed submanifolds in $Y\times_{\cY} U$ which form an
open covering family for $U$. Since being smooth is a local
property, it follows that $T\times_{\cY}U\to U$ is smooth as well.
\end{proof}

\begin{lemma}\label{immersion}
A morphism from a manifold $X$ to a differentiable stack $\cY$ is
an immersion if and only if $X\times_{\cY} U\to U$ is an immersion
for some atlas $U\to\cY$.
\end{lemma}
\begin{proof}
One implication is obvious. If $X\times_{\cY} U\to U$ is an
immersion, let $T\to\cY$ be any submersion from a manifold $T$.
The map $X\times_{\cY}U\to U$ is transformed by base-change by a
submersion  $U\times_{\cY}T\to U$ to a map
$X\times_{\cY}U\times_{\cY}T\to U\times_{\cY}T$, which is an
immersion since being an immersion is preserved by base-change.
One can find a collection of locally closed submanifolds $\{T_i\}$
in $U\times_{\cY} T$ which forms a family of charts of $T$.
Moreover $X \times_\cY T$ is a manifold because $T$ is an atlas of
$\cY$. Using base changes, one can see that $X\times T_i \to T_i$
is an immersion and that $\{ X\times_\cY T_i\}$ forms an open
covering family of $X\times_\cY T$. Since being an immersion is a
local property, it follows that $X\times_{\cY}T\to T$ is an
immersion as well.
\end{proof}

Similarly we have
\begin{lemma}\label{closedimmersion}
A morphism $\cX\to\cY$ of differentiable stacks is a closed
immersion if and only if $\cX\times_{\cY} U\to U$ is a closed
immersion for some atlas $U\to\cY$.
\end{lemma}

\begin{defi}[submersions]
A morphism $f:\cX\to\cY$ of differentiable stacks is called a
submersion\footnote{This is different from the definition in
\cite{metzler}, where $N\times_\cY \cX$ is required to be a
manifold. } if for any atlas $M\to\cX$, the composition
$M\to\cX\to\cY$ satisfies the following: for any atlas $N\to\cY$
the induced morphism $M\times_{\cY}N\to N$ is a submersion.
\end{defi}

\begin{remark}
In particular, a representable submersion is a submersion. But the
converse is not true: for example the source and target maps $\bbs$
and $\bbt$ that we will define in Section \ref{wgpd} are
submersions but not representable submersions in general.
\end{remark}

We will use later the following result about the fibred products
using submersions:

\begin{prop}[fibred products]\label{fp}
Let $Z$ be a manifold and $f:\cX\to Z$ and $g:\cY\to Z$ be two morphisms of
differentiable stacks. If either $f$ or $g$ is a
submersion, then $\cX \times_Z\cY$ is a differentiable stack.
\end{prop}
\begin{proof}
Assume that $f:\cX \to Z$ is a submersion. By definition, for any
atlas  $X\to\cX$, the composition $X\to\cX\to Z$ is a submersion.
Let $Y$ be a presentation of $\cY$, then $X\times_Z Y$ is a
manifold. To see that $\cX \times_Z \cY$ is a differentiable
stack, it suffices to show that there exists a representable
surjective submersion from $X\times_Z Y$ to $\cX\times_Z \cY$. By
Lemma \ref{BC}, $X\times_Z Y \to \cX\times_Z Y$ and $\cX\times_Z Y
\to\cX\times_Z\cY$ are representable surjective submersions. By
Lemma \ref{comp}, their composition is also a representable
surjective submersion.
\end{proof}

\begin{lemma}\label{pa}
Let $\cX, \cY$ be stacks with maps $\cX\to Z$ and $\cY\to Z$ to a
manifold $Z$, one of which a submersion, and let $X\to\cX,
Y\to\cY$ be atlases of $\cX$ and $\cY$ respectively. Then
$X\times_Z Y\to\cX\times_Z\cY$ is an atlas of $\cX\times_Z\cY$.
\end{lemma}
\begin{proof}
Note that $X\times_Z Y$ is a manifold because one of $\cX\to Z$
and $\cY\to Z$ is a submersion. $X\times_Z Y\to\cX\times_Z\cY$
factors into $X\times_Z Y\to\cX\times_Z Y\to\cX\times_Z\cY$.
$X\times_Z Y\to\cX\times_Z Y$ is a representable surjective
submersion because $X\to\cX$ is. $\cX\times_Z Y\to\cX\times_Z\cY$
is a representable surjective submersion because $Y\to\cY$ is.
Thus $X\times_Z Y\to\cX\times_Z\cY$ is a representable surjective
submersion.
\end{proof}

\section{Lie groupoids and differentiable stacks}
Next we explain the relationship between stacks and groupoids.
\subsection{From stacks to groupoids}
Let $\cX$ be a differentiable stack. Given an atlas $X_0\to\cX$,
we can form
\[X_1:=( X\times_{\cX}X)\rightrightarrows X\] with the two maps
being projections from the first and second factors onto $X$. By the
definition of an atlas,
$X_1$ is a manifold. Moreover it has a natural groupoid structure with source and target maps the
two maps above. We call this groupoid a presentation of $\cX$.
Different atlases give rise to different presentations (see for
example the appendix to \cite{v1}). An \'etale differential stack
will have a presentation by an \'etale groupoid.

\begin{ep}
In Example \ref{mfd} we have the stack $\underline{M}$ with the atlas
$M\to\underline{M}$. $M\times_{\underline{M}}M$ is just the
diagonal in $M\times M$, thus is isomorphic to $M$. Hence we have
a groupoid $M\rightrightarrows M$ with two maps both equal to the
identity. This is clearly isomorphic (as a groupoid) to the
transformation groupoid $\{id\}\times M\rightrightarrows M$, where
$\{id\}$ represents the group with only one element.
\end{ep}

\begin{ep}
In Example \ref{bg}, an atlas of the stack $BG$ is a point $pt$.
The fibre-product $pt\times_{BG} pt$ is $G$. So a groupoid
presenting $BG$ is simply $G\rightrightarrows pt$.
\end{ep}

\begin{ep}\label{pa2}
In the situation of Lemma \ref{pa}, put $X_1=X\times_{\cX}X$ and
$Y_1=Y\times_{\cY} Y$, then $\cX\times_Z\cY$ is presented by the
groupoid $(X_1\times_Z Y_1\rightrightarrows X\times_Z Y)$. This
follows from the fact that $(X\times_Z
Y)\times_{\cX\times_Z\cY}(X\times_Z Y)\cong
(X\times_{\cX}X)\times_Z (Y\times_{\cY} Y)$.
\end{ep}
\subsection{From groupoids to stacks}
Conversely, given a
groupoid \gpoid, one can associate a quotient stack $\cX$ with an
atlas $G_0\to\cX$ such that $G_1=G_0\times_{\cX} G_0$. Here we recall the
construction given in \cite{bx} for differentiable stacks. We
begin with several definitions.

\begin{defi}[groupoid action]
A Lie groupoid \gpoid right (resp. left) action on a manifold $M$  consists of the following data: a moment map $J: M\to G_0$
and a smooth map $\Phi: M \times_{J,\bt} G_1 (\text{resp.}
G_1\times_{\bs, J} M) \to M$ such that
\begin{enumerate}
\item $J(\Phi(m, g))=\bs(g)$  (resp. $J(\Phi(g, m)) =\bt(g)$);
\item $\Phi(\Phi(m, g), h)=\Phi(m, gh)$ (resp. $\Phi(h, \Phi(g, m))=\Phi(hg, m)$);
\item $\Phi(m, J(m))=m$ (or $\Phi(J(m), m)=m$).
\end{enumerate}
Here we identify $G_0$ as the identity section of $G_1$. The action $\Phi$ is also denoted by ``$\cdot$'' for simplicity.
\end{defi}

\begin{defi}[Lie groupoid principal bundles, or torsors]
A manifold $P$ is a {\bf right (resp. left) principal bundle} of a
Lie groupoid $H$ over a manifold $S$, if
\begin{enumerate}
\item there is a surjective submersion $\pi: P \to S$;
\item $H$  acts from the right (resp. left) on $P$ fibrewise with respect to $\pi$, that is, $\pi(p\cdot h)=\pi(p)$ for all $(p, h)\in P\times_{J, \bt} H_1$ (resp. $\pi(h\cdot p)=\pi(p)$ for all $(h, p)\in H_1\times_{\bs, J} P$);
\item  the $H$ action is  free
and transitive on each fiber of $\pi$, that is the map $$(pr_1, \Phi): P\times_{J, \bt} H_1 \to P \times{\pi, S, \pi}P, \quad (p, h)\mapsto (p, ph)$$ is a diffeomorphism (resp. $(\Phi, pr_2): H_1\times_{\bs, H_0, J} P \to  P\times{\pi, S, \pi}P, \quad (p, h)\mapsto (hp, p)$ is a diffeomorphism);

A right principal $H$ bundle is called a $H$-torsor too.
\end{enumerate}

\end{defi}

\begin{remark}
Since the action is  free and transitive, one can see
that $\pi$ descends to a diffeomorphism $\bar{\pi}: P/H\cong S$. Thus an $H$-orbit is an embedded submanifold $\bar{\pi}^{-1}(x)$. It is not hard to see that given a free $H$-action on $P$, if the $H$-orbits are embedded submanifolds, then the $H$-action is proper. Thus we obtain that the $H$ action is free and proper for an $H$-principal
bundle. Thus when $H$ is a group, this gives us the usual notion of $H$-principal bundle. On the other hand,  by the groupoid-version of   the groupoid-version of slice theorem (see \cite[Lemma 3.11]{zz}), if $H$ action is free and proper,
then the quotient $P/H$ inherits a manifold structure. A more precise statement is that the quotient stack $[P/H]$ is representible if and only if the $H$-action is free and proper.
\end{remark}

Let $G=$(\gpoid) be a Lie groupoid. Denote by\footnote{Or
alternatively by $[G_0/G_1]$.} $BG$ the category of right
$G$-principal bundles. We now show that $BG$ is moreover a
differentiable stack. An object $Q$ of $BG$ over $S\in \cC$ is a
right $G$-principal bundle over $S$. A morphism between two $G$
torsors $\pi_1:Q_1\to S_1$ and $\pi_2:Q_2\to S_2$ is a smooth map
$\Psi$ lifting the morphism $\psi$ between the base manifolds
$S_1$ and $S_2$ (i.e. $\psi\circ\pi_1=\Psi\circ\pi_2$) such that
$\Psi$ is $G_1$-equivariant, i.e. $\Psi(q_1 \cdot g) =\Psi(q_1)
\cdot g$ for $(q_1, g) \in M \times_{J_1,\bt} G_1$, where  $\pi_i$
are projections of torsors $Q_i$ onto their bases and $J_i$ are
the moment maps, $i=1, 2$.

{\bf Note}: the above condition implies that $J_2\circ\Psi = J_1$.

This makes $BG$ a category over $\cC$. According to \cite{bx} it
is a differentiable stack presented by the Lie groupoid \gpoid: an
atlas $\phi: G_0\to BG$ can be constructed as follows: for $f:
S\to G_0$, we assign the manifold $Q=S \times_{f,\bt} G_1$.  ($Q$
is a manifold because $\bt$ is a submersion). The projection $\pi:
Q\to S$ is given by the first projection and the moment map $J: Q
\to G_0$ is the second projection composed with $\bs$. The
groupoid action is defined by
\[ (s, g) \cdot h= (s, gh), \text{ for all possible choices of } \; (s, g) \in S
\times_{f,\bt} G_1, h \in G_1. \] The $\pi$-fiber is simply a copy
of the $\bt$-fiber, therefore the action of $G_1$ is free and
transitive. $\phi$ is a representable surjective
submersion and $G_1= G_0 \times_{\phi,\phi} G_0$ fits in the
following diagram:
\[
\begin{CD}
  G_1 @>\bt>> G_0 \\
  @V \bs VV @V \phi VV \\
  G_0@>\phi>> BG.
\end{CD}
\]
We refer to \cite{bx} for more details.\\
\begin{ep}
In the case of the trivial transformation groupoid $\{id\}\times
M\rightrightarrows M$ it is easy to see that the stack constructed
above is $\underline{M}$.
\end{ep}
\subsection{Morita equivalence}
To further explore the correspondence between stacks and
groupoids, we need the following definition.

\begin{defi}[Morita equivalence \cite{mrw}]
Two Lie groupoids $G=$(\gpoid) and $H=$(\hpoid) are {\bf Morita
equivalent} if there exists a manifold $E$, such that
\begin{enumerate}
\item $G$ and $H$ act on $E$ from the left and right respectively
 with moment maps $J_G$ and $J_H$ and the two actions commute;
\item The moment maps are surjective submersions;
\item The groupoid actions on the fibre of the moment maps are
free and transitive.
\end{enumerate}
 Such an $E$ is called a {\bf Morita bibundle}
of $G$ and $H$.
\end{defi}

\begin{remark}
If the Morita bibundle is given by an honest groupoid morphism
$\phi: G\to H$, then $\phi$ is called a {\em strong equivalence}
\cite{moerdijk} from $G$ to $H$, and we say that $G$ is strongly
equivalent to $H$.  For any two Morita equivalent groupoids $G$
and $H$, there exists a third groupoid $K$ which is strongly
equivalent to both of them \cite{moerdijk}.
\end{remark}

\begin{prop}(\cite{bx} \cite{metzler} \cite{pronk}) Different presentations of a stack arising from different
atlases are Morita equivalent. Two Lie groupoids present isomorphic differential stack if and only if they are Morita
equivalent.
\end{prop}

\section{Morphisms and 2-morphisms---in the world of groupoids}
\subsection{Morphisms}
(1-)morphisms between stacks can be realized on the level of
groupoids.

\begin{defi}[HS morphisms \cite{hs}]
A  {\bf Hilsum-Skandalis (HS) morphism} of Lie groupoids from $G$
to $H$ is a triple $(E, J_G, J_H)$ such that:
\begin{enumerate}
\item The bundle $J_G: E\to G_0$ is a right $H$-principal bundle with moment map $J_H$;
\item $G$ acts on $E$ from left with moment map $J_G$;
\item The actions of $G$ and $H$ commute, i.e. $(g\cdot x) \cdot
  h=g\cdot (x\cdot h)$.
\end{enumerate}
We call $E$ an {\bf HS bibundle}.
\end{defi}

\begin{remark}
\;\; \item[i)] In the above definition, (3) implies that $J_H$ is
$G$
  invariant and $J_G$ is $H$ invariant.

\;\; \item[ii)] For a homomorphism of Lie groupoids
$f:$$(G_1\underset{\mathbf{t}}{\overset{\mathbf{s}}{\rightrightarrows}}G_0)$$\to$$(H_1\underset{\mathbf{t}}{\overset{\mathbf{s}}{\rightrightarrows}}H_0)$,
one can form an HS morphism via the bibundle $G_0 \times_{f,H_0,
\bt} H_1$ \cite{hm}. Thus the notion of HS morphisms generalizes
the notion of Lie groupoid morphisms.

\;\; \item[iii)]  The identity HS morphism of
$G_1\underset{\mathbf{t}}{\overset{\mathbf{s}}{\rightrightarrows}}G_0$
is given by
 $G_0 \times_{\bt} G_1 \times_{\bs} G_0$. An HS
  morphism is invertible if the bibundle is not only right
  principal but also left principal. In other words, it is a Morita equivalence.

\;\; \item[iv)] Two HS morphisms $E$: $(G_1\rightrightarrows
G_0)$$\to$$(H_1\rightrightarrows H_0)$ and $F$:
$(H_1\rightrightarrows H_0)$$\to$$(K_1\rightrightarrows K_0)$ can
be composed to obtain an HS morphism $(G_1\rightrightarrows
G_0)$$\to$$(K_1\rightrightarrows K_0)$ with the bibundle
$E\times_{H_0} F/H$, where $H$ acts on $E\times_{H_0} F$ by $(x,
y)\cdot h=(x h, h^{-1} y)$ ($G$ and $K$ still have left-over
actions on it). The composition is not strictly associative, but
it is associative up to 2-morphism which we will introduce later.
For this subtlety, please see \cite{hm}.
\end{remark}

\begin{prop}[HS and smooth morphism of stacks]\label{hs}
HS morphisms of Lie groupoids correspond to smooth morphisms of
differentiable stacks. More precisely, an HS morphism $E$:
(\gpoid)$\to$(\hpoid) induces a smooth morphism of differentiable
stacks $\phi_E: BG_1 \to BH_1$. On the other hand, given a smooth
morphism $\phi$: $\cX \to \cY$ and atlases $G_0\to\cX, H_0\to\cY$,
$\phi$ induces an HS morphism $E_\phi$: (\gpoid)$\to$(\hpoid),
where $(G_0\times_{\cX}G_0)=$\gpoid and
$(H_0\times_{\cY}H_0)=$\hpoid present $\cX$ and $\cY$
respectively.
\end{prop}
\begin{proof}
Suppose that  $(E, J_G, J_H)$ is an HS morphism. Given a right
$G$-principal bundle  $P$ over $S$ with the moment map $J_P$, we
form
\[Q=P\times_{G_0}E/G,\] where the $G$-action is given by
\[(p, x)\cdot g=(pg, g^{-1}x), \;\text{if}\;
J_P(p)=\bt(g)=J_G(x).\] Since the action of $G$ is free and proper
on $P$, the $G$-action on $P\times_{G_0}E$ is also free and
proper. So $Q$ is a manifold. In the following steps, we will show
that $Q$ is a $H$-torsor, then we can define $\phi_E$ by
$\phi_E(P)=Q$ on the level of objects.

\begin{enumerate}
\item Define $\pi_Q: Q\to S$ by $\pi_Q([(p, x)])=\pi_P (p)$. Since
$\pi_P: P\to S$ is $G$ invariant, $\pi_Q$ is a well-defined smooth
  map. Since any small enough curve $\gamma(t)$ in $S$ can be pulled back by
  $\pi_P$ as $\tilde{\gamma}(t)$ in $P$, $\gamma(t)$ can be pulled back by
  $\pi_Q$ to $Q=P\times_{G_0}E$ as $[(\tilde{\gamma}(t),
    x)]$. Therefore $\pi_Q$ is a surjective submersion.

\item Define $J_Q: Q\to H_0$ by $J_Q([p, x])=J_H(x)$. Since $J_H$
is $G$ invariant, $J_Q$ is well-defined and smooth.

\item Define $H$ action on $Q$ by $[(p, x)]\cdot h=[(p, xh)]$. It
is well defined since the actions of $G$ and $H$  commute. If
$[(p, x)]\cdot h=[(p, x)]$, then there exists a $g\in G_1$, such
that $(pg, g^{-1}xh)=(p, x)$. Since the $G$ action is free on $P$
and the $H$ action is free on $E$, we must have $g=1$ and $h=1$.
Therefore the $H$ action on $Q$ is free.

\item If $[(p, x)]$ and $[(p', x')]$ belong to the same fibre of
$\pi_Q$, i.e. $\pi_P(p)=\pi_P(p')$, then there exists a $g\in
G_1$, such that $p'=pg$. So $[(p, x)]=[(p', g^{-1} x)]$. Since
$J_G(x')=J_P(p')=\bs(g)=J_G(g^{-1}x)$, there exists an $h \in
H_1$, such that $x'h=g^{-1}x$. So $[(p', x')] h=[(p, x)]$, i.e.
the $H$ action on $Q$ is transitive.
\end{enumerate}

On the level of morphisms, we define a map which takes a morphism
of right $G$ principal bundles $f: P_1\to P_2$ to a morphism of
right $H$ principal bundles
\[ \tilde{f}: P_1\times E/G \to P_2 \times E /G, \; \text{given by} \; [(p,
  x)]\mapsto [(f(p), x)]. \]
Therefore $\phi_E$ is a map between stacks. The smoothness of
$\phi_E$ follows
from the following claim and Lemma \ref{sm}. \\
{\em Claim}: As a manifold, $E$ is isomorphic to $H_0 \times_\cY
G_0$, where the map $G_0 \to \cY$ is the composition of the atlas
projection $\pi_G: G_0\to \cX$ and $\phi_E$. Under this
isomorphism, the two moment maps $J_H$ and $J_G$ coincide with the
projections from $H_0\times_\cY G_0$ to
$H_0$ and $G_0$ respectively. \\
{\em Proof of the Claim:} Since  the category of manifolds is a
full subcategory of the category of stacks, it suffices to show
$E$ and $G_0\times_\cY G_0$ are isomorphic as stacks.

Examining the definition of fibre product of stacks (Definition
\ref{def-fp}, we see that an object in $H_0\times_\cY G_0$ over a
manifold $S\in \cC$ is $(f_H, f, f_G)$ where $f_H: S\to H_0$,
$f_G: S\to G_0$ and $f$ is an $H$ equivariant map fitting inside
the following diagram:
\[
\begin{CD}
S\times_{f_H, H_0, \bt} H_1 @>f>> S\times_{f_G, G_0, J_G} E \\
@VVV  @VVV \\
S @>id>> S.
\end{CD}
\]
Here we use $(x, e)\mapsto [(x, 1_x, e)]$ to identify
$S\times_{f_G, G_0, J_G} E$ with $(S\times_{f_G,
\bt}G_1\times_{\bs, J_G}E)/G$ which is the image of the trivial
torsor $S\times_{f_G, \bt}G_1$ under the map $\phi_E\circ \pi_G$.
Then by $x\mapsto pr_E\circ f(x, 1_x)$, $f$ gives a map $\psi_f:
S\to E$, which is an object of the stack $E$.

On the other hand for any $\psi: S\to E$, one can construct a map
$f:S\times_{f_H, H_0, \bt} H_1 \to S\times_{f_G, G_0, J_G} E $ by
$f(x, h)=(x, \psi(x)\cdot h)$. Moreover, $f_H$ and $f_G$ are
simply the compositions of $\psi$ with the moment maps of $E$.

One can verify that this is a 1-1 correspondence on the level of
objects of these two stacks. The correspondence on the level of
morphisms is also easy to check.

Finally, from the construction above, it is not hard to see that
the moment maps are exactly the projections from $H_0\times_\cY
G_0$ to $H_0$ and $G_0$.  \hfill $\btd$

We sketch the proof of the second statement (which is not used in
the remainder of this paper). We have morphisms
$G_0\to\cX\overset{\phi}{\lra}\cY$ and $H_0\to\cY$. Take the
bibundle $E_\phi$ to be $G_0\times_{\cY} H_0$. It is not hard to
check that $E_\phi$ satisfies the required properties.
\end{proof}

In view of Proposition \ref {hs}, the fact that the composition of
HS morphisms is not associative can be understood by the fact that
compositions of 1-morphisms of stacks are associative only up to
2-morphisms of stacks.

\subsection{2-morphisms}
As HS morphisms correspond to morphisms in stacks, 2-morphisms
also have their exact counterparts in the language of Lie
groupoids. Recall that morphisms of stacks are just functors
between categories, and a 2-morphism of stacks between two
morphisms is  a natural transformation between these two morphisms
viewed as functors. We have 2-morphisms of groupoids defined as
following:
\begin{defi}[2-morphisms]
Let $(E^i, J^i_G, J^i_H)$ be two HS morphisms from the Lie
groupoid $G$ to $H$. A {\bf 2-morphism} from $(E^1, J^1_G, J^1_H)$
to $(E^2, J^2_G, J^2_H)$ is a bi-equivariant isomorphism from
$E_1$ to $E_2$.
\end{defi}
\begin{remark}\label{2morp}
\;\; \item[i)] If the two HS morphisms are given by groupoid
homomorphisms $f$ and $g$ between $G$ and $H$, then a 2-morphism
from $f$ to $g$ is just a smooth map $\alpha: G_0 \to H_1$ so that
$f(x)=g(x)\cdot \alpha(x) $ and $\alpha(\gamma x)=g(\gamma) \alpha
(x) f(\gamma)^{-1}$, where $x\in G_0$ and $\gamma\in G_1$. So it
is easy to see that not every two morphisms can be connected by a
2-morphism and when they can, the 2-morphism may not be unique
(for example, this happens when the isotropy group is nontrivial
and commutative).

\;\; \item[ii)] From the proof of Proposition \ref{hs}, one can
see that a 2-morphism between HS morphisms corresponds to a
2-morphism between the corresponding (1)-morphisms on the level of
stacks.
\end{remark}

\subsection{Invariant maps}
Invariant maps are a convenient way to produce maps between
stacks that we will use later in the construction of the Weinstein
groupoid.

\begin{lemma} \label{im} 
Given a Lie groupoid
$G_1\underset{\mathbf{t}}{\overset{\mathbf{s}}{\rightrightarrows}}G_0$
and a manifold $M$, any $G$-invariant map $f: G_0\to M$ induces a
morphism $\bar{f}: BG \to M$ between the differentiable stacks
such that $f=\bar{f}\circ \phi$, where $\phi: G_0 \to BG$ is the
covering map of atlases.
\end{lemma}
\begin{proof}
Since $f$ is $G$ invariant, $f$ introduces a morphism between Lie
groupoids:
($G_1\underset{\mathbf{t}}{\overset{\mathbf{s}}{\rightrightarrows}}G_0$)$\to
(M\rightrightarrows M)$. By Proposition \ref{hs} it gives a smooth
morphism between differentiable stacks. More precisely, let $Q\to
S$ be a $G$ torsor over $S$ with moment map $J_1$ and projection
$\pi_1$. Since the $G$ action on the $\pi_1$-fibre is free and
transitive, we have $S=Q\times_{f\circ J_1,id} M /G_1$. Notice
that a ($M\rightrightarrows M$)-torsor is simply a manifold $S$
with a smooth map to $M$. Then $\bar{f}(Q)$ is the morphism $J_2:
S\to M$ given by $J_2(s)=f\circ J_1(q)$, where $q$ is any preimage
of $s$ by $\pi$ (it is well defined since $f$ is $G$-invariant).
For any map $a: S\to G_0$, the image under $\phi$ is $Q_a=S
\times_{a,\bt} G_1$, and $\bar{f}(Q_a)$ is the map $f\circ a$
since $f$ is $G$-invariant. Therefore $f=\bar{f}\circ \phi$.
\end{proof}

\begin{lemma}\label{ssts} 
If a $G$ invariant map $f: G_0\to M$ is a submersion, then the
induced map $\bar{f}: BG \to M $ is a submersion of differentiable
stacks.
\end{lemma}

\begin{proof}
Let $U\to M$ be a morphism of manifolds. Using  base change of the
representable surjective submersion $G_0 \to BG$ by the projection
$BG\times_{\bar{f}, M} U \to BG$, we can see that $BG\times_M U$
is a differentiable stack with the atlas $G_0 \times_M U$. Note
that the composition $G_0\times_M U\to BG\times_M U\to U$ is a
submersion because it's the base change of $f: G_0\to M$ by $U\to
M$. Now take an atlas $V\to BG\times_M U$ which is a representable
surjective submersion.
$$\begin{diagram}
\node{G_0\times_M U\times_{BG\times_M U}V}\arrow[2]{e}\arrow{se}\node[2]{G_0\times_M U}\arrow[2]{s}\arrow{se}\arrow{ese}\\
\node[2]{V}\arrow[2]{e}\node[2]{BG\times_M U}\arrow{e}\arrow[2]{s}\node{U}\arrow[2]{s}\\
\node[3]{G_0}\arrow{se}\arrow{ese}\\
\node[4]{BG}\arrow{e}\node{M}
\end{diagram}
$$
We see that $G_0\times_M U\times_{BG\times_M U} V$ is a manifold
and the projections to $G_0\times_M U$ and $V$ are submersions.
The composition
$$G_0\times_M U\times_{BG\times_M U} V\to V\to U$$ coincides with
$$G_0\times_M U\times_{BG\times_M U} V\to G_0\times_M U\to
BG\times_M U\to U, $$ which is a surjective submersion. Hence
$V\to U$ is
a submersion.\\
\end{proof}

\begin{remark} \label{banach} It is not hard to see that the construction of stacks in the category of
smooth manifolds can be extended to the category of Banach
manifolds, yielding the notion of Banach stacks. Many properties
of differentiable stacks, including the ones discussed here, are
shared by Banach stacks as well. Also, the category of
differentiable stacks can be obtained from the category of Banach
stacks by restricting the base category.
\end{remark}

\section{Vector bundles}
\subsection{Vector bundles over stacks}

\begin{ep}\label{cat-vb}
Let $\cE$ be the category of all vector bundles. The morphisms in
$\cE$ are bundle maps. Then $\cE$ is a category over $\cC$. The
set of objects over $S\in \cC$ is the set of all vector bundles
over $S$. Then $\cE$ is a stack over $\cC$ whose the pullbacks (as in
the definition of stacks) are the ordinary pullbacks for vector
bundles. More precisely, given a morphism $f:S\to T$ and a vector
bundle $V\to T$, then $f^*T=S\times_T V$.
\end{ep}

\begin{defi}[Bundle functor] \label{def-vb}
A contravariant functor $\cF: \cX \to \cE$ is called a bundle
functor if
\begin{enumerate}
\item for every morphism $f: S\to T$, there is an isomorphism $\alpha_f$
from $\cF(T)\circ f^*$ to $S\times_T \cdot \circ \cF(S)$, i.e. the
following diagram is commute up to $\alpha_f$:
\[
\begin{CD}
\cX_T @>f>> \cX_S \\
@V\cF_T VV @VV \cF_S V \\
\cE_U @> S\times_T \cdot >> \cE_S.
\end{CD}
\]
\item for every two morphisms $f: S\to T$, $g: T\to R$, the
cocycle condition
\[ \alpha_{g\circ f} = (R \times_{g, T} \alpha_f) \circ g(f^*), \]
where we identify $R\times_{g, T}(T\times{f, S} \cdot )$ with $R
\times_{ g\circ f, S} \cdot $.
\end{enumerate}
\end{defi}

\begin{pdef}[vector bundles over stacks]
A vector bundle $\cV$ over a stack $\cX$ is a stack over $\cC$ along
with a  bundle functor $\cF: \cX\to \cE$ such that,
\begin{enumerate}
\item the set of objects over $S\in\cC$ is
\[ Obj(\cV)|_S=\{ (x, y): x\in \cX, y\; \text{is a global section
of}\; \cF(x)\}; \]
\item a morphism over $f: S\to T$ is an arrow from $(f^*x, y')$
to $(x, y)$ where $y'$ is determined by $y$ via the inverse of the
following isomorphism:
\[ \alpha_f (x): \cF(f^*(x))\to S\times_T \cF(x) .\]
\end{enumerate}
There is a projection $F: \cV \to \cX$ given by,
\[ (x, y)\mapsto x, \quad  ((f^*x, y') \to (x, y)) \mapsto (f^*x \to
x). \] With this projection $\cV$ is also a stack over $\cX$.
\end{pdef}
\begin{proof}
See \cite{lmb}, Chapter 14.
\end{proof}

\begin{pdef}[pull-backs] Let $ \cV$ be a vector bundle over a stack
  $\cY$ given by a bundle functor $\cF$ and $\phi: \cX \to \cY$ a morphism of stacks. Then $\phi\circ\cF$ is also a bundle functor $\cX\to \cE$
and the vector bundle given by it is called the pull-back vector
bundle $\phi^*\cV$ by the map $\phi$.
\end{pdef}
\begin{proof}
For an $f:S\to T$, we can choose $f^*$ (of $\cX$ according to that
of $\cY$) such that the following diagram commutes,
\[\begin{CD}
\cX_T @> f^* >> \cX_S \\
@V \phi VV @ VV\phi V \\
\cY_T @> f^* >> \cY_T.
\end{CD}
\]The rest of the proof follows by composition of diagrams.
\end{proof}

\subsection{Vector bundles over groupoids}

Let us first recall the following definition \cite{m-orbi}:
\begin{defi}[vector bundles over Lie groupoids]  \label{vbgpd}
A  vector bundle over a Lie groupoid $G$ is an equivariant $G$
vector bundle $V$ over $G_0$ such that the $G$ action (from the
right) is linear.
\end{defi}

\begin{pdef}[pull-backs via HS morphisms] \label{pullbackhs}
Let $(E, J_G, J_H)$ be an HS morphism from the Lie groupoid $G$ to
$H$. Let $V_H$ be a vector bundle over $H$ in the sense of
Definition \ref{vbgpd}. Then $V_H\times_{H_0} E/ H$ is a vector
bundle over $G$, where the $H$ action on $V_H\times_{H_0} E$ is
given by $(v, e)\cdot h=(vh,  eh)$. We define it as the pull back
of $V_H$ via $E$ and denote it by $E^*V_H$.
\end{pdef}
\begin{proof}
It is easy to see that $V_H \times_{H_0}E=J_H^*V_H$ is an
$H$-equivariant vector bundle over $E$. Since the $H$ action on
$E$ is free and proper, its action on  $V_H \times_{H_0} E$ is
free and proper too. Hence $V_H\times_{H_0} E/H$ is a manifold and
furthermore a vector bundle over $G_0=E/H$. Moreover $G$ acts on
it from right by $[( v, e)]\cdot g= [(v, g^{-1} e)]$. Clearly this
action is linear since $[(\lambda v, g^{-1} e)]=\lambda [(v,
g^{-1} e)]$.
\end{proof}
\begin{remark}
If the HS morphism is actually given by a groupoid homomorphism
$\phi: G\to H$, then it is easy to check that the pull-back by
$\phi$ viewed as an HS morphism is the same as the usual pull-back
of  vector bundles via $\phi|_{G_0}: G_0\to H_0$.
\end{remark}

\begin{lemma}\label{vb-moeq}
In the above setting, $V_H \times_{H_0} E$ with obvious
projections induces an HS morphism from the groupoid $E^*V_H
\times_{G_0} G_1 \rightrightarrows E^*V_H$ to $V_H \times_{H_0}
H_1 \rightrightarrows V_H$. Moreover, if $E$ is a Morita bibundle,
then $V_H \times_{H_0} E$ is also a Morita bibundle.
\end{lemma}
\begin{proof} It is easy to check that the fibre-wise free and transitive
  action of $H$ (resp. $G$) on $E$ gives the fibre-wise free and transitive
  action of $V_H \times_{H_0} H_1 \rightrightarrows V_H$ (resp. $E^*V_H \times_{G_0}
G_1 \rightrightarrows E^*V_H$) on $V_H \times_{H_0} E$.
\end{proof}

\subsection{Vector bundles over differentiable stacks}
Given a vector bundle $F: \cV\to \cX$, if the base stack $\cX$ is
a differentiable stack, one should have some finer requirements
for $\cV$ to be a vector bundle in a  ``differentiable'' fashion.

\begin{defi}[vector bundles over differentiable stacks]
A vector bundle over a differentiable stack $\cX$ is a vector
bundle $F: \cV \to \cX$ in the sense of stacks such that the map
$F$ is a representable surjective submersion.
\end{defi}

\begin{lemma}
A vector bundle over a differentiable stack is a differentiable
stack.
\end{lemma}
\begin{proof}
Let $F: \cV \to \cX$ be a vector bundle over the differentiable
stack $\cX$. Choose an atlas  $X_0$ of $\cX$. Then $\cV\times_\cX
X_0$ is a manifold since $F$ is a representable submersion and
$\cV\times_\cX X_0 \to \cV$ is a representable surjective
submersion because $X_0\to \cX$ is so.
\[
\begin{CD}
\cV\times_\cX X_0 @>>> X_0 \\
@VVV @VV  V \\
\cV @> F >> \cX.
\end{CD}
\]
\end{proof}

We have an alternative and more direct way to define the above
concept if we look more carefully into the definition of the
vector bundles over stacks. The new definition allows us to link
the vector bundles over differentiable stacks and  the vector
bundles over Lie groupoids.

\begin{defi}[vector bundles over  differentiable stacks]\label{vb-skst}
Let $\cX$ be an  differentiable stack. A vector
bundle $\cV$ on $\cX$ consists of the following set of data:
\begin{itemize}
\item for each  groupoid presentation $G$ of $X$, a vector bundle $V_G$ over
$G$,
\item for each commutative diagram
\begin{equation}\label{commute}
\begin{diagram}
\node{G_0}\arrow{se}\arrow[2]{e,t}{\varphi}\node[2]{H_0}\arrow{sw}\\
\node{}\node{\cX}\node{}
\end{diagram}
\end{equation}
with $G$ and $H$ groupoid presentations and $\varphi$ a strong
equivalence,  an isomorphism
$\alpha_{\varphi}:V_G\to\varphi^*V_H$.
\end{itemize}
The isomorphisms $\alpha_{\varphi}$ are required to satisfy
\begin{itemize}
\item the cocycle condition: for any three groupoid
  presentations: $G$, $H$, and $K$, and strong equivalences $\varphi$
  and $\psi$ which fit into a commutative
diagram
$$\begin{diagram}
\node{G_0}\arrow{e,t}{\varphi}\arrow{se}\node{H_0}\arrow{s}\arrow{e,t}{\psi}\node{K_0}\arrow{sw}\\
\node[2]{\cX}
\end{diagram}$$
we have
$$\alpha_{\psi\circ\varphi}=\varphi^*\alpha_{\psi}\circ\alpha_{\varphi}: V_G\to
(\psi\circ\varphi)^*V_K=\varphi^*(\psi^*V_K).$$
\end{itemize}
\end{defi}
\begin{remark} This definition is more like a definition for
``differentiable'' vector functors, nevertheless we don't make
distinguish between bundle functors and vector bundles here. One might wonder why $\alpha$ being an
isomorphism of vector bundles in this definition is enough to
encode it being an isomorphism of functors in Definition
\ref{def-vb}. The reason is that all the objects with the form
$(x, y)$ such that $x: P\to S$ is a $G$ torsor can be recovered by
$V_G$. Please see Proposition \ref{vb-gpd-sk} for more details.
\end{remark}





\begin{prop}\label{vb-gpd-sk}
Given a vector bundle $V_G$ over a Lie groupoid $G$, one can
construct a vector bundle $\cV$ over the differentiable stack
$\cX$ that $G$ presents, such that $\cV$ is presented by $V_G
\times_{G_0} G_1 \rightrightarrows V_G$.
\end{prop}
\begin{proof}
Let $\cV$ be the vector bundle constructed by the contra-variant
functor \[\cF: \cX \to \cE, \quad (x: P\to S ) \mapsto P\times_{J,
G_0, \pi} V_G / G_1,  \]where $J$ is the moment map of the $G$
torsor $x: P\to S$ and $\pi: V_G \to G_0$ is the projection. Then
on the level of objects, $\cV\times_\cX G_0$ consists of the elements
with the form $((x, y), \eta, S\to G_0)$ where $x: P\to S$ is a
$G$-torsor, $y$ is a global section of the vector bundle
$P\times_{J, G_0} V_G / G_1$ and $\eta$ is an isomorphism from $x$
to the trivial principal $G$-bundle $S \times_{G_0} G_1 \to S$.
Hence $x:P\to S$ is also a trivial principal $G$-bundle. Therefore
$\eta(y)$ is given by a map $S\to V_G$. Then
\[ \phi:  ((x, y), \eta, S\to G_0) \mapsto (\eta(y): S\to V_G), \]
and
\[ \varphi: (f: S\to V_G) \mapsto ((x=S\times_{\pi\circ f, G_0} V_G, y=(s,
f(s)), id, \pi\circ f S\to G_0 ), \] give the isomorphisms between
$\cV \times_\cX G_0$ and $V_G$. Notice that $\varphi\circ \phi$
and $id$  differ by a natural transformation given by the ``$\eta$
part'' of an element.
\end{proof}
\begin{remark}
Given any other presentation $H$ of $\cX$, let $V_H=E^*V_G$, where
$E$ is the Morita bibundle from $H$ to $G$. It is easy to see they
determine the same differentiable stack.
\end{remark}

\begin{prop}
Given a differentiable stack $\cX$ with the groupoid presentation
$G$, there is a 1-1 correspondence between the set of the vector
bundles over $\cX$ with the set of vector bundles over $G$.
\end{prop}
\begin{proof}
We adapt the notation in the previous proposition. One direction
is assured by the previous proposition. The other direction is
also true by taking $V_G=\cF(\bs: G_1\to G_0)$. It is not hard to
check that such an $V_G$ gives the same vector bundle as $\cF$.
\end{proof}

The following proposition tells us the relation between the
concept of pull-backs in the setting of stacks and groupoids.

\begin{prop}[pull-backs] Let $\cV$ be a vector bundle over a differentiable
stack $\cY$. Let $\phi: \cX\to \cY$ be a smooth morphism between
differentiable stacks. For every presentation $G$ of $\cX$, let
\[V_G := E^*V_H, \]where $H$ is a groupoid presentation
of $\cY$ and $E$ is the HS bibundle corresponding to $\phi$. Then
the vector bundle given by $V_G$ as above is $\phi^*\cV$
 the pull-back of $\cV$ via $\phi$.
\end{prop}
\begin{proof}
Here we show it for a strong morphism $\phi: G\to H$ since from this the general
case follows. The key of the proof is the observation that
\[\phi^*(V_H) = G_0 \times_{H_0} V_H \] is also \[\cF(\phi(\bs: G_1\to G_0)).\]\end{proof}



\subsection{Tangent ``bundles'' and tangent groupoids}

Here we put ``bundle'' in quotation marks because of the following
reason: the first obvious try to define the tangent bundle of a
stack is to construct the bundle functor by $\cF(x)=TS$ for any
element $x$ over $S$. However given a map $f: S\to T$, it is easy
to see that the pull-back of the tangent bundle of $T$ is not
always $TS$. However one can give the definition of tangent
bundles via the usage of groupoids, though it is not a
``chart-independent'' method. Via this method, one can see that a
tangent bundle is not always a vector bundle.


\begin{pdef} Let $\cX$ be a differentiable stack presented by a Lie
  groupoid $G$. The tangent bundle of $\cX$ is the differentiable
  stack presented by the tangent groupoid $TG$.
\end{pdef}
\begin{proof}
We have to show the definition does not depend on the choice of
groupoid representations. Let $H$ be another groupoid presentation
of $\cX$. Then $H$ and $G$ are Morita equivalent through a Morita
bibundle $(E, J_G, J_H)$. Then we claim that $(TE, TJ_G, TJ_H)$
gives the Morita equivalence between the tangent groupoids $TG$
and $TH$. To see this, notice that the groupoid action $\Phi$ of
$G$ on $E$,
\[\Phi: G\times E \to E, \]
lifts to the tangent groupoid by taking derivatives,
\[T\Phi: TG\times TE \to TE. \]
The fact that the action $\Phi$ is free and transitive fibre-wise
on $E$, i.e. $E\to H_0$ is a left principal $G$ bundle, is
equivalent to the fact that the map
\[ pr_2 \times \Phi: G \times E \to E\times E \]
is an isomorphism, where $pr_2$ is the projection to the second
factor. Then it is easy to see that
\[ \tilde{pr}_2 \times T\Phi: TG \times TE \to TE \times TE \]
is also an isomorphism, where $\tilde{pr}_2$ the projection to the
second factor. Hence $TE$ is a left principal $TG$-bundle with
moment map $TJ_G$. Similarly, $TE$ is also a right principal
$TH$-bundle with moment map $TJ_H$.
\end{proof}

In the case when $\cX$ is \'etale, take an \'etale presentation
$G$ of $\cX$. Then $TG=(TG_1\rightrightarrows TG_0)$ is simply the
action groupoid $TG_0 \times_{G_0} G_1 \rightrightarrows TG_0$. By
Proposition \ref{vb-gpd-sk}, it is a vector bundle over $\cX$. On
the other hand, if $TG=(TG_1\rightrightarrows TG_0)$ is in the
form of $V_G \times G_1 \rightrightarrows V_G$, then it has to be
the action groupoid $TG_0 \times_{G_0} G_1 \rightrightarrows TG_0$
which is true iff $G$ is an \'etale presentation of $\cX$, in
particular $\cX$ has to be \'etale.

\chapter{Weinstein groupoids}\label{w}
In this chapter, we will introduce the  new concept of Weinstein
groupoids and solve the Lie's third theorem for Lie algebroids.

\section{Path spaces} \label{pathspace}
We define the $A_0$-path space, which is very
similar to\footnote{Actually it is a submanifold of the $A$-path
space.} the $A$-paths defined in \cite{cf}.

Let us first recall the definition of a Lie algebroid.

\begin{defi}
\label{def:liealg} A {\em Lie algebroid} is a vector bundle $A\to
M$, together with a Lie algebra bracket $\br$ on the space of
sections $\Gamma(A)$ and a bundle map $\rho:A\to TM$, called the
{\em anchor}, satisfying
\[[a_1, fa_2]=f[a_1, a_2]+(\rho(a_1)f)a_2 \]
for any $a_1, a_2\in\Gamma(A), f\in\ci(M)$.
\end{defi}
\begin{remark}From this condition it follows that
for any $a_1, a_2\in\Gamma(A)$, $\rho [a_1, a_2]=[\rho a_1,\rho
a_2]$  \cite{kosmann}.
\end{remark}

\begin{defi} [$A_0$-paths and $A$-paths] \label{apath}
Given a Lie algebroid $A \overset{\pi} {\lra} M$ with  anchor
$\rho : A \to TM$, a $C^1$ map $a$: $I=[0, 1] \to A$ is an $A_0$-path if
\[ \rho (a(t)) = \frac{\di}{\di t} \left( \pi \circ a(t) \right) , \]
and it satisfies the following boundary conditions,
\[ a(0)=0,\; a(1)=0,\; \dot{a}(0)=0,\; \dot{a}(1)=0.\]
We often denote by $\gamma(t)$  the base path  $\pi\circ a(t)$ in
$M$. We denote $P_0(A)$ the set of all $A_0$-paths of $A$. It is a
topological space with topology given by uniform convergence of
maps. Omitting the boundary condition above, one gets the
definition of  $A$-paths, and we denote the space of
$A$-paths by $P_aA$.
\end{defi}

We can equip $P_0A$ with the structure of a smooth (Banach)
manifold using a Riemannian structure on $A$. On the total space
of $C^1$ path  $PA= C^1(I, A )$, there is a $\ci$-structure as
follows: at every point $a: I \to A$ in $PA$, let $a^*TA \to I$ be
the pull-back of the tangent bundle to $I$. For $\epsilon>0$, let
$T_{\epsilon}\subset a^*TA$ be the open set consisting of tangent
vectors of length less than $\epsilon$. For sufficiently small
$\epsilon$, we have the exponential map $\exp$: $T_{\epsilon} \to
A$, $(t, v) \mapsto \exp_{a(t)} v$. It maps $T_{\epsilon}$ to an
open subset of $A$. Using this map we can identify
$PT_{\epsilon}$, the  $C^1$-sections of $T_{\epsilon}$, with an
open subset of $PA$. The oriented vector bundle $a^*TA$ over $I$
is trivial. Let $\varphi : a^*TA \to I\times\R^n$ be a
trivialization where $n$ is the dimension of $A$. Then $\varphi$
induces a mapping from $PT_{\epsilon}$ to $P\R^n=C^1(I, \R^n)$.
Since $C^1(I, \R^n)$ is a Banach space with norm
$\|f\|^2=\sup\{|f|^2+ |f'|^2\}$, $PT_{\epsilon}$ can be used as a
typical Banach chart for the Banach manifold structure of $PA$.
$P_0A$ is defined by equations on $PA$ which, in the  local
charts above, can be written as
$$ \dot{\gamma}^k(t)= \sum_{j=1}^{m-n} \rho^k_j (\gamma(t)) a^j(t),
\;\;\; a^j(0)=a^j(1)=0, \; \dot{a}^j(0)=\dot{a}^j(1)=0, $$ for
$j=1, ..., n=\text{rank} A, k=1,..., m=\dim M.$ The space of
solutions is a closed subspace of $P(\R^n)$, hence is also a
Banach space and it gives a typical chart of $P_0A$. In this way,
$P_0A$ inherits the structure of a Banach manifold from $PA$.
We refer to \cite{lang} for the definition and further properties
of Banach manifolds.

\begin{pdef}\label{homotopy} Let $a(\epsilon, t)$ be a family of
$A_0$-paths of class $C^2$ in $\epsilon$ and assume that their
base paths $\gamma(\epsilon,t)$ have fixed end points. Let
$\nabla$ be a connection on $A$ with torsion $T_{\nabla}$ defined
as
\begin{equation} \label{eq:homotopy}
T_{\nabla} ( \alpha ,
\beta) = \nabla_{\rho(\beta)} \alpha - \nabla_{\rho(\alpha)} \beta
+ [\alpha, \beta].
\end{equation} Then the solution $b=b(\epsilon,
t)$ of the differential equation\footnote{Here, $T_{\nabla} (a,
b)$ is not quite well defined. We need to extend $a$ and $b$ by
sections of $A$, $\alpha$ and $\beta$, such that
$a(t)=\alpha(\gamma(t), t)$ and the same for $b$. Then $T_{\nabla}
(a, b)|_{\gamma(t)}:=$ $ T_{\nabla} (\alpha, \beta)|_{\gamma(t)}$
at every point on the base path. However, the choice of extending
sections does not affect the result.}
\begin{equation}\label{homotopyeq}
\partial_t b -\partial_{\epsilon }a = T_{\nabla} (a, b),
\;\;\;\;\; b(\epsilon ,0)=0 \end{equation} does not depend on the
choice of connection $\nabla$. Furthermore, $b(\cdot, t)$ is an
$A$-path for every fixed $t$, i.e. $\rho (b(\epsilon, t)) =
\frac{\di}{\di \epsilon} \gamma(\epsilon, t)$. If the solution $b$
satisfies $b(\epsilon,1)=0$, for all $\epsilon$, then $a_0$ and
$a_1$ are said to be  {\bf equivalent} and we write $a_0 \sim
a_1$.
\end{pdef}
\begin{remark} A homotopy of $A$-paths \cite{cf}  is
  defined by replacing $A_0$ by $A$ in the definition above. A similar result as above holds for $A$-paths \cite{cf}.  So the above statement holds viewing $A_0$-paths as $A$-paths.
\end{remark}

This flow of $A_0$-paths $a(\epsilon, t)$ generates a foliation
$\cF$.  The $A_0$-path
space is a Banach submanifold of the $A$-path space and $\cF$ is the
restricted foliation of the foliation defined in Section 4 of
\cite{cf}. For any foliation, there is an associated  {\bf
monodromy groupoid}\cite{moerdijk} (or {\bf fundamental groupoid}
as in \cite{cw}) : the objects are points in the manifold and the
arrows are paths within a leaf up to leaf-wise homotopy with fixed end
points. The source and target maps associate the
equivalence class of paths to the starting and ending points
respectively. It is a Lie groupoid in the sense of \cite{cf}, for
any regular foliation on a smooth manifold. In our case, it is an
infinite dimensional groupoid equipped with a Banach manifold
structure. Here, we slightly generalized the definition of Lie
groupoids to the category of Banach manifolds by requiring exactly
the same conditions but in the sense of Banach manifolds. Denote the
monodromy groupoid of $\cF$ by $Mon(P_0 A)
\underset{\mathbf{t}_M}{ \overset{\mathbf{s}_M}
{\rightrightarrows}} P_0 A$. In a very similar way, one can also
define the {\bf holonomy groupoid} $Hol(\cF)$ of $\cF$
\cite{moerdijk}: the objects are points in the manifold and the
arrows are equivalence classes of paths with the same holonomy.

Since $P_0 A$ is second countable, we can take an open cover $\{
U_i\}$ of $P_0 A$ which consists of countably many small enough
open sets so that in each chart $U_i$, one can choose a
transversal $P_i$ of the foliation $\cF$. By Proposition 4.8 in
\cite{cf}, each $P_i$ is a smooth manifold with dimension equal to
that of $A$. Let $P=\coprod P_i$ be the smooth immersed
submanifold of $P$. We can choose $\{U_i\}$ and transversal
$\{P_i\}$ to satisfy the following conditions:
\begin{enumerate}
\item If $U_i$ contains the constant path $0_x$ for some $x\in M$, then
$U_i$ has the transversal $P_i$ containing all constant paths
$0_y$ in $U_i$ for  $y\in M$.
\item If $a(t)\in P_i$ for some $i$, then $a(1-t)\in P_j$ for some
$j$.
\end{enumerate}
It is possible to meet the above two conditions: for (1) we refer
readers to Proposition 4.8 in \cite{cf}. There the result is for
$P_aA$. For $P_0A$, one has to use a smooth reparameterization
$\tau$ with the properties:
\begin{enumerate}
\item $\tau(t)=1$ for all $t \geq 1$ and $\tau(t)=0$ for all $t
  \leq 0$;
\item $\tau'(t)>0$ for all $t\in (0, 1)$.
\end{enumerate}
Then $a^\tau(t):= \tau(t)' a(\tau(t))$ is in $P_0A$ for all $a\in
P_aA$. $\phi_\tau: a\mapsto a^\tau$ defines an injective bounded
linear map from $P_aA \to P_0A$. Therefore, we can adapt the
construction for $P_a A$ to our case by using $\phi_\tau$.
For (2), we define a map $inv:\; P_0 A \to P_0 A $ by
$inv(a(t))=a(1-t)$. Obviously $inv$ is an isomorphism. In
particular, it is open. So we can add $inv(U_i)$ and $inv(P_i)$ to
the collection of open sets and transversals. The new collection
will have the desired property.

Restrict $Mon(P_0 A)$ to $P$. Then $Mon(P_0 A)|_P$ is a finite
dimensional \'etale Lie groupoid\footnote{An \'etale Lie groupoid
is a Lie groupoid such that the source (hence the target) map is a
local diffeomorphism.} \cite{mrw}, denoted by \gammapoid. If we
choose a different transversal $P'$, the restriction $\Gamma'$ of
$Mon(P_0 A)$ to $P'$ will be another finite dimensional \'etale
Lie groupoid. All these groupoids are related by Morita
equivalence: $\Gamma'$ is Morita equivalent to $\Gamma$ through
the  finite dimensional bibundle $\bs_M^{-1}(P) \cap
\bt_M^{-1}(P')$, where $\bs_M$ and $\bt_M$ are the source and target maps
of $Mon(P_0 A)$; $Mon(P_0 A)\rightrightarrows P_0 A $ is Morita
equivalent to \gammapoid through the Banach bibundle
$\bs_M^{-1}(P)$. One can do the same to $Hol(P_0 A)$ and get
a finite dimensional \'etale Lie groupoid, which we denote by
\hgammapoid. However, these groupoids are Morita equivalent to
each other in a similar way as their monodromy counterpart, but
not to the groupoids induced from $Mon(P_0 A)$.

We will build a Weinstein groupoid of
$A$ based on this path space $P_0 A$. One can interpret the ``identity section'' as the embedding
obtained by taking constant paths $0_x$, for all $x\in M$,  the
``inverse'' of a path $a(t)$ as $a(1-t)$, and the source and
target map $\bs$ and $\bt$ as taking the end points of the base
path $\gamma(t)$.  According to the two conditions above, these
maps are well-defined on the finite dimensional space $P$ as well.
Since reparameterization and projection are bounded linear
operators in Banach space $\ci (I, \R^n)$, the maps defined above
are smooth maps in $P_0 A$, hence in $P$. So we
 could almost make $P$ or $P_0A$ into a Lie groupoid, except that the
 multiplication has not been defined yet.

To define multiplication, notice that for any $A$-paths $a_1$,
$a_0$ in $P_0 A$ such that the base paths satisfy
$\gamma_0(1)=\gamma_1(0)$, one can define a ``concatenation''
\cite{cf}:
\[  a_1\odot a_0 =\left\{
 \begin{array}{rr}
 2a_0(2t), \quad &\mbox{$0\leq t\leq \half$}\\
 2a_1(2t-1),\quad &\mbox{$\half < t\leq 1$}
 \end{array}\right.   \]

Concatenation is a bounded linear operator in the local charts,
hence is a smooth map. However it is not associative. Moreover it is
not well-defined on $P$. If we quotient out by the  equivalence
relation induced by $\cF$, concatenation is associative and well-defined.
However, after quotienting out by the equivalence, we may not end up
with a smooth manifold any more. To overcome the difficulty, our
solution is to pass to the world of differentiable stacks.

\section{Construction}

Recall that in Section \ref{pathspace}, given a Lie algebroid $A$, we
constructed an \'etale groupoid \gammapoid. Hence  we can construct an \'etale differentiable stack $\cG(A)$
presented by \gammapoid. If we choose a different transversal
$P'$, the restriction of $Mon(P_0 A)$ on $P'$, $\Gamma'$, is
Morita equivalent to $\Gamma$. As we have seen, this implies that
they present isomorphic differentiable stacks. Therefore, we
might as well base our discussion on \gammapoid.

Moreover, $Mon(P_0 A)\rightrightarrows P_0 A $ is Morita
equivalent to \gammapoid. So $\cG(A)$ can also be presented by
$Mon(P_0A)$  as a Banach stack.

In this section, we will construct two Weinstein groupoids
$\cG(A)$ and $\cH(A)$ for every Lie algebroid $A$ and  prove
Theorem \ref{integ}.

\noindent {\bf Theorem \ref{integ}.} {\em
A Lie algebroid $A$ is integrable in the classical sense iff
$\cH(A)$ is representable, i.e. it is an ordinary manifold. In
this case $\cH(A)$ is the source-simply connected Lie groupoid of
$A$ (it is also called the Weinstein groupoid of $A$ in
\cite{cf}).}

We begin with $\cG(A)$. We first define the inverse, identity
section, source and target maps on the level of groupoids.

\begin{defi}
Define
\begin{itemize}
  \item $i:$ (\gammapoid)$\to$(\gammapoid) by
  $g=[a(\epsilon, t)]\mapsto [a(\epsilon, 1-t)]$, where $[\cdot]$ denotes the homotopy class in $Mon(P_0A)$;
  \item $e: M \to$(\gammapoid) by $x\mapsto 1_{0_x}$, where
  $1_{0_x}$ denotes the identity homotopy of the constant path
  $0_x$;
  \item $\bs:$(\gammapoid)$\to M$ by $g=[a(\epsilon,
  t)]\mapsto \gamma(0, 0)(=\gamma(\epsilon, 0),\, \forall \epsilon)$,
  where $\gamma$ is the base path of $a$;
  \item $\bt:$(\gammapoid)$\to M$ by $g=[a(\epsilon,
  t)]\mapsto \gamma(0, 1)(=\gamma(\epsilon, 1),\, \forall \epsilon)$;
\end{itemize}
These maps can be defined similarly on $Mon(P_0A)\rightrightarrows
P_0A $. These maps are all bounded linear maps in the local charts
of $Mon(P_0A)$. Therefore they are smooth homomorphisms between
Lie groupoids. Hence, they define smooth morphisms between
differentiable stacks. We denote the maps corresponding to $i$,
$\epsilon$, $\bs$, $\bt$ on the stack level by $\bar{i}$,
$\bar{\epsilon}$, $\bar{\bs}$ and $\bar{\bt}$.
\end{defi}

\begin{lemma}\label{e-immersion}
The maps $\bbs$ and $\bbt$ are surjective submersions. The map
$\bar{e}$ is a monomorphic immersion. The map $\bar{i}$ is an
isomorphism.
\end{lemma}
\begin{proof}
$\bs$ and $\bt$ restricted to $P$ are $\Gamma$-invariant
submersions because any path through $x$ in $M$ can be lifted to a
path in $P$ passing through any given preimage of $x$. According
to Lemma \ref{im} and \ref{ssts}, the induced maps $\bar{\bs}$ and
$\bbt$ are  submersions.

Denote by $e_0: M \to P$ the restricted map of $e$ on the level of objects. Notice that $e_0$ fits into the following diagram
(which is not  commutative):
\begin{equation}\label{eq: immersion}\begin{diagram}
\node{M\times_{\cG(A)}P}\arrow{e,t}{pr_2}\arrow{s,l}{pr_1}
\node{P}\arrow{s,r}{\pi} \\
\node{M}\arrow{e,t}{\bar{e}}\arrow{ne,t}{e_0}
\node{\cG(A)}\end{diagram}
\end{equation}
Given $x=(f:U\to M)\in M$, as a   $G$-torsor we have $\bar{e}(x)=U\times_{e_0\circ f,
G_0} G_1$ and $e_0(x)=(e_0\circ f:U\to G_0)\in
G_0$. Given $y=(g: U\to G_0)\in G_0$, we have
$\pi(y)=U\times_{g, G_0} G_1$. A typical object of $M_i
\times_{\cG} G_0$ is $(x, \eta, y)$ where $\eta$ is a morphism of
$G$-torsors from $\bar{e}(x)$ to $\pi(y)$ over $id_U$ of $U$. Then
by the equivariancy of $\eta$, we have a map $\phi$: $U\to G_1$,
such that $e_0\circ f=g\cdot \phi$. Therefore, we have a map
$\alpha: M\times_{\cG(A)} G_0\to G_1$ given by $\alpha (x, \eta,
y)=\phi$, such that \[e_0 \circ pr_1 =pr_2 \cdot \alpha.\] Since
$\pi$ is \'etale, so is $pr_1$. Moreover, since $e_0$ is an
embedding, $pr_2$ must be an immersion. Therefore, by Lemma
\ref{immersion}, $\bar{e}$ is an immersion.

As $\bs \circ e = \bt \circ e = id$ holds on the level of
groupoids, this identity passes to an identity on the level of
differentiable stacks too. Since $\bbs \circ \bar{e}=\bbt \circ
\bar{e}=id$, it is easy to see that $\bar{e}$ must be monomorphic
and $\bbs$ (and $\bbt$) must be epimorphic.

The map $i$ is an isomorphism of groupoids, hence it induces an
isomorphism at the level of stacks.

\end{proof}

We define the multiplication for the infinite dimensional
presentation $Mon(P_0 A)$ with source and target maps $bs_M$ and
$\bt_M$. First we extend ``concatenation'' to $Mon(P_0A)$.
Consider two elements $g_1, g_0 \in Mon(P_0 A)$ whose base paths
on $M$ are connected at the end points. Suppose $g_i$ is
represented by $a_i(\epsilon, t)$. Define
\[ g_1\odot g_0=[a_1(\epsilon, t) \odot_t a_0
  (\epsilon , t)], \]
where $\odot_t$ means concatenation with respect to  the parameter
$t$ and $[\cdot]$ denotes the equivalence class of homotopies.

Notice that $\bs\circ\bs_M =\bs\circ\bt_M$ and $\bt\circ\bs_M =
\bt\circ\bt_M$ are surjective submersions. By Lemma \ref{pa} and Example \ref{pa2},
\[Mon(P_0A)\times_{\bs\circ\bs_M, M, \bt\circ\bt_M}Mon(P_0A)\rightrightarrows
P_0A \] with source and target maps $\bs_M \times \bs_M$ and
$\bt_M\times\bt_M$ is a Lie groupoid and it presents the stack
$$\cG
\times_{\bbs,M, \bbt} \cG.$$

Finally let $m$ be the following smooth homomorphism between
Lie groupoids:
$$\label{multbanach}
\xymatrix@=50pt{Mon(P_0 A)\mathop\times\limits_{\bs\circ\bs_M, M, \bt\circ\bt_M}Mon(P_0 A)\ar[d]^{\bs_M\times \bs_M}\ar@<-1ex>[d]_{\bt_M\times\bt_M}\ar[r]^-{\odot}& Mon(P_0 A)\ar[d]^{\bs_M}\ar@<-1ex>[d]_{\bt_M} \\
P_0 A \times P_0 A\ar[r]^{\odot}& P_0 A \\
}
$$

Multiplication is less obvious obvious for the \'etale presentation
$\Gamma\rightrightarrows P$. We will have to define the multiplication
through an HS morphism.

Viewing $P$ as a submanifold of $P_0A$, let $E=\bs_M^{-1}(P)\cap
\bt_M^{-1}(m(P\times_M P))\subset Mon(P_0A)$. Since $\bs_M$ and
$\bt_M$ are surjective submersions and $m(P\times_M P)\cong
P\times_M P$ is a submanifold of $P_0A$, $E$ is a smooth manifold.
Since $P$ is a transversal, $\bt_M: E\to m(P\times_M P)$ is
\'etale. Moreover $\dim m(P\times_M P)= 2\dim P -\dim M$.  So $E$
is finite dimensional. Further notice that $m: P_0A\times P_0A \to
P_0A $ is injective and its ``inverse'' $m^{-1}$ defined on the
image of $m$ is given by
\[ m^{-1}: b(t)\mapsto (b(2t_1), b(1-2t_2)) \;\;\; t_1\in [0,
\frac{1}{2}],\; t_2\in [\frac{1}{2}, 1] \] which is bounded linear in
a local chart. Let $\pi_1=m^{-1}\circ \bt_M: E\to P\times_M P$ and
$\pi_2=\bs_M: E\to P $. Then it is routine to check that $(E,
\pi_1, \pi_2)$ is an HS morphism from $\Gamma\times_M \Gamma
\rightrightarrows P\times_M P$ to $\Gamma \rightrightarrows P$. It
is not hard to verify that on the level of stacks $(E, \pi_1, \pi_2)$ and $m$ give two
1-morphisms differed by a 2-morphism. Thus,
after modifying $E$ by this 2-morphism, we get another HS-morphism
$(E_m, \pi'_1, pi'_2)$ which presents the same map as $m$.
Moreover, $E_m\cong E$ as bibundles.

Therefore, we have the following definition:

\begin{defi} Define $\bm: \cG(A)\times_{\bbs,\bbt} \cG(A) \to \cG(A)$  to be
the smooth morphism between \'etale stacks presented by $(E_m,
\pi'_1, \pi'_2) $.
\end{defi}

\begin{remark}
If we use $Mon(P_0A)$ as the presentation, $\bm$ is also presented
by $m$.
\end{remark}

\begin{lemma} The multiplication $\bm: \cG(A) \times_{\bs, \bt} \cG(A) \to \cG(A)$ is
a smooth morphism between \'etale stacks and is associative up to
a 2-morphism, that is, the diagram
$$\begin{diagram}
\node{\cG(A)\mathop\times\limits_{\bs,\bt}\cG(A)\mathop\times\limits_{\bs,\bt}\cG(A)}\arrow{e,t}{id\times\bar{m}}\arrow{s,b}{\bar{m}\times id}\node{\cG(A)\mathop\times\limits_{\bs,\bt}\cG(A)}\arrow{s,r}{\bar{m}}\\
\node{\cG(A)\mathop\times\limits_{\bs,\bt}\cG(A)}\arrow{e,t}{\bar{m}}\node{\cG(A)}
\end{diagram}
$$
is 2-commutative, i.e. there exists a 2-morphism $\alpha:
\bm\circ(\bm\times id) \rightarrow \bm\circ(id\times \bm)$.
\end{lemma}
\begin{proof}
We will establish  the 2-morphism on the level of Banach stacks.
Notice that a smooth morphism in the category of Banach manifolds
between finite dimensional manifolds is a smooth morphism in the
category of finite dimensional smooth manifolds. Therefore, the
2-morphism we will establish gives a 2-morphism for the \'etale
stacks.

Take the Banach presentation $Mon(P_0 A)$, then $\bm$ can simply
be presented as a homomorphism between groupoids as in
\eqref{multbanach}. According to Remark \ref{2morp}, we now
construct a 2-morphism $\alpha: P_0 A\times_M P_0 A \times_M P_0 A
\to Mon(P_0 A)$ in the following diagram
\[
\xymatrix@=50pt{Mon(P_0 A)\mathop\times\limits_{M}Mon(P_0 A)\mathop\times\limits_{M}Mon(P_0 A)\ar[r]^-{m \circ (m \times id)}_-{m \circ (id \times m)}\ar[d]^{\bs_M\times\bs_M\times\bs_M}\ar@<-1ex>[d]_{\bt_M\times\bt_M\times\bt_M}& Mon(P_0 A)\ar[d]^{\bs_M}\ar@<-1ex>[d]_{\bt_M}\\
P_0 A\times_M P_0 A \times_M P_0 A \ar[r]& P_0 A }
\]

Let $\alpha(a_1, a_2, a_3)$ be the natural rescaling between $a_1
\odot (a_2\odot a_3)$ and $(a_1\odot a_2)\odot a_3$. Namely,
$\alpha(a_1, a_2, a_3)$ is the homotopy class
represented by \begin{equation}\label{rp}
a(\epsilon, t) =
((1-\epsilon)+\epsilon\sigma'(t))a((1-\epsilon)t+\epsilon\sigma(t)),
\end{equation}
where $\sigma(t)$ is a smooth reparameterization such that
$\sigma(1/4)=1/2, \; \sigma(1/2)=3/4.$ In local charts, $\alpha$
is a bounded linear operator, therefore, it is a smooth morphism
between Banach spaces. Moreover, $ m \circ (m \times id)=  m \circ
(id \times m)\cdot\alpha $. Therefore $\alpha$ serves as the
desired 2-morphism.
\end{proof}

One might be curious about whether there is any further obstruction to
associativity. There are six ways to  multiply four elements in
$\cG(A)$. Let
\[
\begin{split}
F_1&=\bm \circ \bm \times id \circ \bm \times id \times id, \\
F_2&=\bm \circ id \times \bm \circ \bm \times id \times id, \\
F_3&=\bm \circ \bm \times id \circ id \times id \times \bm, \\
F_4&=\bm \circ id \times \bm \circ id \times id \times \bm, \\
F_5&=\bm \circ id \times \bm \circ id \times \bm \times id, \\
F_6&=\bm \circ \bm \times id \circ id \times \bm \times id.
\end{split}
\]
They are morphisms fitting into the following commutative cube.
$$
\xymatrix@=5pt{
     & & \cG(A)\mathop\times\limits_{M}\cG(A)\mathop\times\limits_{M}\cG(A) \ar[dr]^{\bar{m}\times id} \ar[ddd]^{id\times\bar{m}}& \\
  \cG(A)\mathop\times\limits_{M}\cG(A)\mathop\times\limits_{M}\cG(A)\mathop\times\limits_{M}\cG(A) \ar[urr]^{id\times id\times\bar{m}} \ar[dr]^{\bar{m}\times id\times id} \ar[ddd]_{id\times\bar{m}\times id} & & & \cG(A)\mathop\times\limits_{M}\cG(A)\ar[ddd]^{\bar{m}} \\
     & \cG(A)\mathop\times\limits_{M}\cG(A)\mathop\times\limits_{M}\cG(A)\ar[urr]^{id\times\bar{m}} \ar[ddd]^{\bar{m}\times id} & & & \\
     & & \cG(A)\mathop\times\limits_{M}\cG(A) \ar[dr]^{\bar{m}} & \\
  \cG(A)\mathop\times\limits_{M}\cG(A)\mathop\times\limits_{M}\cG(A) \ar[urr]^{id\times\bar{m}} \ar[dr]^{\bar{m}\times id} & & & \cG(A) \\
     & \cG(A)\mathop\times\limits_{M}\cG(A) \ar[urr]^{\bar{m}} & &
       }
$$

There is a 2-morphism on each face of the cube to connect $F_i$
and $F_{i+1}$ ($F_7=F_1$), constructed as in the last Lemma. Let
$\alpha_i: F_i\to F_{i+1}$. Will the composition $\alpha_6 \circ
\alpha_6 \circ ... \circ \alpha_1$ be the identity 2-morphism? If
so, given any two different ways of multiplying four (hence any
number of) elements, different methods to obtain 2-morphisms
between them will give rise to the same 2-morphism. Since
2-morphisms between two 1-morphisms are not unique if our
differential stacks are not honest manifolds,  it is necessary to
study the further obstruction.

\begin{prop} \label{3asso} There is no further obstruction to the associativity of
$\bm$ in $\cG(A)$.\end{prop}
\begin{proof} Choose $Mon(P_0A)$ as the presentation
of $\cG(A)$, then the $\alpha_i$'s constructed above can be
explicitly written out as a smooth morphism: $P_0A\times_M
P_0A\times_M P_0A\times_M P_0A \to Mon(P_0A)$. More precisely,
according to the lemma above, $\alpha_i(a_1, a_2, a_3, a_4)$ is
the natural rescaling between $F_i(a_1, a_2, a_3, a_4)$ and
$F_{i+1}(a_1, a_2, a_3, a_4)$. Here by abuse of notation, we
denote by $F_i$ also the homomorphism on the groupoid level. It is
not hard to see that $\alpha_6 \circ \alpha_6 \circ ... \circ
\alpha_1$ is represented by a rescaling that is homotopic to the
identity homotopy between $A_0$-paths.

Therefore, the composed 2-morphism is actually the identity since
$Mon(P_0A)$ is made up by the homotopy of homotopy  of
$A_0$-paths. We also notice that the identity morphism in the
category of Banach manifolds between two finite dimensional
manifolds is the identity morphism in the category of finite
dimensional smooth manifolds. Therefore, there is no further
obstruction even for 2-morphisms  of  \'etale stacks.

\end{proof}

Now to show that $\cG(A)$ is a Weinstein groupoid that we have
defined in the introduction, we only have to show that the
identities in item (4) and (5) in Definition \ref{wgpd} hold and
the 2-morphisms in these identities are identity 2-morphisms when
restricted to $M$. Notice that for any $A_0$-path $a(t)$, we have
\[ a(t) \odot_t 1_{\gamma(0)} \sim a(t), \;\; a(1-t) \odot_t a(t)
\sim \gamma(0) ,\] where $\gamma$ is the base path of $a(t)$.
Using i) in Remark \ref{2morp}, we can hence see that on the
groupoid level $m\circ ((e\circ \bt) \times id)$ and $id$ only
differ by a 2-morphism, and the same for $m\circ (i \times id)$
and $e\circ \bs$. Therefore the corresponding identities hold on
the level of differentiable stacks.  Moreover, the 2-morphisms (in
all presentations of $\cG(A)$ we have described above) are formed
by rescalings. But when they restrict to constant paths in $M$,
they are just $id$.

Summing up what we have discussed above, $\cG(A)$ with all of the
structures we have given is a Weinstein groupoid over $M$.

We further comment that one can construct another natural
Weinstein groupoid $\cH(A)$ associated to $A$ exactly in the same
way as $\cG(A)$ by the Lie groupoid $Hol(P_0 A)$ or
$\Gamma^h\underset{\mathbf{t}_1}{\overset{\mathbf{s}_1}{\rightrightarrows}}P$
since they are Morita equivalent by a similar reason as for their
monodromy counterparts. One can establish the identity section,
the inverse, etc., even the multiplication in exactly the same
way. One only has to notice that in the construction of the
multiplication, the 2-morphism in the associativity diagram is the
holonomy class (not the homotopy class) of the reparameterization
\eqref{rp}. One can do so because homotopic paths have the same
holonomy. Moreover, for the same reason, there is no further
obstruction for the multiplication on $\cH(A)$, too.

The
integrability of $A$ and the representability of $\cG(A)$ are not
exactly the same, due to the presence of isotropy groups. But,
since holonomy groupoids are always effective \cite{moerdijk}, the
integrability of $A$ is equivalent to the representability of
$\cH(A)$ (see Theorem \ref{integ-g} and Theorem \ref{integ}).

\begin{pdef}[orbit spaces]
Let $\cX$ be a differentiable stack presented by Lie groupoid
$X=(X_1\rightrightarrows X_0)$. The  orbit space of $\cX$ is
defined as the topological quotient $X_0/X_1$. Throughout the
paper, when we say that the orbit space is a smooth manifold, we
mean that it has the natural smooth manifold structure induced
from $X_0$.
\end{pdef}

\begin{proof}
We have to show that  the topological quotient is independent of the
choice of presentations. Suppose that there is another
presentation $Y$ which is Morita equivalent to $X$ through $(E,
J_X, J_Y)$. Let $O_x$ be the orbit of $X_1$ in $X_0$ through point
$x$. By the fact that both groupoid actions are free and
transitive fiber-wise, $J_Y\circ J_X^{-1}(O_x)$ is another orbit
$O_y$ of $Y$. In this way, there is a 1-1 correspondence between
orbits of $X$ and $Y$. Hence, $Y_0/Y_1$ understood as the space of
orbits is the same as $X_0/X_1$ (i.e. the projection $X_0\to X_0/X_1$ is
smooth).
\end{proof}

\begin{thm} \label{integ-g} A Lie algebroid $A$ is integrable in the classical sense,
  i.e. there is a Lie groupoid whose Lie algebroid is $A$, iff the
  orbit space of $\cG(A)$ is a smooth manifold. Moreover, in this
  case, the orbit space of $\cG(A)$ is the unique source-simply
  connected Lie groupoid integrating $A$.
\end{thm}
\begin{proof}First  let $Mon(P_a A)$ be the
monodromy groupoid of the foliation induced by homotopy of
$A$-paths in Section \ref{apath}. We will show that $Mon(P_a A)$
is Morita equivalent to $Mon(P_0 A)$. Notice that $P_0A$ is a
submanifold of $P_a A$, so there is another groupoid $Mon(P_a
A)|_{P_o A}$ over $P_0 A$. We claim it is the same as $Mon(P_0
A)$. Namely, an A-homotopy $a(\epsilon, t)$ between two $A_0$
paths $a_0$ and $a_1$ is homotopic to an $A_0$-homotopy
$\tilde{a} (\epsilon, t)$ between $a_0$ and $a_1$. The idea is to
divide $\tilde{a}$ into three parts:
\;\;\item[i)] First deform $a_0$ to $a_0^\tau$ through $a_0(\epsilon, t)$,
  which is defined as $$(1-\epsilon
  +\epsilon \tau' (t)) a_0((1-\epsilon) t +\epsilon \tau(t)), $$ where $\tau$ is the
  reparameterization induced in Section \ref{apath};
\;\;\item[ii)] Then, deform $a_0^\tau$ to $a_1^{\tau}$ through $a(\epsilon,
  t)^\tau$;
\;\;\item[iii)] At last, connect $a_1^{\tau}$ to $a_1$ through
  $a_1(\epsilon, t)$, which is defined as $a_1((1-\epsilon)
  \tau'(t) +\epsilon) a_1 (\epsilon t +(1-\epsilon) \tau (t))$.
Then connect these three pieces by a similar method as in the
  construction of concatenations (though it might be piecewise smooth
  at the joints). Obviously, $\tilde{a}$ is a homotopy
  inside $A_0$-paths and it is homotopic to a rescaling (over $\epsilon$) of $a(\epsilon, t)$ through
the concatenation of $a_0((1-\lambda)\epsilon, t)$ and
  $(\lambda+(1-\lambda)\tau'(t))a(\epsilon,
  \lambda+(1-\lambda)\tau'(t))$ and $a_1((1-\lambda) \epsilon +
  \lambda, t)$. Eventually, we can smooth out everything to
  make the homotopy and  the homotopy of homotopies both smooth so that
  they are as desired.

It is routine to check that  $Mon(P_a A)|_{P_0 A}$ is Morita
equivalent to $Mon(P_a A)$ through the bibundle $\bt^{-1} (P_0
A)$, where $\bt$ is the target of the new groupoid $Mon(P_aA)$.

So the orbit space of $\cG(A)$ can be realized as $P_0 A/ Mon(P_0
A)$ which is isomorphic to $P_a A/ Mon(P_a
A)$. The rest of the proof follows from the main result in \cite{cf}, $P_a A/ Mon(P_a
A)$ is a smooth manifold iff $A$ is integrable and if so, $P_a A/
Mon(P_a A)$ is the unique source-simply connected Lie groupoid
integrating $A$.

\end{proof}

\noindent {\bf Theorem \ref{integ}.} {\em
A Lie algebroid $A$ is integrable in the classical sense iff
$\cH(A)$ is representable, i.e. it is an ordinary manifold. In
this case $\cH(A)$ is the source-simply connected Lie groupoid of
$A$ (it is also called the Weinstein groupoid of $A$ in
\cite{cf}).}

\begin{proof} [Proof for Theorem \ref{integ}]

According to \cite{moerdijk}, if the orbit space of a holonomy groupoid  is a manifold
then it is Morita equivalent to the holonomy groupoid itself.

Hence a differentiable stack $\cX=BG$  presented by a holonomy
groupoid $G$ is representable if and only if the orbit space $G_0/
G_1$ is a smooth manifold. One direction is obvious because
$G_0/G_1 \rightrightarrows G_0/G_1 $ is Morita equivalent to $G=
(G_1 \rightrightarrows G_0)$ if the orbit space is a manifold. The
other implication is not hard either, if one examines the Morita
equivalence diagram of $G$ and $\cX \rightrightarrows \cX$. The
Morita bibundle has to be $G_0$ since $\cX$ is a manifold.
Therefore $G_0$ is a principal $G$ bundle over $\cX$. This implies
that $G_0/G_1$ is the manifold $\cX$.

Notice that in general the orbit spaces of monodromy groupoids
and holonomy groupoids of a foliation are the same. By Theorem
\ref{integ-g} and the argument above, we conclude that  $A$ is
integrable iff $\cH(A)$ is representable and in this case,
$\cH(A)$ is $$P_0A/Hol(P_0 A)=P_0A/Mon(P_0 A)=P_aA/Mon(P_a A),$$ the unique source-simply connected
Lie groupoid integrating $A$.
\end{proof}

Recall Theorem \ref{cf}:

\noindent {\bf Theorem \ref{cf}.}{\em
As topological spaces, the orbit spaces of $\cH(A)$ and $\cG(A)$
are both isomorphic to the universal topological groupoid of $A$
constructed in \cite{cf}.}

Combining the proofs of Theorem \ref{integ-g} and Theorem
\ref{integ}, Theorem \ref{cf} follows naturally.

So far we have constructed $\cG(A)$ and $\cH(A)$ for every Lie
algebroid $A$ and verified that they are Weinstein groupoids. Thus
we have proved half of Theorem \ref{lieIII}. For the other half of
the proof, we will first introduce some properties of Weinstein
groupoids.

\section{Weinstein groupoids and local groupoids} In this
section, we examine the relation between abstract Weinstein
groupoids and local groupoids. A {\bf local Lie groupoid} $G_{loc}$ is an object
satisfying all axioms of a Lie groupoid except that the multiplication
is defined only near the local section. Namely there is a neighborhood
of $U$ the identity section $M$ such that the multiplication is defined
from $U\times_M U \to G_{loc}$.


Let us first show a useful lemma.

\begin{lemma}\label{eie} 
Given any \'etale atlas $G_0$ of $\cG$, there exists an open
covering $\{ M_l\}$ of $M$ such that the immersion $\bar{e}: M \to
\cG$ can be lifted to  embeddings $e_l: M_l \to G_0$. On the
overlap $M_l\cap M_j$, there exists an isomorphism $\varphi_{lj}$:
$e_j(M_j\cap M_l)\to e_l(M_j\cap M_l)$, such that $\varphi_{lj}
\circ e_j =e_l$ and $\varphi_{lj}$'s satisfy cocycle conditions.
\end{lemma}
\begin{proof}
Let $(E_e, J_M, J_G)$ be the HS-bibundle presenting the immersion
$\bar{e}: M \to \cG$. As  a right $G$-principal bundle over $M$,
$E_e$ is locally trivial, i.e. we can pick an open covering
$\{M_l\}$ so that $J_M$ has a section $\tau_l: M_l \to E_e$ when
restricted to $M_l$. Since $\bar{e}_l:=\bar{e}|_{M_l}$ is an
immersion (the composition of immersions $M_l\to M$ and $\bar{e}$
is still an immersion), it is not hard to see that $pr_2: M_l
\times_{\cG} G_0 \to G_0$ transformed by base change $G_0 \to \cG$
is an immersion. Notice that $e_l=J_G \tau_l: M_l \to G_0$ fits
into a diagram similar to  \eqref{eq: immersion}:
\[
\begin{diagram}
\node{M_l\times_{\cG}G_0}\arrow{e,t}{pr_2}\arrow{s,t}{pr_1}
\node{G_0}\arrow{s,r}{\pi} \\
\node{M_l}\arrow{e,t}{\bar{e}_l}\arrow{ne,t}{e_l} \node{\cG}
\end{diagram}
\]
Following a similar argument as in the proof of Lemma
\ref{e-immersion}, we can find a map $\alpha: M_l\times_{\cG}G\to
G_1$ such that
\[ e_l\circ pr_1= pr_2 \cdot \alpha.\]
Since $\pi$ is \'etale, so is $pr_1$. Therefore $e_l$ is an
immersion.

Since an immersion is locally an embedding, we can choose an open
covering $M_{ik}$ of $\{ M_l\}$ so that $e_l|_{M_{ik}}$ is
actually an embedding. To simplify the notation, we can choose a
finer covering $\{ M_l\}$ at the beginning and make $e_l$ an
embedding. Moreover, using the fact that $G$ acts on $E_e$
transitively (fiber-wise), it is not hard to find a local
bisection $g_{lj}$ of $G_1:=G_0\times_{\cG} G_0$, such that $e_l
\cdot g_{lj} =e_j$. Then $\varphi_{lj} = \cdot g^{-1}_{lj}$
satisfies $\varphi_{lj} \circ e_j = e_l$. Since  the $e_l$'s are
embeddings, the $\phi_{lj}$'s naturally satisfy the cocycle condition.
\end{proof}

Before the proof of Theorem \ref{local}, we need a local
statement.

\begin{thm}\label{local-local}
For every Weinstein groupoid $\cG$, there exists an open covering
$\{ M_l\}$ of $M$ such that one can associate a local Lie groupoid
$U_l$ over each open set $M_l$.
\end{thm}
\begin{proof} Let $\cG$ be presented by $G=(G_1\rightrightarrows G_0)$, and
$\{M_l\}$ be an open covering as in Lemma \ref{eie}.

Let $(E_m, J_1, J_2)$ be the HS bibundle from $G_1\times_M G_1
\rightrightarrows G_0 \times_M G_0$ to $G$ which presents the
stack morphism $\bm: \cG\times_{M} \cG \to \cG$.


Notice that $M$ is the identity section, i.e.
\[
\begin{array}{ccc}
M_l \times_M M_l(=M_l) & \overset{ \bm=id}{\lra} &  M_l \\
\downarrow &\curvearrowright & \downarrow \\
\cG\times_M\cG &\overset{\bm} {\lra} &\cG.
\end{array}
\]
Translate this commutative diagram into groupoids. Then the
composition of HS morphisms
\begin{equation}\label{2hs}
\xymatrix{M_l\times_M M_l(=M_l)\ar[r]\ar[dd]\ar@<-1ex>[dd]& G_1\times_M G_1\ar[dd]\ar@<-1ex>[dd]&  &G_1\ar[dd]\ar@<-1ex>[dd]\\
& & E_m\ar[dr]^{J_2}\ar[dl]_{J_1}& &\\
M_l\times_M M_l\ar[r]^{e_l\times e_l}&G_0\times_M G_0&  &G_0\\
}
\end{equation}
is the same (up to 2-morphism) as $e_l:  M_l \to G_0$.
Therefore, composing the HS maps in \eqref{2hs} gives an HS
bibundle $J^{-1}_1 (e_l\times e_l (M_l\times_M M_l))$ which is
isomorphic (as an HS bibundle) to $M_l\times_{G_0} G_1$ which
presents the embedding $e_l$. Therefore, one can easily find a
global section
\[ \sigma_l: M_l\to  M_l\times_{G_0} G_1 \cong J_1^{-1}(e_l\times e_l
(M_l\times_M M_l))\subset E_m  \]  defined by $ x\mapsto (x,
1_{e_l(x)})$. Furthermore, we have $J_2 \circ \sigma_l (M_l)
=e_l(M_l)$. Since $G$ is an \'etale groupoid, $E_m$ is an \'etale
principal bundle over $G_0\times_M G_0$. Hence $J_1$ is a local
diffeomorphism.
Therefore, one can choose two open neighborhoods $V_l\subset U_l$
of $M_l$ in $G_0$ such that there exists a unique section
$\sigma'_l$ extending $\sigma_l$ over $(M_l=M_l\times_M M_l
\subset) V_l\times_{M_l} V_l$ in $E_m$ and the image of $J_2\circ
\sigma'_l$ is $U_l$. The restriction of $\sigma'_l$ on $M_l$ is
exactly $\sigma_l$. Since $U_l\rightrightarrows U_l$ acts freely
and transitively fiber-wise on $\sigma'_l(V_l\times_{M_l} V_l)$
from the right,  $\sigma'_l(V_l\times_{M_l} V_l)$ can serve as an
HS bibundle from $V_l\times_{M_l} V_l$ to $U_l$. (Here, we view
manifolds as groupoids.) In fact, it is the same as the morphism
\[m_l:=J_2\circ \sigma'_l: V_l \times_{M_l} V_l \to U_l.\]

 By a similar method, we can define the inverse as
follows. By $\bar{i}\circ \bar{e}_l=\bar{e}_l$, we have the
following commutative diagram:
\[
\begin{array}{ccc}
M_l  & \overset{ \bm=id}{\lra} &  M_l \\
\downarrow &\curvearrowright & \downarrow \\
\cG &\overset{\bar{i}} {\lra} &\cG.
\end{array}
\]
Suppose $(E_i, J_1, J_2)$ is the HS bibundle representing
$\bar{i}$.  Translate the above diagram into groupoids. Then we
have the composition of the following HS morphisms:
\begin{equation} \label{2hs-i}
\xymatrix{M_l\ar[r]\ar[dd]\ar@<-1ex>[dd]& G_1\ar[dd]\ar@<-1ex>[dd]&  &G_1\ar[dd]\ar@<-1ex>[dd]\\
& & E_i\ar[dr]^{J_2}\ar[dl]_{J_1}& &\\
M_l\ar[r]^{e_l}&G_0&  &G_0\\
}
\end{equation}
is the same (up to a 2-morphism) as  $e_l:  M_l \to G_0$.
Therefore, composing the HS maps in \eqref{2hs-i} gives an HS
bibundle $J^{-1}_1 (e_l(M_l))$, which is isomorphic (as an HS
bibundle) to $M_l\times_{G_0} G_1$ which represents the embedding
$e_l$. Therefore, one can easily find a global section
\[ \tau_l: M_l\to  M_l\times_{G_0} G_1 \cong J_1^{-1}(e_l(M_l))\subset
E_i  \] defined by $ x\mapsto (x, 1_{e_l(x)})$. Furthermore, we
have $J_2 \circ \sigma_l (M_l) =e_l(M_l)$. Since $G$ is an \'etale
groupoid, $E_i$ is an \'etale principal bundle over $G_0$. Hence
$J_1$ is a local diffeomorphism. Therefore, one can choose an open
neighborhood of $M_l$ in $G_0$, which we might assume to be $U_l$
as well, such that there exists a unique section $\tau'_l$
extending $\tau_l$ over $(M_l \subset) U_l$ in $E_i$ and the image
of $J_2\circ \tau'_l$ is in $U_l$. The restriction of $\tau'_l$ on
$M_l$ is exactly $\tau_l$. So we can define
\[i_l:=J_2\circ \tau'_l: U_l \to U_l.\]

Since $M$ is a manifold, examining the groupoid picture of maps
$\bbs$ and $\bbt$, one finds that they actually come from two maps
$\bs$ and $\bt$ from $G_0$ to $M$. Hence, we define source and
target maps of $U_l$ as the restriction of $\bs$ and $\bt$ on
$U_l$ and denote them by $\bs_l$ and $\bt_l$ respectively.

The 2-associative diagram of $\bm$ tells us that $m_l\circ
(m_l\times id)$ and $m_l \circ (id\times m_l)$ differ in the
following way: there exists a smooth map from an open subset of
$V_l\times_{M_l} V_l \times_{M_l} V_l$, where both of the above
maps are defined, to $G_1$, such that, \[m_l\circ (m_l\times id)=
m_l \circ (id\times m_l) \cdot \alpha. \] Since the 2-morphism in
the associative diagram restricting to $M$ is $id$, we have
\[\alpha(x, x, x) =1_{e_l(x)}.\] Since $G$ is \'etale and $\alpha$ is
smooth, the image of $\alpha$ is inside the identity section of
$G_1$. Therefore $m_l$ is associative.

It is not hard to verify the other groupoid properties in a similar
way by translating corresponding properties on $\cG$ to $U_l$.
Therefore, $U_l$ with maps defined above is a local Lie groupoid
over $M_l$.
\end{proof}

To prove the global result, we need the following proposition:

\begin{prop}\label{glue}
Given $U_l$ and $U_j$ constructed as above (one can shrink them if
necessary), there exists an isomorphism of local Lie groupoids
$\tilde{\varphi}_{lj}: U_j \to U_l$ extending the isomorphism
$\varphi_{lj}$ in Lemma \ref{eie}. Moreover
the $\tilde{\varphi}_{lj}$'s also satisfy cocycle conditions.
\end{prop}
\begin{proof}
Since we will restrict the discussion to $M_l\cap M_j$, we might
assume that $M_l=M_j$. Then according to Lemma \ref{eie}, there is
a local bisection $g_{lj}$ of $G_1$ such that $e_l \cdot g_{lj}
=e_j$. Extend the bisection $g_{lj}$ to $U_l$ (we denote the
extension still by $g_{lj}$, and shrink $V_k$ and $U_k$ if
necessary for $k=l, j$) so that \[ (V_l \times_{M_l} V_l) \cdot
(g_{lj} \times g_{lj} ) = V_j \times_{M_j} V_j \quad \text{and} \;
U_l \cdot g_{lj} = U_j.\] Notice that since $G_1$ is \'etale, the
source map is a local isomorphism. Therefore, by choosing small
enough neighborhoods of the $M_l$'s, the extension of $g_{lj}$ is
unique. Let $\tilde{\varphi}_{lj}=\cdot g^{-1}_{lj}$. Then it is
naturally an extension of $\varphi_{lj}$. Moreover, by uniqueness
of the extension, the $\tilde{\varphi}_{lj}$'s satisfy cocycle
conditions as the $\varphi_{lj}$'s do.

Now we show that $\tilde{\varphi}_{lj}=\cdot g_{lj}$ is a morphism
of local groupoids. It is not hard to see that $\cdot g_{lj}$
preserves the source, target and identity maps. So we only have
to show that
\[ i_l \cdot g_{lj} =i_j,  \qquad m_l\cdot g_{lj}=m_j. \] For this purpose, we have
to recall the construction of these two maps. $i_l$ is defined as
$J_2 \circ \tau'_l$. Since there is a global section of $J_1$ over
$U_l$ in $E_i$, we have $J_1^{-1}(U_l)\cong U_l\times_{i_l, G_0}
G_1$ as $G$ torsors. Under this isomorphism, we can write
$\tau'_l$ as
\[ \tau'_l (x) = (x, 1_{e_l(x)}).\]
The $G$ action on $U_l\times_{i_l, G_0} G_1$ gives  $(x,
1_{e_l(x)})\cdot g_{lj}=(x, g_{lj})$. Moreover, we have
\[J_2((x, g_{lj}))=J_2(x, 1_{e_j(x)})=\bs_G(g_{lj}), \]
where $\bs_G$ is the source map of $G$. Combining all these, we
have shown that $i_l \cdot g_{lj} =i_j$. The other identity for
multiplications follows in a similar way.
\end{proof}

Recall Theorem \ref{local}:

\noindent {\bf Theorem \ref{local}.}
{\em Given a Weinstein groupoid $\cG$, there is an\footnote{It is canonical up to isomorphism near the
identity section.} associated local Lie groupoid $G_{loc}$ which
has the same Lie algebroid as $\cG$.}

\begin{proof} [Proof of Theorem \ref{local}]
Now it is easy to construct $G_{loc}$ as in the statement of the
theorem. Notice that the set of $\{U_l\}$ with isomorphisms the
$\varphi_{lj}$'s which satisfy cocycle conditions serve as a chart
system. Therefore, gluing them together, we arrive at a global
object $G_{loc}$. Since the $\varphi_{lj}$'s are isomorphisms of local
Lie groupoids, the local groupoid structures also glue together.
Therefore $G_{loc}$ is a local Lie groupoid.

If we choose two different open coverings $\{ M_l\}$ and
$\{M'_l\}$ of $M$ for the same \'etale atlas $G_0$ of $\cG$, we
will arrive at two systems of local groupoids $\{U_l\}$ and
$\{U'_l\}$. Since $\{M_l\}$ and $\{M'_l\}$ are compatible chart
systems for $M$, combining them and using Proposition \ref{glue},
$\{U_l\}$ and $\{U'_l\}$ are compatible chart systems as well.
Therefore they glue into the same global object up to isomorphism
near the identity section.

If we choose two different \'etale atlases $G'_0$ and $G''_0$ of
$\cG$, we can take their refinement $G_0=G'_0\times_{\cG} G''_0$
and we can take a fine enough open covering $\{M_l\}$ so that it
embeds into all three atlases. Since $G_0\to G'_0$ is an \'etale
covering, we can choose the $U_l$'s in $G'_0$ small enough so that
they still embed into $G_0$.  So the groupoid constructed from the presentation $G_0$ with the
covering $U_l$ is the same as the groupoid constructed
from the presentation $G'_0$  with the covering $U_l$'s.   The same
is true for $G''_0$ and $G_0$. Therefore our local groupoid
$G_{loc}$ is canonical.

We will finish the proof of the part containing the
Lie algebroid in the next section.

\end{proof}

\section{Weinstein groupoids and Lie algebroids}

In this section, we define the Lie algebroid of a Weinstein
groupoid $\cG$. One obvious way to do it is to define the Lie algebroid of
$\cG$  as the Lie algebroid of the local Lie groupoid $G_{loc}$.
We now give an equivalent and more direct definition.

\begin{pdef} Given a Weinstein groupoid $\cG$ over $M$, there is
a canonically associated Lie algebroid $A$ over $M$.
\end{pdef}
\begin{proof}
We just have to examine  the second part of the
proof of Theorem \ref{local}. Choose an \'etale groupoid
presentation $G$ of $\cG$ and an open covering $\{M_l\}$'s as in
Lemma \ref{eie}. According to Theorem \ref{local}, we have a local
groupoid $U_l$ and its Lie algebroid $A_l$ over each $M_l$.
Differentiating the $\tilde{\varphi}_{lj}$'s in Proposition
\ref{glue}, we can achieve the algebroid isomorphisms
$T\tilde{\varphi}_{lj}$'s which also satisfy cocycle conditions.
Therefore, using these data, we can glue $A_l$'s into a vector
bundle $A$. Moreover, since the $T\tilde{\varphi}_{lj}$'s are Lie
algebroid isomorphisms, we can also glue the Lie algebroid
structures. Therefore $A$ is a Lie algebroid.

Following the same arguments as in the proof of Theorem
\ref{local}, we can show uniqueness.  For a different
presentation $G'$ and a different open covering $M_l$, we
choose the refinement of these two systems and will arrive at a
Lie algebroid which is glued from a refinement of both systems.
Therefore this Lie algebroid  is isomorphic to both Lie algebroids
constructed from these two systems. Hence the construction is
canonical.

In the language of stacks, what we have just constructed is
actually $\bar{e}^*\ker T\bbs$. As a differential stack, $\ker
T\bbs$, is presented by $\ker T\bs \times_{G_0} G_1
\rightrightarrows \ker T\bs$ for an \'etale presentation $G$ of
$\cG$. Following a similar argument as for the tangent stack, the
above definition for $\ker T\bbs$ is atlas-independent. Its
pull-back to $M$ will be a vector bundle (in the category of
\'etale differentiable stacks) over a manifold. By definition it
is an honest vector bundle over the manifold $M$. Moreover, let
$q_i: M_i \to M$, then $\bar{e}\circ q_i=e_i \circ \pi$, where
$\pi: G_0\to \cG$ is the projection from atlases. Notice that
$TG_0=\pi^* T\cG$, hence $\ker T\bs = \pi^* \ker T\bbs$. So we
have $ q_i^* \bar{e}^*\ker T\bbs =e_i^* \ker T\bs=A_i$ and this
shows that $\bar{e}^*\ker T\bbs$ is $A$.
\end{proof}

Now it is easy to see that the following proposition holds:
\begin{prop}
A Weinstein groupoid $\cG$ has the same Lie algebroid as its
associated local Lie groupoid $G_{loc}$.
\end{prop}

Together with the Weinstein groupoid $\cG(A)$ we have constructed
in Section \ref{w},  we are now ready to finish the proof of
Theorem \ref{lieIII}.


\noindent {\bf Theorem \ref{lieIII}. (Lie's third theorem).}
{\em To each Weinstein groupoid one can associate a Lie algebroid. For
every Lie algebroid $A$, there are naturally two Weinstein
groupoids $\cG(A)$ and $\cH(A)$ with Lie algebroid $A$. }

\begin{proof}[Proof of the second half of Theorem \ref{lieIII}]
We take the \'etale presentation $P$ of $\cG(A)$  and $\cH(A)$ as we
constructed in Section \ref{apath}. Let us first recall how we construct
local groupoids from $\cG(A)$ and $\cH(A)$.

In our case, the HS morphism corresponding to $\bar{m}$ is $$(E:=
\bt_M^{-1}( m(P\times_M P)\cap \bs_M^{-1}(P), m^{-1}\circ \bt_M,
\bs_M).$$ The section $\sigma: M\to E$ is given by $x\mapsto
1_{0_x}$. Therefore if we choose two small enough open
neighborhoods $V\subset U$ of $M$ in $P$, the bibundle
representing the multiplication $m_V$ is a section $\sigma'$ in
$E$ over $V\times_M V$ of the map $m^{-1}\circ \bt_M$.

Since the foliation $\cF$ intersects each transversal slice only
once, we can choose an open neighborhood $O$ of $M$ inside $P_0A$
so that the leaves of the restricted foliation $\cF|_O$ intersect
$U$ only once. We denote the homotopy induced by $\cF|_O$ as
$\sim_O$ and the holonomy induced by $\cF_O$ by $\sim_O^{hol}$.
Then there is a unique element $a \in U$ such that $a\sim_O
a_1\odot a_2$. Since the source map of $\Gamma$ is \'etale, there
exists a unique arrow  $g:
a_1\odot a_2 \curvearrowleft a$ between $ a_1\odot a_2$ and $ a$ in $\Gamma$ near the identity
arrows at $1_{0_x}$'s.

Then we can choose the section $\sigma'$ near $\sigma$ to be
\[\sigma': (a_1, a_2) \mapsto g.\]
So the multiplication $m_V$ on $U$ is
\[ m_V(a_1, a_2) = a(\sim_O a_1\odot a_2). \]

Because the leaves of $\cF$ intersect $U$ only once, $a$ has to be
the unique element in $U$ such that $a\sim^{hol}_O a_1\odot a_2$. It is not hard to verify that both Weinstein groupoids give
the same local Lie groupoid structure on $U$.

Moreover, $U= O/\sim_O$ is exactly the local groupoid constructed
in Section 5 of \cite{cf}, which has Lie algebroid $A$. Therefore,
$\cG(A)$ and $\cH(A)$ have the same local Lie groupoid and the same
Lie algebroid $A$.
\end{proof}

\chapter{Application to integration of Jacobi manifolds}\label{app-ja}
In this chapter, we apply our Weinstein groupoids
to the integration problem of  Jacobi manifolds. To do this, we
will first introduce {\em symplectic Weinstein groupoids}, which
are the generalization of symplectic groupoids in the sense that
any Poisson manifolds (not just integrable ones) can be integrated
into symplectic Weinstein groupoids.

\section{Jacobi manifolds}\label{jam} A {\em Jacobi manifold} is
a smooth manifold $M$ with a bivector field $\Lambda$ and a vector
field $E$ such that
\begin{equation}
[\Lambda,\Lambda]=2E\wedge\Lambda\; \;  \text{and}  \; \;[\Lambda,E]=0,
\label{jacobi}
\end{equation}
where $[\cdot, \cdot]$ is the usual Schouten-Nijenhuis bracket. A
Jacobi structure on $M$ is equivalent to a ``local Lie algebra''
structure on $C^{\infty} (M)$ in the sense of Kirillov
\cite{kirillov}, with the bracket
\[ \{ f, g\} = \sharp \Lambda (df, dg) +fE(g)-gE(f)\quad \forall f, g
\in C^{\infty}(M).\] We call this bracket a Jacobi bracket on
$\ci(M)$. It is a Lie bracket satisfying the following equation
(instead of the Leibniz rule, as for Poisson brackets):
\begin{equation} \label{localbk} \{f_1f_2, g\}=f_1\{ f_2,
g\}+f_2\{ f_1, g\}-f_1f_2\{1, g\},
\end{equation} i.e. it is a first order differential operator in
each of its arguments. If $E=0$, $(M,\Lambda)$ is a Poisson
manifold.

Recall that a \emph{contact manifold}\footnote{ A related concept
is the following: a \emph{contact structure} on the manifold $M$
is a choice of hyperplane $\cH \subset TM$ such that locally
$\cH=\ker(\theta)$ for some 1-form $\theta$ satisfying $\theta
\wedge (d \theta)^n \neq 0$. In this thesis all contact structures
will be co-orientable, so that $\cH$ will be the kernel of some
globally defined contact 1-form $\theta$.}  is a
$2n+1$-dimensional manifold equipped with a 1-form $\theta$ such
that $\theta \wedge (d \theta)^n$ is a volume form. If
$(M,\Lambda,E)$ is a Jacobi manifold such that $\Lambda^n \wedge
E$ is nowhere 0, then M is a contact manifold with the contact
1-form $\theta$ determined by
\[ \iota(\theta)\Lambda=0, \; \; \;\; \;\iota(E) \theta =1 ,\]
where $\iota$ is the contraction between differential forms and
vector fields. On the other hand, given a contact manifold $(M,
\theta)$, let $E$ be the Reeb vector field of $\theta$, i.e. the
unique vector field satisfying
\[ \iota(E) d \theta =0, \; \;\;\; \;\iota(E) \theta =1. \]
Let $\mu$ be the map $TM \to T^*M$ defined by $\mu (X) = -\iota(X)
d \theta$. Then $\mu$ is an isomorphism between $\ker( \theta)$
and $\ker(E)$, and can be extended to their exterior algebras. Let
$\Lambda = \mu^{-1} (d \theta)$. (Note that if $\iota(E) d \theta
=0$, then $d \theta$ can be written as $\alpha \wedge \beta$ and
$\iota(E)\alpha=\iota(E)\beta=0$.) Then $E$ and $\Lambda$ satisfy
\eqref{jacobi}. So a contact manifold is always a Jacobi manifold
\cite{l}. Notice that in this case the map $\sharp \Lambda:T^*M
\lra TM$ given by $\sharp \Lambda(X)=\Lambda(X,\cdot)$ and the map
$\mu$ above are inverses when restricted to $\ker(\theta)$ and
$\ker(E)$.

A {\em locally conformal symplectic manifold} (\emph{l.c.s.
manifold} for short) is a $2n$-dimensional manifold equipped with
a non-degenerate 2-form $\Omega$ and a closed 1-form $\omega$ such
that $d\Omega=\omega \wedge \Omega$. To justify the terminology
notice that locally $\omega=df$ for some function $f$, and that
the local conformal change $\Omega \mapsto e^{-f} \Omega$ produces
a symplectic form. If $(M,\Lambda,E)$ is a Jacobi manifold such
that $\Lambda^n$ is nowhere 0, then M is a l.c.s. manifold: the
2-form $\Omega$ is defined so that the corresponding map $TM\lra
T^*M$ is the negative inverse of $\sharp \Lambda: T^*M \lra TM$,
and the 1-form is given by $\omega=\Omega(E,\cdot)$. Conversely,
if $(\Omega,\omega)$ is an l.c.s. structure on $M$, then defining
$E$ and $\Lambda$ in terms of $\Omega$ and $\omega$ as above,
\eqref{jacobi} will be satisfied.

A Jacobi manifold is always foliated by contact and locally
conformal symplectic (l.c.s.) leaves \cite{marrero}. In fact, like
that of a Poisson manifold, the foliation of a Jacobi manifold is
given by the distribution of the Hamiltonian vector fields
$$X_u:=uE+\sharp \Lambda (du).$$ The leaf through a point will be
a l.c.s. (resp. contact) leaf when $E$ lies (resp. does not lie)
in the image of $\sharp \Lambda$ at that point.

Given a nowhere vanishing smooth function $u $ on a Jacobi
manifold $(M, \Lambda, E)$, a conformal change by $u$ defines a
new Jacobi structure:
\[ \Lambda_u = u \Lambda, \;\;\;\;\; E_u= uE+\sharp \Lambda (d u)
=X_u.\]
 We call two Jacobi structures equivalent if they differ by a conformal change.
 A {\em conformal Jacobi structure} on a manifold
 is an equivalence class of Jacobi structures\footnote{Clearly, a conformal contact manifold is just a manifold
with a coorientable contact structure.}. The relation between the
Jacobi brackets induced by the  $u$-twisted and the original
Jacobi structures is given by
$$\{f,g\}_u=u^{-1}\{uf,ug\}$$
The relation between the Hamiltonian vector fields is given by
$$X^u_f=X_{u\cdot f}.$$
 A smooth map $\phi$ between Jacobi
manifolds $(M_1, \Lambda_1, E_1)$ and $(M_2, \Lambda_2, E_2)$ is a
{\em Jacobi morphism} if
\[ \phi_* \Lambda_1= \Lambda_2, \;\;\;\;\; \phi_* E_1 =E_2, \]
or equivalently if $\phi_*(X_{\phi^*f})=X_f$ for all functions $f$
on $M_2$.
 Given $u\in \ci(M_1)$, a {\em u-
 conformal Jacobi
morphism} from a Jacobi manifold $(M_1, \Lambda_1, E_1)$ to $(M_2,
\Lambda_2, E_2)$ is a Jacobi morphism from  $(M_1, (\Lambda_1)_u,
(E_1)_u)$ to $(M_2, \Lambda_2, E_2)$.

The Lie algebroid  associated to a Jacobi manifold $(M, \Lambda,
E)$ is $T^*M\oplus_M\R$ \cite{ks}, with anchor $\rho$ and bracket
$\br$ on the space of sections $\Omega^1(M) \times \ci(M)$ defined
by
\begin{equation}
\begin{split}
\rho: \Omega^1(M) \times \ci(M) \lra \chi(M) \\
\rho(\omega_1, \omega_0) = \sharp\Lambda (\omega_1) + \omega_0 E ,
\end{split}
\end{equation}
and
\begin{equation}
\begin{split}
[(\omega_1,\omega_0),(\eta_1,\eta_0)]  = & \big(
\call_{\sharp\Lambda(\omega_1)}\eta_1-\call_{\sharp\Lambda(\eta_1)}\omega_1-\di(\Lambda(\omega_1,
\eta_1))  \\
&+\omega_0\call_E\eta_1-\eta_0\call_E\omega_1-i(E)\omega_1 \wedge \eta_1  ,\\
&\Lambda(\eta_1,\omega_1)+\sharp\Lambda(\omega_1)(\eta_0)-\sharp\Lambda(\eta_1)(\omega_0)\\
&+\omega_0 E(\eta_0)-\eta_0 E (\omega_0) \big),
\end{split}
\end{equation} where $ \sharp \Lambda : T^*M \lra TM $ is the
bundle map defined by $\langle\sharp \Lambda (\omega),
\eta\rangle= \Lambda (\omega, \eta)$ for all $\omega, \eta \in
\Omega^1(M) $.

At first sight, this seems rather complicated, but the bracket is
determined by the following two conditions:
\begin{align}
& [(\omega_1, 0), (\eta_1, 0)]= ([\omega_1, \eta_1]_{\Lambda}, 0)- (i_E(\omega_1\wedge\eta_1), \Lambda(\omega_1, \eta_1))\nonumber\\
& [(0, 1), (\omega_1, 0)]= (\mathcal{L}_{E}(\omega_1), 0),
\nonumber
\end{align}
where the bracket $[\cdot, \cdot ]_\Lambda$ is analogous to that
for Poisson manifolds,
\[
[\omega_1, \eta_1]_{\Lambda}=
\mathcal{L}_{\sharp\Lambda\omega_1}\eta_1-
\mathcal{L}_{\sharp\Lambda\eta_1}\omega_1- d(\Lambda (\omega_1,
\eta_1)).
\]

\section{Homogeneity and Poissonization} \label{homogeneous}
Recall that a {\em homogeneous symplectic manifold} $(M, \omega,
Z)$ is a symplectic manifold with a 2-form $\omega$ and a vector
field $Z$ satisfying $ \call_Z \omega =\omega$.

Given a contact manifold $(M, \theta)$, we can construct its
associated  homogeneous symplectic manifold $(M \times \R ,
\omega, \frac{\partial}{\partial s})$, where $\omega=\di ( e^s
\pi^* \theta) $, $\pi$ is the projection $M \times \R \lra M$, and
$s$ is the coordinate on $\R$. On the other hand, if we have a
homogeneous symplectic manifold of the special form $(M\times \R ,
\omega, \frac{\partial}{\partial s})$ such that
$\call_{\frac{\partial}{\partial s}} \omega = \omega$ (i.e. $Z$
generates a free action of $\R$), then $M$ has a contact 1-form
$\theta = i_0^* \iota (\frac{\partial}{\partial s}) \omega$ and we
will have $\omega = \di (e^s \pi^* \theta)$. Here $i_0$ is the
embedding of $M$ as the 0-section of $M \times \R$. A similar
construction exists for Jacobi and Poisson manifold:

We call a Poisson manifold $(P, \tilde{\Lambda})$ {\em
homogeneous}  if there is a vector field $Z$ such that
\[ \tilde{\Lambda}=-[Z,  \tilde{\Lambda}].\]
If $Z$ never vanishes on $P$, and if $\pi$: $P \lra P/Z=:M$ is a
submersion\footnote{By $P/Z$, we mean the quotient of $P$ by the
flow generated by $Z$.}, then $M$ has a unique Jacobi structure
induced by $P$ \cite{jacobi1}.

Given a Jacobi manifold $(M, \Lambda, E)$, we can construct a
homogeneous Poisson structure on $M \times \R$, by defining
\[ \tilde{\Lambda}= e^{-s} \left(i_*\Lambda+\frac{\partial}{\partial s}\wedge i_* E \right), \]
where $s$ is the coordinate on $\R$ and $i_*: TM \lra T(M
\times\R)$ is the inclusion. Here we view $TM$ as a bundle over $M
\times \R$ using the pullback by the  projection $\pi: M \times \R
\to M$. The bivector field $\tilde{\Lambda}$ is Poisson, i.e. $
[\tilde{\Lambda}, \tilde{\Lambda}]=0$, precisely when
\eqref{jacobi} holds for $\Lambda $ and $ E$, and it is easy to
check that $\tilde{\Lambda}=-\left[\frac{\partial}{\partial
s},\tilde{\Lambda}\right]$. So $(M \times \R, \tilde{\Lambda})$ is
a homogeneous Poisson manifold with vector field
$\frac{\partial}{\partial s}$.

Conversely, given  a homogeneous Poisson manifold $(M \times \R,
\tilde{\Lambda})$ satisfying
\[ \tilde{\Lambda}= -\left[\frac{\partial}{\partial s}, \tilde{\Lambda}\right], \] $M$ inherits a unique Jacobi
structure
\[ \Lambda= \pi_* (e^s \tilde{\Lambda}) , \;\;\;\;\; E= \sharp (e^s \tilde{\Lambda})(\di s). \]
The condition \eqref{jacobi} for $(M, \Lambda, E)$ to be a Jacobi
manifold is exactly equivalent to $ [\tilde{\Lambda},
\tilde{\Lambda}]=0$. Thus we  have the following lemma:

\begin{lemma}
There is a one-to-one correspondence between Jacobi structures on
$M$ and homogeneous Poisson structures on $M\times\R$ with
homogeneous vector field $\dds$.

Moreover, restricting this procedure to contact manifolds, there
is a one-to-one correspondence between contact structures on $M$
and homogeneous symplectic structures on $M\times \R$ with
homogeneous vector field $\dds$.
\end{lemma}

\begin{remark}
Explicitly, the relation between $\tilde{\Lambda}$ and ($\Lambda$,
$E$) is the following:
\begin{equation} \label{tl}
\tilde{\Lambda}( \omega_1 + \omega_0 \di s , \eta_1 + \eta_0 \di
s) =e^{-s} ( \Lambda (\omega_1, \eta_1) + \omega_0 E (\eta_1)
-\eta_0 E(\omega_1)).
\end{equation}
Here, at every point of $M \times \R$, we view $\omega_1$,
$\eta_1$ as 1-forms on $M$ and $\omega_0$, $\eta_0$ as functions
on $M$ after fixing $s$.
\end{remark}



\section{Symplectic (resp. contact) Weinstein groupoids} The
associated Lie algebroid of  a Poisson manifold $P$ is $T^*P$.
Therefore we can associate to $P$ two {\em symplectic Weinstein
groupoids} $\Gamma_s^m(P)$ and $\Gamma_s^h(P)$, which are
$\cG(T^*P)$ and $\cH(T^*P)$ respectively. We will define what a
symplectic Weinstein groupoid is. They are called symplectic
Weinstein groupoids because, when $P$ is integrable,
$\Gamma_s^h(P)$ which is the same as the orbit space of
$\Gamma_s^m(P)$ is the source-simply connected symplectic groupoid
integrating $P$. Similarly, to a Jacobi manifold $M$ we associate
two  Weinstein groupoids $\Gamma_c^m(M)$ and $\Gamma^h_c(M)$
defined as  $\cG(T^*M \oplus \R)$ and $\cH(T^*M\oplus \R)$
respectively. In the next section, we will try to make them into
contact Weinstein groupoids.

Let $\cX$ be a stack over $\cC$. Then a {\bf sheaf of differential
k-forms} $\cF^k$ is defined as follows: for every $x\in \cX$ over
$U\in \cC$, $\cF(x)=\Omega^k(U)$. It is a contravariant functor:
for every arrow $y\to x$ over $f:V\to U$, there is a map
$\Omega^k(U)\to\Omega^k(V)$ defined by pull back via $f$. It is
moreover a sheaf over $\cX$. As for the definition of sheaves over
stacks and the proof, we refer to \cite{bx} since we will not use
this later in this thesis. Then a {\bf differential $k$-form}
$\omega$ on $\cX$ a map that associates an element $x\in\cX$ over
$U$ a section $\omega(x) \in \Omega^k(U)$ such that the following
compatibility condition holds: if there is an arrow $y\to x$ over
$f: V\to U$, then $\omega(x)$ is the pull back of $\omega(y)$ via
$\cF^k(f)$. Notice that according to the above definition, the
0-forms on $\cX$ are simply morphisms of stacks from $\cX$ to $\R$.

\begin{lemma} \label{etale-form}
When $\cX$ is a \'etale differentiable stack, let $G$ be an \'etale groupoid
presentation. Then there is a 1-1 correspondence between
differential forms on $\cX$ and $G$ invariant forms on $G_0$.
\end{lemma}
\begin{proof}
A $G$ invariant $k$-form $\omega$ on $G_0$ defines a differential
form on $\cX$ in the following way: given a right $G$-principal
bundle $\pi: P\to U$ with moment map $J: P\to G_0$, the pull back
form $J^*\omega$ is $G$ invariant on $P$, therefore it induces a
$k$-form $\pi_* J^* \omega$ on $U$ and this is what $P$ associates
to via $\omega$. Notice here we use the fact that $\pi$ is \'etale
to show that $G$-invariant form is a basic form. On the other
hand, given any $k$-form $\omega$ on $\cX$, consider $\bs:G_1\to
G_0$ as a right $G$-principal bundle. Then $\omega(G_1)$ is a
$k$-form on $G_0$. Notice that $ g\cdot : G_1\to G_1$ is a
morphism of $G$-principal bundles. Using the compatibility
condition of $\omega$, we can see that $\omega(G_1)$ is
$G$-invariant.
\end{proof}

Use this correspondence in the special case of \'etale
differentiable stacks, we can make the following definitions:

\begin{defi}[symplectic (resp. contact) forms on \'etale differentiable stacks]
A symplectic (resp. contact) form on an \'etale differentiable
stack $\cX$ is a $G$ invariant symplectic (resp. contact) form on
$G_0$, where $G$ is an \'etale presentation of $\cX$.
\end{defi}

\begin{pdef}[pull back of forms on stacks]
Let $\phi:\cY\to\cX$ be a map between stacks and $\omega$ a form
on $\cX$. Then $\phi^*\omega$ is a form on $\cY$ defined by
associating $y\in \cY$ to $\omega(\phi(y))$.
\end{pdef}
\begin{proof}
It is not hard to verify the compatibility condition for
$\phi^*\omega$.
\end{proof}
\begin{remark}
Using Lemma \ref{etale-form}, it is not hard to see that pull
backs of forms on \'etale differentiable stacks corresponds to the
ordinary pull backs on the \'etale atlases.
\end{remark}

\begin{defi}[symplectic Weinstein groupoids] A Weinstein groupoid
$\cG$ over a manifold $M$ is a symplectic Weinstein groupoid if
there is a symplectic form $\omega$ on $\cG$ satisfying the
following multiplicative condition:
\[ \bm^* \omega = pr_1^* \omega + pr_2^* \omega, \]on $\cG\times_{\bbs, M,
\bbt}\cG$, where $pr_i$ is the projection onto the $i$-th factor.
\end{defi}

\begin{defi}[contact Weinstein groupoids] A  Weinstein
groupoid $\cG$ over a manifold $M$ is a contact Weinstein groupoid
if there are a contact 1-form $\theta$ and a function $f$, such
that the following twisted multiplicative condition hold on
$\cG\times_{\bbs, M, \bbt}\cG$:
\[ \bm^* \theta =pr_2^* f \cdot pr_1^*\theta +pr_2^*\theta.\]
\end{defi}
\begin{remark}
When the Weinstein groupoid $\cG$ is a Lie groupoid, the above
definitions coincide with the definitions of sympletic groupoids
and contact groupoids \cite{ks} respectively.
\end{remark}

\begin{thm} Let $N$ be a Poisson manifold. Then $\Gamma_s^m(N)$
and $\Gamma_s^h(N)$ are symplectic Weinstein groupoids over $N$.
\end{thm}
\begin{proof}
We prove it for $\Gamma_s^m(N)$ and the proof for $\Gamma_s^h(N)$
is similar. Let $\omega_c$ be the canonical symplectic form on
$T^*M$. Then according to \cite{cafe}, $\omega_c$ induces a
symplectic form on the path space $P T^*M$. This symplectic form
restricted to $P_a T^*M$ has kernel exactly the tangent space of
the foliaction $\cF$ and invariant along the foliation. Consider
the \'etale presentation $\Gamma\rightrightarrows P$ of
$\Gamma_s^m(N)$. $P$ is  the transversal of the foliation $\cF$,
hence the restricted form is a $\Gamma$-invariant symplectic form.
This form induces a symplectic form $\omega$ on $\Gamma_s^m(N)$.
The multiplicativity of $\omega$ follows from the additivity of
the integrals after examining the definition of $\omega_c$.
\end{proof}

\begin{remark} We can not extend the proof to the contact case.
The obstruction is that the contact 1-form $\theta_c + ds$ on
$T^*M \oplus \R$ can not induce a contact form on the path space
$P T^*M \oplus \R$, where $M$ is a Jacobi manifold. That is why we
still need the  following results about the relations between
$\Gamma_s^m(M\times \R)$ and $\Gamma_s^m(M)$ (resp.
$\Gamma_s^h(M\times \R)$ and $\Gamma_s^h(M)$).
\end{remark}

\section{An integration theorem}
\begin{defi} A multiplicative function $r$ on a Weinstein groupoid
$\cG$ is a smooth function $r: \cG\to \R$ such that $$r\circ
\bar{m} = r\circ pr_1 + r \circ pr_2,$$ where $pr_i$ is the
$i$-th projection $\cG\times_{\bs, \bt} \cG \to \cG$.
\end{defi}

\begin{pdef} Given a multiplicative function $r$ on a
Weinstein groupoid $\cG$, one can form a new Weinstein groupoid
$\cG\times_r \R$ over $M \times \R$ which is $\cG\times \R$ as a
stack and has the following groupoid structure:
\begin{alignat*}{3}
\bbs & = & (\bbs\circ pr_{\cG}, \; pr_{\R}), \qquad \bbt &=&
(\bbt\circ pr_{\cG}, \; pr_{\R} - r\circ pr_{\R}), \\ \bar{e} &=&
(\bar{e}, \; id), \qquad \bar{i} & = & (\bar{i}, \; id-r\circ pr_1), \\
\bar{m} & = & (\bar{m} \circ (pr^1_\cG \times pr^2_\cG), \;
pr^2_\R)
\end{alignat*} where $pr^i_\cG$ (or $pr^i_\R$) is the projection
from  $(\cG \times_r \R) \times_{M\times \R}(\cG \times_r \R)$
onto the $i$-th copy of $\cG$ (or $\R$).
\end{pdef}
\begin{proof}
The multiplicativity of $r$ and the 2-associativity of the
multiplication on $\cG$ imply  property (3) in the definition of
Weinstein groupoid. It is routine to check that  the other properties
also hold.
\end{proof}

\begin{thm} \label{wgpd-c-s}Let $M$ be a Jacobi manifold, $M\times \R$ its
Poissonization. Then
\begin{enumerate}
\item[i)] there are two well-defined multiplicative functions $r_m$ and $r_h$ on
$\Gamma_c^m(M)$ and $\Gamma_c^h(M)$ respectively, such that
\[ \Gamma_s^m(M\times\R)\cong \Gamma_c^m(M)\times_{r_m} \R, \quad \Gamma_s^h(M\times\R)\cong
\Gamma_c^h(M)\times_{r_h} \R,\]as Weinstein groupoids;
\item[ii)] $M$ is integrable as a Jacobi manifold iff $M\times
\R$ is integrable as a Poisson manifold.
\end{enumerate}
\end{thm}

Before proving this, we introduce a useful lemma.

\begin{lemma}\label{basic correspondence} Let $\ta(t) $ be an $A$-path  over
$\tilde{\gamma}$ in $T^*(M \times \R)$. Let $s$ be the coordinate
on $\R$. We can decompose $\ta(t)=\ta_1(t)+\ta_0(t) ds$ and
$\tilde{\gamma}=(\gamma_1, \gamma_0)$. Here $\ta_1(t)$ is the part
that does not contain $\di s$; $\gamma_1$ and $\gamma_0$ are paths
in $M$ and $\R$ respectively. Let
\begin{equation}
a_i(t) =e^{-\gamma_0(t)} \ta_i (t), \;\;\;\;\; (i=0,1).
\end{equation}
Then $(a_1(t), a_0(t))$ will be an $A$-path over $\gamma_1(t)$ in
$T^*M \oplus_M \R$.

Conversely, if $(a_1(t), a_0(t)) $ is an $A$-path over
$\gamma_1(t)$ in $T^*M \oplus_M \R$ and $s$ is any real number,
let
\begin{equation}
\begin{split}
\ta_i(t) &=e^{\gamma_0(t)} a_i(t),\\
\gamma_0(t) &= -\int_0^t \iota(E)a_1(t)\di t + s.
\end{split}
\end{equation}
Then $\ta_1(t)+\ta_0(t) \di s $ will be an $A$-path over
$(\gamma_1(t), \gamma_0(t))$ in $T^*(M \times \R)$. In other
words, there is a 1-1 correspondence:
\[
\begin{array}{ccc}
P_a(T^*(M \times \R)) & \os{\phi_a}{\lra} &  P_a(T^*M \oplus_M \R) \times \R \\
 \ta_1+ \ta_0 \di s   & \longmapsto          &  ((a_1, a_0),
 \gamma_{0}(0)).
\end{array}
\]
Furthermore, this correspondence extends to the level of
homotopies of $A$-paths
\[ \ta(\epsilon, t) \mapsto ( a(\epsilon, t), \gamma_0(\epsilon,
0)=\gamma_0(0, 0)). \]
\end{lemma}
\begin{proof}
$\ta(t)$ being an $A$-path is equivalent to saying that
\begin{equation} \label{1}
\rho(\ta) = \frac{ \di} {\di t} \tilde{\gamma} (t) =
\frac{\di}{\di t} \gamma_1 + \frac{\di}{\di t} \gamma_0
\frac{\partial}{\partial s} .
\end{equation}
By \eqref{tl},
\begin{equation}
\begin{split}
\tilde{\rho} ( \omega_1 + \omega_0 \di s) &= \sharp
\tilde{\Lambda} (\omega_1
+ \omega_0 \di s) \\
&= e^{-s} \left( \rho( \omega_1, \omega_0) - E(\omega_1)
\frac{\partial}{\partial s} \right).
\end{split}
\end{equation}
So \eqref{1} is equivalent to
\begin{equation}
e^{-\gamma_0} \left( \rho(\ta_1, \ta_0)- \ta_1(E)
\frac{\partial}{\partial s} \right) = \frac{\di}{\di t} \gamma_1 +
\left( \frac{\di}{\di t} \gamma_0 \right) \frac{\partial}{\partial
s}.
\end{equation}
Given that $\ta_i(t)=e^{-\gamma_0(t)} a_i(t)$, this is equivalent to \\
\[
\begin{cases}
    e^{\gamma_0(t)} = e^{\gamma_0(0)}- \int_0^t \iota(E)\ta_1(t) \di t   &
    \text{or} \; \;
    \gamma_0(t) =  \gamma_0(0) - \int_0^t \iota(E) a_1(t) \di t, \\
    \rho(a_1, a_0) = \frac{\di}{\di t} \gamma_1 ,
\end{cases}
\]
which shows that $(a_1(t), a_0(t) ) $ is an $A$-path over
$\gamma_1$ in $T^*M \oplus_M \R$.

On the other hand, if $(a_1(t), a_0(t))$ is an $A$-path, reversing
the above reasoning shows that $\ta_1(t)+ \ta_0(t) \di s $ will
be an $A$-path too.

We use Proposition-Definition \ref{homotopy}, to see that the 1-1
correspondence preserves the equivalence classes, let
$\ta(\epsilon , t) $ be a family of $A$-paths such that the
solution of
\begin{equation} \label{2}
\partial_t \tilde{b} (\epsilon, t) -\partial_{\epsilon} \ta
(\epsilon, t) = T_{\tilde{\nabla}} (\ta, \tilde{b}) , \;\;\;\;\;
\tilde{b} (\epsilon, 0) =0,
\end{equation}
satisfies $\tilde{b} (\epsilon , 1) =0$. Here $\tilde{\nabla}$ is
the product of a connection $\nabla$ on $M$ and the trivial
connection on $\R$. Straightforward calculation shows that
\begin{equation}
\begin{split}
e^s T_{\tilde{\nabla}} (\ta, \tilde{b}) & =  \left(T_{\nabla}
(\ta,
\tilde{b} )\right)_1+ \ta_1(E) \tb_1 -\tb_1 (E) \ta_1 \\
&\quad + \left(-\Lambda ( \ta_1, \tb_1) -\ta_0 \tb_1(E) + \tb_0
\ta_1(E)\right) \di s.
\end{split}
\end{equation}
Here $\tb = \tb_1 + \tb_0 \di s$ and $T_{\nabla} (\ta, \tilde{b}
)=\left(T_{\nabla} (\ta, \tilde{b} )\right)_1 + \left(T_{\nabla}
(\ta, \tilde{b} )\right)_0 \di s$. So \eqref{2} is equivalent to
\begin{equation} \label{3}
  \begin{cases}
    \partial_t \tb_1 -\partial_{\epsilon} \ta_1 =  e^{-s} \left( \left(T_{\nabla} (\ta,
\tilde{b} )\right)_1+ \ta_1(E) \tb_1 -\tb_1 (E) \ta_1\right), \\
    \partial_t \tb_0 -\partial_{\epsilon} \ta_0 = e^{-s}\left(-\Lambda (
\ta_1, \tb_1) -\ta_0 \tb_1(E) + \tb_0 \ta_1(E)\right).
  \end{cases}
\end{equation}
On the other hand, $a(0, t) \sim a(1, t) $ iff there is a family of $A$-paths $a(\epsilon, t ) $ such that
the solution of
\begin{equation} \label{4}
\partial_t b -\partial_\epsilon a = T_{\nabla} (a, b),\;\;\;\;\; b(\epsilon
,0)=0
\end{equation}
satisfies $ b(\epsilon, 1)=0$.

Let $b_i( \epsilon, t) = \tb_i(\epsilon, t) e^{-\gamma_0(\epsilon,
t)}$ and $b=(b_1, b_0)$. Then \eqref{2} implies \eqref{4}, and
$b(\epsilon, 1) =0$.

Let \[\tb_i (\epsilon , t)= b_i(\epsilon, t) e^{\gamma_0(\epsilon,
t) } , \] where $\gamma_0(\epsilon, t)  = -\int_0^ \epsilon
\iota(E) b_1(\epsilon, t) \di \epsilon - \iota(E) \int_0^t a_1(0,
t) \di t + \gamma_0(0, 0)$. Then \eqref{4} implies \eqref{2} and
$\tb_i(\epsilon, 1)=0$.

So $\ta(0, t) \sim \ta(1, t)$ if and only if $a(0, t) \sim a(1,
t). $
\end{proof}

Now we are ready to prove the theorem.

\begin{proof} [Proof of Theorem \ref{wgpd-c-s}]
We adapt the notation of Lemma \ref{basic correspondence}. Write
an $A_0$-path $a$ of $T^*M \oplus \R$ as $(a_1, a_0)$.  Let
$r(a)=-\int_0^1 \iota(E)a_1(t) \di t$ be a function on $P_0 (T^*M
\oplus \R)$. From the calculation in Lemma \ref{basic
correspondence}, $-\int_0^1 \iota(E) a_1(t) \di t = \gamma_0(1)
-\gamma_0(0)$. Since the base paths of equivalent $A_0$-paths all
have the same end points, $-\int_0^1\iota(E)a_1(t) \di t $ does
not depend on the choice of $(a_1, a_0)$ within an equivalence class.
Therefore $r(a)$ is invariant under the action of the  monodromy
groupoid of $P_0 (T^*M \oplus \R)$. By Lemma \ref{im}, we obtain a
smooth map $r_m$ on $\Gamma_c^m(M)$. Moreover, by the definition
of $r$, \[r\circ m- r\circ pr_1-r\circ pr_2 =0, \quad \text{on}\;
P_0(T^*M \oplus \R)\times_{M} P_0(T^*M \oplus \R).  \] This means
that
\[r_m \circ \bm -r_m\circ pr_1-r_m\circ pr_2 =0, \]
on the level of stacks since this function composed with $\pi:
P_0(T^*M \oplus \R)\times_{M} P_0(T^*M \oplus \R) \to
\Gamma_c^m(M)\times_M \Gamma_c^m(M)$ is $r\circ m- r\circ
pr_1-r\circ pr_2 =0$. Hence $r_m$ is multiplicative.

Lemma \ref{basic correspondence} gives the correspondence between
$P_0(T^*M \oplus \R)\times \R$ with foliation $\cF\times \R$ and
$P_0(T^*(M\times \R)$ with foliation $\cF$, where $\cF$ is the
foliation we defined on any $A_0$-path spaces in Section
\ref{apath}. Moreover, the condition $\gamma_0(\epsilon,
0)=\gamma_0(0, 0)$ tells us that the monodromy groupoid of the
foliation $\cF\times \R$ on $P_0(T^*M \oplus \R)\times \R$ splits
into a product of two groupoids: $Mon_\cF (P_0(T^*M \oplus \R))$
and $\R \rightrightarrows \R$. On the
level of stacks, we have
\[ \Gamma_c^m (M)\times \R \cong \Gamma_s^m(M\times \R).\]
The groupoid structure of the right hand side carries over to the
left hand side and we obtain  exactly $\Gamma_c^m(M)\times_r\R$ by
the calculation below,
\[ \bbs\big([\ta_1 + \ta_0 \di s]\big) = \big( \gamma_1(0), \gamma_0(0) \big)
=\big(\bbs ( [(a_1, a_0)]),\gamma_0(0)\big), \] and
\[
\bbt \big([\ta_1 + \ta_0 \di s]\big) = \big( \gamma_1(1),
\gamma_0(1)\big) =\left(\bbt ( [(a_1, a_0)]),\gamma_0(0)-\int_0^1
\iota(E)a_1(t) \di t \right).
\]

The same argument applies to $\Gamma_c^h(M)$. By Theorem \ref{integ},
$M$ is integrable iff $\Gamma_c^h(M)$ is representable, so (iii)
follows easily.
\end{proof}

\section{Contact groupoids and Jacobi manifolds} \label{jc}
Due to lack of knowledge on differentiable stacks, we will only
present the construction of contact groupoids in the integrable
case and leave the rest as a conjucture, which will be straight
forward to carry out when people know more about differentiable
stacks. In this section, we assume that  $\Gamma_c^h(M)$ is
representable (i.e. $M$ is integrable)  explore the
geometric structures on $\Gamma_c^h(M)$.

Let us first recall:
\begin{defi}
A {\em contact groupoid} \cite{ks} is a Lie groupoid $\Gamma
\underset{ \bt} {\overset{\bs} {\rightrightarrows}} \Gamma_0 $
equipped with a contact 1-form $\theta$ and a smooth function $f$,
such that on the space of multipliable pairs $\Gamma_2$ we have
\begin{equation} \label{contact gpd}
m^* \theta = pr_2^* f \cdot pr_1^* \theta + pr^*_2 \theta ,
\end{equation}
where $pr_j$ is the projection from  $\Gamma_2 \subset \Gamma
\times \Gamma$ onto the $j$-th factor.
\end{defi}
\begin{remark}
Contact groupoids can also defined without refering to a 1-form but
just contact structures. We refer the readers to \cite{zz} for a
detailed presentation of the
definition and the relation between the two definitions.
\end{remark}
\begin{remark}
If $\Gamma$ is a Lie groupoid and $\theta$ is a 1-form such that
$(\Gamma, \theta, f)$ is a contact groupoid for some function $f$,
then $f$ is {\em unique}.  If  $f_1$ and $f_2$ are two such
functions, then by \eqref{contact gpd},
\[ pr_2^* (f_1 -f_2) pr_1^* \theta=0. \]
So $(f_1(y) -f_2(y)) \theta(x)=0$ if $\bt(y)= \bs(x)$. Since
$\theta$ is nowhere 0, $f_1(y)=f_2(y)$ for all $ y \in \Gamma$.
\end{remark}

It is known that given a contact groupoid $( \Gamma \underset{
\bt}{ \overset{\bs} {\rightrightarrows}} M, \theta, f)$, the
manifold $M$ of units is a Legendrian submanifold of $\Gamma$ and
there is a unique Jacobi structure on $M$ so that $\bs$ is a
Jacobi morphism. Then $\bt$ is a $-f$ conformal Jacobi morphism
and the Lie algebroid of $\Gamma$ is isomorphic to $T^*M \oplus_M
\R$---the Lie algebroid associated to $M$ \cite{ks}. In this case,
we call $\Gamma$ the {\em contact groupoid of the Jacobi manifold}
$M$.

\begin{thm} \label{congpd} For an integrable Jacobi manifold $M$, $\Gamma_c^h(M)$ is the
unique source-simply connected contact groupoid over $M$ such that
$\bs$ is a Jacobi map.
\end{thm}

From now on in this section, we assume the integrability of $M$.
Let us first prove some propositions and lemmas.

\begin{prop}\label{prop2}
The groupoid $\Gamma_c^h(M) \times \R \underset{ \bbt}
{\overset{\bbs} {\rightrightarrows} }M \times \R$, with the
symplectic form $\omega$ induced by the isomorphism
$\Gamma_s^h(M\times\R) \cong \Gamma_c^h(M) \times \R$, is a
homogeneous symplectic groupoid with vector field
$-\frac{\partial}{\partial s}$.
\end{prop}

\begin{proof}
We only have to show that $\call_{-\frac{\partial}{\partial s}}
\omega =\omega$. For any $u \in \R$, let
$$\phi_u: M\times \R \lra M\times \R,\;\;\;\; (p, s) \mapsto (p, s+u),$$
$$\Phi_u: \Gamma_c^h(M) \times \R \lra \Gamma_c^h(M) \times \R, \;\;\;\;
\left([(a_1, a_0)], s\right) \mapsto \left([(a_1, a_0)],
s+u\right).$$ Then $\Phi_u$ is an automorphism of the groupoid
$\Gamma_c^h(M) \times \R$ and we have the  commutative diagram
\[
\begin{array} {ccc}
 \Gamma_c^h(M) \times \R  & \os{ \Phi_u }{\lra} & \Gamma_c^h(M)
\times \R \\
 \bar{\bs}\! \downarrow \! \downarrow \! \bar{\bt} & &
 \bar{\bs} \! \downarrow \!\downarrow \bar{\bt} \\
M \times \R & \os{\phi_u}{\lra} & M\times \R.
\end{array}
\]

Let $(\Gamma\times\R)$ be the source-simply connected groupoid of
$(M\times\R, \tilde{\Lambda})$. Then $(\Gamma_c^h(M) \times \R,
e^{-u} \omega) $ is the source-simply connected groupoid of
$(M\times\R, e^u\tilde{\Lambda})$. Since $\phi_u$ is an
isomorphism between the Poisson manifolds $(M\times\R,
\tilde{\Lambda})$ and $(M\times\R, e^u\tilde{\Lambda})$, by the
uniqueness of the symplectic groupoid, $\Phi_u$ must be an
isomorphism of symplectic groupoids, i.e. the following diagram
commutes:
\[
\begin{array} {ccc}
( \Gamma_c^h(M) \times \R , \omega) & \os{ \Phi_u }{\lra} &
(\Gamma_c^h(M)
\times \R, e^{-u} \omega) \\
 \bar{\bs}\! \downarrow \! \downarrow \! \bar{\bt} & &
 \bar{\bs} \! \downarrow \!\downarrow \bar{\bt}\\
(M \times \R, \Lambda) & \os{\phi_u}{\lra} & (M\times \R, e^u
\Lambda).
\end{array}
\]
Since $\Phi_u$ is the flow generated by the vector field
$\frac{\partial}{\partial s}$, and $\phi_u^* \omega =e^{-u}
\omega$, we immediately have $\call_{-\frac{\partial}{\partial s}}
\omega =\omega$.
\end{proof}
\begin{remark} By the explanation in Section \ref{homogeneous}, there is a contact 1-form
$\theta= -i_0^* \iota( \frac{\partial}{\partial s} ) \omega$ on
$\Gamma_c^h(M)$ and $\omega$ has the form  $\omega= \di (e^{-s}
\pi^* \theta)$.
\end{remark}

\begin{prop} \label{prop3}
The groupoid $(\Gamma_c^h(M) \times \R, \omega) $, with the
induced groupoid structure given in Theorem \ref{wgpd-c-s}, is a
symplectic groupoid over $M \times \R$ with $\omega= \di (e^{-s}
\pi^* \theta)$ iff $(\Gamma_c^h(M), \theta, e^{-c} )$ is a contact
groupoid over $M$, where $\pi: \Gamma_c^h(M) \times \R \lra
\Gamma_c^h(M)$ is the projection.
\end{prop}
\begin{proof} Before, we didn't distinguish  structure
  maps on Weinstein groupoids. For clarity, we use
  the following notations only for the proof of this proposition:
for $\Gamma_c^h(M) \times \R$, suppose that $$ \tilde{m}: \;
\tilde{\Gamma}_2 := \{ ( (x', s'), (x, s)):
\bar{\bt}((x',s'))=\bar{\bs}((x,s)) \} \lra \Gamma_c^h(M) \times
\R$$ is the multiplication and $\tilde{pr}_1$, $\tilde{pr}_2$:
$\tilde{\Gamma}_2 \lra \Gamma_c^h(M) \times \R$ are the
projections onto the first and second components respectively.

Similarly, for $\Gamma( M)$, suppose that $m$: $\Gamma_2 := \{ (
x', x): \bbt(x')= \bbs(x) \} \lra \Gamma_c^h(M)$ is the
multiplication and $pr_1$, $pr_2$: $\Gamma_2 \lra \Gamma_c^h(M)$
are the projections onto the first and second components
respectively.

Given $\omega =\di (e^{-s} \pi^* \theta)= e^{-s} \di (\pi^*
\theta) - e^{-s} \di s \wedge \pi^* \theta $, we only have to
establish the equivalence between the two equations:
\begin{equation}\label{*}
\tilde{m}^* \omega =\tilde{pr}_1^* \omega + \tilde{pr}_2^* \omega
\end{equation}
and
\begin{equation}\label{**} m^*\theta=\pr_2^* (e^{-c}) \pr_1^*
\theta + \pr_2^* \theta.
\end{equation}

Note that $$\tilde{\Gamma}_2 = \{ ((x',s'), (x, s)): (x', x) \in
\Gamma_2 , s + c(x) =s' \},$$ so $\tilde{\Gamma}_2 \cong \Gamma_2
\times \R$ by $(x', x, s) \mapsto \left((x', s+c(x)), (x,
s)\right)$. Let $\pi_2: \tilde{\Gamma}_2 \cong \Gamma_2 \times \R
\lra \Gamma_2$ be the projection. Then $\pi \circ \tilde{m} = m
\circ \pi_2$, $\pi \circ \tilde{pr}_1 = pr_1 \circ \pi_2$, and
$\pi \circ \tilde{pr}_2 = pr_2 \circ \pi_2$. Let $s_2$ be the
coordinate on $\R$ in $\Gamma_2 \times \R$. Then $s \circ
\tilde{m} = s_2$, $s \circ
\tilde{pr}_1=s_2+c\circ\pi\circ\tilde{pr}_2$, and $s \circ
\tilde{pr}_2 = s_2$.

So, \eqref{*} implies
\[
\begin{split}
& \;\;\;\; (e^{-s} \circ \tilde{m}) \cdot \di (\tilde{m}^* \pi^* \theta) \\
&\quad - (e^{-s} \circ \tilde{pr}_1) \cdot \di (\tilde{pr}_1^*
\pi^* \theta) - (e^{-s} \circ \tilde{pr}_2) \cdot \di
(\tilde{pr}_2^* \pi^*
\theta) \\
&=(e^{-s} \circ \tilde{m})  \cdot \di(s\circ \tilde{m}) \wedge
(\tilde{m}^* \pi^*)
\theta \\
&\quad -(e^{-s}\circ \tilde{pr}_1) \cdot \di (s \circ pr_1) \wedge
\tilde{pr}_1^* \pi^* \theta - (e^{-s} \circ \tilde{pr}_2) \cdot
\di (s \circ pr_2) \wedge \tilde{pr}_2^* \pi^* \theta,
\end{split}
\]
which implies
\[
\begin{split}
& \quad  e^{-s_2} \di (\pi_2^* m^*  \theta)\\
&\quad - e^{-s_2-c\circ pr_2\circ\pi_2} \di (\pi_2^* pr_1^*
\theta) - e^{-s_2} \di ( \pi_2^* pr_2^*
\theta) \\
&=e^{-s_2} \di s_2 \wedge (\pi_2^* m^*  \theta) \\
&\quad - e^{-s_2-c\circ pr_2\circ\pi_2} \di (s_2+c\circ
pr_2\circ\pi_2) \wedge (\pi_2^* pr_1^* \theta) - e^{-s_2} \di
s_2\wedge ( \pi_2^* pr_2^* \theta).
\end{split}
\]
Looking at the part that contains $\di s_2$, we have
\[ e^{-s_2} \di s_2 \wedge ( \pi_2^* m^* \theta - e^{-c\circ pr_2\circ\pi_2}\pi_2^* pr_1^*
\theta -  \pi_2^* pr_2^* \theta)=0, \] so
\[ \pi_2^* \left( m^* \theta - pr_2^*(e^{-c})pr_1^* \theta - pr_2^* \theta \right) =0. \] Since $\pi_2^*$ is
injective, we have
\[ m^* \theta = pr_2^* (e^{-c}) pr_1^*
\theta + pr_2^* \theta. \]

On the other hand,  \eqref{**} implies
\[ \di \left( e^{-s_2} \pi_2^* m^* \theta\right)= \di \left( e^{-s_2} \pi_2^* ( pr_2^* \theta +
pr_2^* (e^{-c} ) pr_1^* \theta) \right), \] which implies
\[\di\left( \tilde{m}^*(e^{-s} \pi^* \theta) \right) =\di\left(\tilde{pr}_1^*(e^{-s} \pi^*
\theta) + \tilde{pr}_2^*(e^{-s} \pi^* \theta)\right),\] so
\[\tilde{m}^* \omega = \tilde{pr}_1^* \omega + \tilde{pr}_2^* \omega.\]
\end{proof}

\begin{remark}
This Proposition only tells us that $\Gamma_c^h(M)$ is a contact
groupoid and that it integrates $T^*M \oplus \R$. To see whether
it is the contact groupoid of $M$, we have to check that it
induces  the same Jacobi structure on $M$ as the one  we started
with. To see this, we only have to check that $\bs$ is a Jacobi
map, since $\bs$ together with the contact structure on
$\Gamma_c^h(M)$ determines the Jacobi structure on $M$.
\end{remark}

\begin{lemma} \label{forms}
With the same notation as in Proposition \ref{prop3}, let us compare the
two source maps: \[\bbs: \Gamma_s^h(M\times\R) \to M\times\R,\quad
\bbs: \Gamma_c^h(M) \to M. \] $\Lambda_{\omega}$and
$\tilde{\Lambda}$ are $\bar{\bs}$--related iff $ \Lambda_{\theta}
$ and $\Lambda$, $E_{\theta}$ and  $E$ are $\bbs$--related.  Here
$\Lambda_{\omega}$ is the Poisson bivector corresponding to
$\omega$ and $(\Lambda_{\theta}, E_{\theta})$ is the corresponding
Jacobi structure of $\theta$.
\end{lemma}
\begin{proof}
Since $\Gamma_c^h(M) \times \R$ is the homogeneous symplectic
manifold of $\Gamma_c^h(M)$, $\Lambda_{\omega}=e^{-s} (
i_{\Gamma*} \Lambda_{\theta} + \frac{\partial}{\partial s} \wedge
i_{\Gamma *} E_{\theta})$, where $(i_{\Gamma})_*$ is defined
analogously to $i_*$. Considering the difference of the two bivector
fields at point $(x, s) \in \Gamma_c^h(M) \times \R$, we have
\[
\begin{split}
&\;\;\;\;\wedge^2 T_{(x, s)}\bar{\bs} \Lambda_{\omega}(x, s) -
\tilde{\Lambda}\left(\bar{\bs}(x,s)\right)\\
&= e^{-s} \left( \wedge^2 T_{(x,s)}\bar{\bs} (i_{\Gamma*}
\Lambda_{\theta})(x,s) + \frac{\partial}{\partial s} \wedge
T_{(x,s)}\bar{\bs} (i_{\Gamma *}
E_{\theta})(x,s)\right)\\
&\quad -e^{-s} \left( i_* \Lambda + \frac{\partial}{\partial s} \wedge i_* E\right)(\bbs(x), s) \\
&= e^{-s} \left(i_*(\wedge^2T_x\bbs \Lambda_{\theta}
-\Lambda)(\bbs(x),s) +
\frac{\partial}{\partial s} \wedge i_* ( T_x \bbs E_{\theta} -E)(\bbs(x),s)\right) \\
&= e^{-s} \left( (\wedge^2T_x \bbs \Lambda_\theta (x) -
\Lambda(\bbs(x)), 0, 0)+ \frac{\partial}{\partial s} \wedge
\left((T_x \bbs E_{\theta} -E) (\bbs(x)), 0 \right)\right).
\end{split}
\]
Since  $\wedge^2T_x \bs \Lambda_\theta (x) - \Lambda(x)$ is a
bivector field which does not contain $\frac{\partial}{\partial
s}$, we have
\[\bar{\bs}_*(\Lambda_{\omega}) - \tilde{\Lambda} =0 \quad\text{iff}\quad \bs_* \Lambda_{\theta}
-\Lambda=\bs_* E_{\theta} -E=0. \]
\end{proof}

So by Proposition \ref{prop3} and Lemma \ref{forms}, we have
\begin{prop} \label{prop5}
The symplectic groupoid $( \Gamma_c^h(M) \times \R , \omega,
\frac{\partial}{\partial s})$ is the homogeneous symplectic
groupoid of $(M\times \R, \tilde{\Lambda})$ iff $(\Gamma_c^h(M),
\theta, e^{-c})$ is a contact groupoid of $(M, \Lambda, E)$.
\end{prop}

\begin{proof} [Proof of Theorem \ref{congpd}]
By Theorem \ref{wgpd-c-s}, we see that $\Gamma_c^h(M)$ is smooth
if and only if $\Gamma_s^h(M\times\R)$ is smooth. So a Jacobi
manifold $M$ is integrable if and only if its homogeneous Poisson
manifold is integrable.  By Proposition \ref{prop2}, the
homogeneous symplectic groupoid $(\Gamma_c^h(M)\times \R, \omega,
\frac{\partial}{\partial s})$ is the unique source-simply
connected symplectic groupoid of $(M \times \R, \tilde{\Lambda})$.
So $(\Gamma_c^h(M), \theta, e^{-c})$ is a contact groupoid of $(M,
\Lambda, E)$ by Proposition \ref{prop5}. Since $\bs^{-1} (x)
\times s = \bar{\bs}^{-1} (x,s)$, it is clear that $\Gamma_c^h(M)$
being $\bbs$-simply connected is equivalent to $\Gamma_s^h(M)$
being $\bar{\bs}$-simply connected. Moreover, by our construction,
$\Gamma_c^h(M)$ integrates $T^*M \oplus_M \R$.

So we only need to show uniqueness. If there is another 1-form
$\theta_1$ and another function $f_1$ that makes $(\Gamma_c^h(M),
\theta_1, f_1)\underset{\bt}{\overset{\bs}{\rightrightarrows}}M$
into a contact groupoid of $(M, \Gamma , E)$, then by uniqueness
of the source-simply connected symplectic groupoid over $(M \times \R, \tilde{\Lambda})$,
we must have an automorphism $\tilde{F}$ of $\Gamma_c^h(M) \times
\R$ that preserves its structure as a symplectic groupoid. Since
$\bar{\bs}(\tilde{F}(x,s))=\bar{\bs}(x,s)=(\bbs(x),s)$,
$\tilde{F}$ must preserve the $\R$-component; i.e. $\tilde{F}(x,s)
=(F(x), s)$ where $F$ is an automorphism of $\Gamma_c^h(M)$ with
$F \circ \pi = \pi \circ \tilde{F}$. Then we have
\[\tilde{F}^* \di (e^{-s} \pi^* \theta_1) = \di(e^{-s}
\pi^* \theta),\] which implies
\[ \di (e^{-s} \pi^* F^* \theta_1) =\di(e^{-s} \pi^* \theta),\]
so
\[e^{-s} \di (\pi^*(F^* \theta_1 -\theta))- e^{-s} \di s \wedge
\pi^* ( F^* \theta_1 - \theta)=0 .\] Since $\di (\pi^*(F^*
\theta_1 -\theta))$ does not contain $\di s$ and $\pi^*$ is
injective, we must have $ F^* \theta_1 = \theta$. By uniqueness of
$f$, $f_1 = F^* e^{-c} $. This implies that $F$ is an isomorphism
preserving the contact groupoid structure.
\end{proof}

Combining Theorem \ref{wgpd-c-s} and Theorem \ref{congpd}, we have
proved Theorem \ref{main} from the introduction.

\chapter{A Further application---Poisson manifolds from the Jacobi point of view}\label{fur-app}
In this chapter, we always assume that $(M, \Lambda)$ is a Poisson
manifold, $\Gamma_s(M)$ the orbit space of $\Gamma_s^m(M\times
\R)$ (or\footnote{According to Theorem \ref{cf}, these two groupoids
have the same orbit space.} $\Gamma_s^h(M\times \R)$) and
$\Gamma_c(M)$ the orbit space of $\Gamma_c^m(M)$ (or
$\Gamma_c^h(M)$). When $\Gamma_s(M)$ is a smooth manifold, it
becomes the source-simply connected symplectic groupoid of $M$.
Similarly, when $\Gamma_c(M)$ is a smooth manifold, it is the
source-simply connected contact groupoid of $M$. We will
study the integrability of Poisson
bivectors  by viewing a Poisson manifold $M$ as a
Jacobi one.  Moreover we will study the relation between the
integrability of $M$ as a Poisson manifold and that as a Jacobi
one. Finally, we will apply the above theory to the  prequantization of
symplectic groupoids.

\section{ Poisson bivectors}
\subsection{Relation between $\Gamma_s(M)$ and $\Gamma_c(M)$
via the Poisson bivector}

When the Poisson bivector $\Lambda$ is integrable as a Lie
algebroid 2-cocycle on $T^*M$, the two groupoids $\Gamma_s(M)$ and
$\Gamma_c(M)$ are related through $\Lambda$. To study the relation
between the two groupoids, we begin with a short exact sequence
containing their Lie algebroids:
\[ 0\lra \R \lra T^*M \oplus_M \R \os{\pi}{\lra} T^*M \lra 0. \]
Here, the Lie bracket on $T^*M \oplus_M \R$ is the one induced  by
the Jacobi structure $(M, \Lambda, 0)$, i.e. for all $(a, u),
(b,v) \in \Omega^1(M) \times \ci(M)$,
\begin{equation}\label{bracket} [(a, u), (b, v) ]= \left([a, b],
\Lambda(a, b) + \sharp \Lambda (a) (v) - \sharp \Lambda(b)
(u)\right),\end{equation} and the natural projection $\pi$ is a
Lie algebroid morphism.

\begin{prop}\label{6.1}
If the symplectic groupoid $(\Gamma_s(M), \Omega)$ has source
fibre with trivial second homology group, then $M$ is also
integrable as a Jacobi manifold. Moreover, the contact groupoid
$\Gamma_c(M)$ is isomorphic as a groupoid to the twisted
semi-direct product Lie groupoid $\Gamma_s(M) \ltimes_c \R$ for
some  Lie groupoid 2-cocycle $c \in C^2(\Gamma_s(M), \R)$.
\end{prop}
\begin{proof}
Let us first recall some general results about the cohomology of Lie
algebroids and Lie groupoids.

1. In general, for any closed $\omega \in C^2(A, \R)$, we can construct
a Lie algebroid structure on the direct sum $A\oplus \R$
\cite{Mac}.   For all $X, Y \in A$ and $x, y \in \R$, the new
bracket is defined by
\[ [(X, x), (Y, y)]=([X, Y]_A , \omega(X,Y) + \call_{\rho(X)}
y-\call_{\rho(Y)} x) . \] It is a Lie bracket exactly because
$\omega$ is closed. The new anchor is the composition of the
anchor of $A$ and the natural projection from $A\oplus\R$ onto
$A$. We denote this Lie algebroid by $A \ltimes_{\omega} \R$. It
turns out that the isomorphism class of $A\ltimes_{\omega} \R$
only depends on the cohomology class of $\omega$ in $H^2(A, \R)$,
i.e. if $\omega_1$ and $\omega_2$ differ by an exact form, then
$A\ltimes_{\omega_1} \R$ and  $A\ltimes_{\omega_2} \R$ are
isomorphic as Lie algebroids.

2. Conversely, a short exact sequence of Lie algebroids over a
certain manifold $M$,
\begin{equation}\label{ses} 0 \lra \R \lra \tA \os{\pi}{\lra} A
\lra 0 \end{equation} gives \cite{Mac} an element in the Lie
algebroid cohomology $H^2(A, \R)$ in the following way: pick any
splitting of \eqref{ses} of vector bundles $\sigma: A\lra \tA$.
For all $X, Y \in A$, let
\[ \omega(X, Y) = [\sigma(X), \sigma(Y)]_{\tA} - \sigma ([X,
Y]_A).
\]
The image of $\omega$ lies in the kernel of $\pi$. So, we can view
$\omega$ as a real-valued 2-form in $C^2(A, \R)$. Furthermore,
$\omega$ is closed because the brackets on $A$  and
$\tA$ are Lie brackets. In fact, it is not hard to see that
$\tA$ is isomorphic to $A \ltimes_{\omega} \R$ as a Lie algebroid.
Different choices of splitting won't change the cohomology class
of $\omega$. Combining this with result 1, we can see that the Lie
algebroid structures of $\tA$ that make \eqref{ses} into a short
exact sequence of Lie algebroids are characterized by $H^2(A,
\R)$.

3. Suppose $A$ can be integrated into a source-simply connected
Lie groupoid $G$. If $G$ has source fibres with trivial second
homology group, then by Theorem 4 in \cite{m-von-est}, $\omega$ can be
integrated into a 2-cocycle $c$ on the groupoid $G$. So,
$A\ltimes_{\omega}\R$ is automatically integrable and its unique
source-simply connected Lie groupoid is $G\ltimes_c\R$, with
multiplication given by
\[ (g, x) \cdot (h, y) = ( gh, x+y+ c(g, h) ). \]

The proof of the theorem is now straightforward. Notice that
$\Lambda \in \wedge^2 TM$ is closed in the Lie algebroid complex
$(C^n(T^*M, \R), d_{\Lambda})$. Have a closer look at
\eqref{bracket}, then we can see that $T^*M \oplus_M \R \cong T^*M
\ltimes_{\Lambda} \R$. Therefore, under the conditions stated in
this theorem, $T^*M \oplus_M \R$ is integrable and the contact
groupoid integrating it is $\Gamma_s(M) \ltimes_c \R$ for some
closed 2-cocycle $c$ integrating $\Lambda$.
\end{proof}

\begin{remark}
In this case, the symplectic form on $\Gamma_s(M)$ is exact, and
$\Gamma_c(M)$ devided by some suitable $\Z$-action will be
the pre-quantization of $\Gamma_s(M)$.  Please refer to the
subsection 6.3 for details.
\end{remark}

\subsection{General case---without assuming integrability}
When exactly will $\Lambda$ be integrable? To answer this
question, we should look more carefully into the contact groupoid
$\Gamma_c(M)$. There is a natural projection $pr: \Gamma_c(M) \lra
\Gamma_s(M)$, by $[(a_1, a_0)] \mapsto [a_1]$. It is well-defined
because if $(a_1, a_0) $ is an $A$-path of $T^*M \oplus_M \R$, $a_1$ is also an $A$-path of $T^*M$. So there is a
short exact sequence of groupoids,
\begin{equation}\label{ses2}
 1\lra \G \lra \Gamma_c(M) \os{\pi}{\lra} \Gamma_s(M) \lra
 1.
\end{equation}

In fact, $\R$ acts on $\Gamma_c(M)$  by
\[s \cdot [a_1, a_0] = [a_1, a_0 +s]. \]
It is well defined because  we use only differentiation in defining
``$\sim$''. Therefore  $$(a_1(1, t), a_0(1, t)) \sim (a_1(0, t) ,
a_0(0, t)) $$ is equivalent to $$ (a_1(1, t), a_0(1, t) +s) \sim
(a_1(0, t) , a_0(0, t) +s).$$ However, this action is not always
free. When it is free, $\G$ will simply be the trivial groupoid
$\R \times M$ over $M$ and $\Gamma_c(M)$ will be isomorphic to
$\Gamma_s(M) \times \R$ as in the last theorem. It turns out that
$\G$ is closely related to the monodromy groups of the two Lie
algebroids.

Let us first recall some facts from \cite{cf} and \cite{cf2} about
monodromy groups.

\begin{defi} \cite{cf}
Let $A$ be a Lie algebroid over $X$ with anchor $\rho$ and $\g_x(A)$
the isotropy Lie algebra $\ker_x(\rho)$. The \emph{monodromy
group} $N_x(A)$ of $A$ at a point $x \in X$ consists of those
elements in the center of $\g_x$  which, as  constant $A$-paths,
are homotopic to the trivial $A$-path $0_x$.
\end{defi}

Let $L\subset X $ be a leaf through the  point $x\in X$ and
$\G(\g_x(A)) $ the Lie group integrating $\g_x(A)$. Then there is
a homomorphism $\partial: \pi_2(L, x) \to \G(\g_x)$ defined as
follows \cite{cf}: let $[\gamma] \in \pi_2(L, x) $ be represented
by a smooth map $\gamma: I\times I \to L$ which maps the boundary
into $x$. One can always choose $A$-paths $a(\epsilon, \cdot)$ and
$A$-paths $b(\cdot, t) $  over $\gamma$ in $A|_L$
satisfying \eqref{homotopy} and the boundary conditions $a(0,
t)=b(\epsilon, 0)=b(\epsilon, 1)=0$. For example, we can ask that
$b(\epsilon, t)= \sigma(\frac{\di}{\di \epsilon} \gamma(\epsilon,
t))$ where $\sigma: TL \to A|_L$ is any splitting of the anchor,
and take $a$ to be the unique solution of \eqref{homotopy} with
initial condition $a(0, t)=0$. Since $\gamma(1, t)$ is the
constant path $0_x$, $a(1, t)$ must lie in $\g_x(A)$ entirely. As
a path in the Lie algebra $\g_x(A)$, $a(1,t)$ can be integrated
into a path $g(1, t)$ in $\G(\g_x)$ \cite{lie-group} or \cite{cf}.
Then $\partial([\gamma])$ is defined as
$\partial([\gamma])=g(1,1)$.

The map $\partial$ fits into the exact sequence:
\begin{equation} \label{ses3}
 \pi_2(L, x) \os{\partial}{\lra} \G(\g_x(A)) \lra \G(A)_x
\lra \pi_1(L, x),
\end{equation}
where $\G(A)_x:= \bs^{-1} (x) \cap \bt^{-1} (x)$. The map from
$\G(\g_x(A))$ to $\G(A)_x$ is defined by mapping each equivalence
class $[a]$ in $\G(\g_x(A))$ to the equivalence class $[a]$ in
$\G(A)$. Since every two $A$-paths equivalent as $A$-paths in
$\g_x(A)$ must be equivalent as $A$-paths in $A$, this map is well
defined. The map from $\G(A)_x$ to $\pi_1(L, x)$ is simply defined
by sending the equivalence classes of $A$-paths to the equivalence
classes of their base paths.

The image of $\partial$ in $\G(\g_x)$ is defined as
$\tilde{N}_x(A)$ by Crainic and Fernandes in \cite{cf}; it is
closely related to $N_x$. Actually, $\tilde{N}_x(A)$ is a
subgroup of $Z(\G(\g_x(A)))$---the center of $\G(\g_x(A))$, and
it intersection with the connected component of the identity of
$Z(\G(\g_x(A)))$ is isomorphic to $N_x(A)$ by the exponential map
on the Lie algebra $\g_x$.

Returning to our case, where $(M, \Lambda)$ is a Poisson manifold,
the two Lie algebroids $T^*M$ and $T^*M\oplus_M\R$ induce the same
leaves on $M$, namely the symplectic  leaves of $M$. On a leaf $L$
through a point $x\in M$, let $\omega_L$ be the symplectic form
induced by $\Lambda$, and $\partial_c$ and $\partial_s$ the
homomorphisms from $\pi_2(L, x)$ to $\G(\g_x(T^*M))$ and
$\G(\g_x(T^*M\oplus_M\R))$ respectively. Define the group
\[ Per_0(\omega_L):=\{ \int_{\gamma} \omega_L : [\gamma] \in
\pi_2(L,x) \; \text{and} \; \partial_s\gamma=1_x\}.\] It is a
subgroup of the period group of $\omega_L$
\[
Per(\omega_L):=\{\int_{\gamma} \omega_L : [\gamma] \in \pi_2(L,x)
\}.\]

In general, even without assuming the integrability of $T^*M$ or
$T^*M \oplus_M \R$, we have the following theorem:

\begin{thm}\label{relation2}
Let $pr$ be the projection from $\Gamma_c(M)$ to $\Gamma_s(M)$
as defined above. Then $\G$ in \eqref{ses2} is a bundle of groups
over $M$. Furthermore, at each point $x\in M$, $\G_x= \R /
Per_0(\omega_{L_x})$, where $L_x$ is the leaf through $x$.
\end{thm}

Before proving this theorem, let us first prove a useful lemma.

\begin{lemma}\label{partial}
Let $L$ be a leaf through a point $x \in M$.  Then
\[ \partial_c \gamma = ( \partial_s \gamma, - \int_{\gamma}
\omega_L). \] for every $\gamma$  representing $[\gamma] \in
\pi_2(L, x)$.
\end{lemma}

\begin{proof}
Let $(a, u) $ and $(b, v)$ be $A$-paths in $T^*M \oplus_M \R$ over
$\gamma$ satisfying \eqref{homotopy} and the boundary conditions:
$$a(0, t)= b(\epsilon, 0) = b(\epsilon, 1)=0 \in T^*M,$$
and $$u(0, t)=v(\epsilon, 0) =v(\epsilon, 1) =0 \in \R.$$ Writing
out equation \eqref{eq:homotopy} on the $\R$-component more
explicitly, we have
\[\partial_t v -\partial_{\epsilon} u = \Lambda (a, b). \]
Notice that
\[ \sharp \Lambda (a) = \frac{\di}{\di t} \gamma,\;\;\;\;\; \sharp \Lambda (b) = \frac{\di}{\di \epsilon}
\gamma, \] and $\gamma$ stays entirely in the leaf $L$. We have
\[ \partial_t v - \partial_{\epsilon} u = \omega_L ( \frac{\di}{\di t} \gamma,
\frac{\di}{\di \epsilon} \gamma). \] So
\[
\begin{split}
\int_I \di \epsilon \int_I \di t \omega_L( \frac{\di}{\di t}
\gamma, \frac{\di}{\di \epsilon} \gamma)&= \int_I \di \epsilon
\int_I \partial v \di t -\int_I \di
t \int_I \partial_{\epsilon} u \di \epsilon\\
&= \int_I (v(\epsilon , 1)- v(\epsilon , 0))- \int_I \di t (
u(1,t) u(0,t)) \\
&=-\int_I u(1, t) \di t,
\end{split}
\]
i.e.  $\int_{\gamma} \omega_L = -\int_I u(1,t) \di t$.

The brackets on $\g_x(T^*M \oplus_M \R) $ and $\g_x(T^*M)$ are
induced from $T^*M \oplus_M \R$ and $T^*M$ respectively.
$(X,\lambda)$ and $(Y, u) \in \g_x (T^*M \oplus_M\R)$ can be
extended  to sections $(\tilde{X}, \tilde{\lambda})$ and
$(\tilde{Y}, \tilde{\mu})$ in $T^*M \oplus_M\R$ such that
$\tilde{\lambda}$ and $\tilde{\mu}$ are locally constant functions
around point $x$. Then,
\[\begin{split}
[(X, \Lambda), (Y, \mu)]_{\g_x(T^*M\oplus_M\R)} &:= [(\tilde{X},
\tilde{\lambda}),(\tilde{Y}, \tilde{\mu})]_{\T^*M \oplus_M\R}(x)\\
&=([\tilde{X}, \tilde{Y}]_{T^*M}(x), \sharp
(\tilde{X})(\tilde{\mu})-\sharp\Lambda(\tilde{Y})(\tilde{\lambda})+\Lambda(\tilde{X},\tilde{Y})(x))\\
&= ([X,Y]_{\g_x(T^*M)} , 0).
\end{split}
\]
So $\g_x(T^*M\oplus_M\R)$ is isomorphic to $\g_x(T^*M)\oplus \R$
as a Lie algebra. Therefore, as Lie groups,
$\G(\g_x(T^*M\oplus_M\R))=\G(\g_x(T^*M))\times \R$.

Then $\partial_c \gamma$, defined as the end point of the
integration path of $(a(1,\cdot), u(1,\cdot))$, has the first
component the end point of the integration path of $a(1,\cdot)$
and the second component $\int_I u(1,t)\di t=-\int_{\gamma}
\omega_L$. Therefore, we have $\partial_c \gamma=(\partial_s
\gamma, -\int_{\gamma} \omega_L)$.
\end{proof}
\begin{remark}\label{regular}
When $M$ is a regular Poisson  manifold, in the definition of Lie
brackets on $\g_x(T^*M \oplus_M \R) $ and $\g_x(T^*M)$, we can
choose the extension $\tilde{X}$ and $\tilde{Y}$ both lying in
$\g_y(T^*M)$ for all $y$ in a neighborhood of $x$. Therefore, the Lie brackets of
$\g_x(T^*M\oplus_M\R)$ and $\g_x(T^*M)$ are both 0. So the Lie
groups $\G(\g_x(T^*M\oplus_M\R))$ and $\G(\g_x(T^*M))$ are abelian
and isomorphic to their Lie algebras.
\end{remark}

Now we are ready to prove Theorem \ref{relation2}.
\begin{proof} [Proof of Theorem \ref{relation2}]
At a point $x \in M$, by definition, we have
\[ \begin{split}
\G_x = \pi^{-1} ( [1_x])=& \{ [1_x, u]: (1_x, u)\; \text{is an $A$-path in}
T^*M \oplus_M \R \\
 & \text{with the constant path $x$ as its base path}\}.
\end{split}
 \]
Notice that $(1_x, u) \sim (1_x, \int_I u(t)\di t)$ by the
homotopy $(b(\epsilon, t), v(\epsilon, t))=(0_x, -\int^t_0 u(s)\di
s + t\int_I u(s) \di s)$. We can rewrite $\Sigma_x$ as:
\[\G_x=\{ [(1_x,c)]: (1_x, c)\; \text{is a constant $A$-path in}\;
T^*M\oplus_M\R \;\text{over}\; x\}. \] By the definition of
monodromy groups and their close relation to $\tilde{N}_x$, we
have
\[\g_x=\R / \tilde{N}_x(T^*M\oplus_M\R) \cap 1_x \times \R, \]
because $1_x \times \R$ lies in the connected component of the
identity $(1_x, 0)$.

By Lemma \ref{partial},
\[\begin{split}
\tilde{N}_x (T^*M\oplus_M\R)&=\{\partial_c \gamma : \gamma\in
\pi_2(L,
x) \}\\
&=\{(\partial_s \gamma, -\int_{\gamma}\omega_{L_x}), [\gamma] \in
\pi_2(L,x)\}.
\end{split}
\]
So
\[
\begin{split}
\tilde{N}_x(T^*M\oplus_M\R)\cap 1_x \times\R&=1_x \times \{
-\int_{\gamma} \omega_{L_x}:\;\partial_s \gamma=1_x,\; \gamma \in
\pi_2(L,x)\}\\
&=1_x \times Per_0(\omega_{L_x}). \end{split} \] Therefore
$\G_x=\R / Per_0(\omega_{L_x})$.
\end{proof}

\subsection{The integrable case} \label{int-case}
With the same setting as in Theorem \ref{relation2}, if we assume the
integrability of $T^*M$, then $\Gamma_s(M)$ is a symplectic
groupoid with symplectic 2-form $\Omega$. We can express the group
$\G_x$ in terms of the period group of $\Omega|_{\bs^{-1}(x)}$,
which is defined as
\[ Per(\Omega|_{\bs^{-1}(x)})= \{ \int_{g} \Omega: [g] \in
\pi_2(\bs^{-1}(x))\} \]

\begin{cor} \label{per}
If $T^*M$ is integrable, i.e. if  $(\Gamma_s(M), \Omega)$ is a
symplectic groupoid, then the group $\G_x = \R/
Per(\Omega|_{\bs^{-1}(x)})$. \end{cor}
\begin{proof}
On an $\bs$-fibre $\bs^{-1}(x)$ of $\Gamma_x(M)$, $\bt:
\bs^{-1}(x) \to L$ is a submersion. We also know that $\bt$ is an
anti-Poisson map, so $\bt^* \omega_L =-\Omega$. Examining
$\partial $ more carefully, it is not hard to see that $\partial
\gamma =1_x$ means exactly that $\gamma(\epsilon, t) $ can be
lifted to g-paths (i.e. paths inside the source fibre of the
groupoid) $g(\epsilon, t)$ such that $g(0,t)=g(1,t)=g(\epsilon,
0)=g(\epsilon, 1)=1_x$. This tells us that $\gamma$ can be lifted
to a 2-cocycle $g$ in $\bs^{-1} (x) $. They satisfy
\[
\int_{\gamma} \omega_{L}=\int_{g}\bt^* \omega_{L}=-\int_g \Omega.
\]
So we have
\[\begin{split}
\G_x&= \R /  \{ -\int_{\gamma} \omega_{L}:\;\partial_s
\gamma=1_x,\;
\gamma \in \pi_2(L,x)\} \\
&=\R/\{ \int_g \Omega, [g] \in \pi_2(\bs^{-1}(x))\}\\
&=\R/ Per(\Omega|\bs^{-1} (x)).
\end{split}
\]
\end{proof}

When $\Gamma_c(M)$ is a smooth manifold, there is a 1-form
$\theta$ such that $(\Gamma_c(M), \theta ,1)$ is the source-simply
connected contact groupoid of $M$. There is an $\R$-action on
$\Gamma_c(M)$ given by
\[ s\cdot[(a_1, a_0)]=[(a_1, a_0+s)] \]
It is well defined. In fact, we have the following lemma:
\begin{lemma} \label{action}
Any $A$-path $(a_1, a_0) $ in  $T^*M\oplus_M\R$ has the following
property:
\[[(a_1, a_0+s)]=[(a_1, a_0)]\cdot[(0_x, s)]=[(0_y,s)]\cdot[(a_1,
a_0)],\] where $x=\gamma(0)$ and $y=\gamma(1)$ are the end points
of the base path $\gamma$, and $(0_x, s)$, $(0_y, s)$ are constant
paths over $x$ and $y$ respectively.
\end{lemma}
\begin{proof}
Choose a suitable cut-off function $\tau \in \ci (I, I)$ with the
property that $\tau' $ is zero near 0 and 1. For any path $c$,
denote the reparameterization of $c$ by $\tau$ by
$c^{\tau}=\tau'c(\tau(t))$. Note the following facts about
$A$-paths in $T^*M\oplus_M\R$:

\begin{itemize}

\item If  $\int_I a_0 = \int_I a^*_0$, then $(a_1, a_0) \sim
(a_1, a^*_0)$ through $\left(0, \int_0^t \epsilon a_0(t) +
(1-\epsilon) a^*_0(t)\right)$.\\

\item $(a_1, a_0) \sim (a_1^{\tau}, a_0)$ through
$((\tau(t)-t)a_1((1-\epsilon)t+(1-\epsilon)\tau(t)), 0)$.
\end{itemize}

Then
\[
\begin{split}
[(a_1, a_0)]\cdot[(0_x, s)] &= [(a_1^{\tau}, a_0^{\tau})]
\cdot[(0_x^{\tau}, s^{\tau})]\\
&=[(a_1^{\tau}, a_0^{\tau}\odot s^{\tau})].
\end{split}
\]
Since $\int_I a_0^{\tau} \odot s^{\tau} \di t= \int_I a_0^{\tau} +
\int s^{\tau} = \int_I a_0 + s$, we have
\[ (a_1^{\tau}, a_0^{\tau} \odot s^{\tau}) \sim ( a_1^{\tau}, a_0
+s) \sim ( a_1, a_0 +s).\] Therefore $[(a_1, a_0)]\cdot[(0_x,
s)]=[( a_1, a_0 +s)]$.

Similarly $[(0_y,s)]\cdot[(a_1, a_0)]=[(a_1, a_0+s)]$.
\end{proof}

This action is not always free. In fact, with the above lemma, it
is easy to conclude that it is free iff $\G_x=\R$ for all $x\in
M$.
The vector field $\frac{\partial}{\partial s}$ generating this
action always has orbits $\R^1$ or $S^1$ when $M$ is integrable as
a Jacobi manifold. Please refer to subsection 6.3 for details.

By a calculation in local coordinates, we can see that
$\frac{\partial}{\partial s}$ is the Reeb vector field of
$\theta$, i.e. \begin{equation} \label{basic}
 \call_{\frac{\partial}{\partial s}}
\theta =0, \;\;\;\;\; i(\frac{\partial}{\partial s}) \theta =1.
\end{equation}
This tells us that $\di \theta$ is basic, i.e. there is a 2-form
$\omega$ on $\Gamma_s(M)$, such that $\di \theta =\pi^* \omega$.
$\omega$ is obviously closed. Moreover, it is nondegenerate and
multiplicative. This  follows from the nondegeneracy and
multiplicativity of $\theta$. Therefore it is a multiplicative
symplectic 2-form on $\Gamma_s(M)$. It is easy to check that the
source map $\bs_s: (\Gamma_s(M), \omega) \to(M, \Lambda)$ is a
Poisson map because $\bs_c: (\Gamma_c(M),
\theta) \to (M, \Lambda)$ is a Jacobi map. So $(\Gamma_s(M),
\omega)$ is the source-simply connected symplectic groupoid of
$(M, \Lambda)$. By uniqueness, we must have $\omega=\Omega$.
Therefore $\pi^* \Omega =\di \theta$.

\subsection{Integrability of Poisson bivectors}

\begin{thm}\label{main2'}
Suppose that $(M, \Lambda)$ is integrable as a Poisson manifold
and $(\Gamma_s(M), \Omega)$ is the source-simple connected
symplectic groupoid of $M$. Then the following statements are
equivalent:
\begin{enumerate}
\item The symplectic 2-form $\Omega$ is exact; \\
\item The period groupoid $Per(\Omega|_{\bs^{-1}(x)})=0$;\\
\item  The group bundle $\G$ is the trivial line bundle $\R \times M$;\\
\item The Poisson bivector $\Lambda$ is integral as a Lie algebroid 2-cocycle on
$T^*M$;\\
\item $M$ is integrable as a Jacobi manifold and as a groupoid $\Gamma_c(M)=
\Gamma_s(M) \ltimes_c \R$, for some groupoid 2-cocycle $c$ on
$\Gamma_s(M)$.
\end{enumerate}
\end{thm}
\begin{proof}
From Theorem \ref{relation2} and Proposition \ref{6.1}, it's easy
to see that (1) $\Rightarrow$ (2) $\Leftrightarrow$ (3)
$\Leftarrow$(5)$\Leftarrow$(4). So we only have to show that
(3)$\Rightarrow$(1) and (2)$\Rightarrow$(4).

``(3)$\Rightarrow$(1)'': Since $\G_x=\R$, the $\R$-action we
constructed earlier is free. So $\Gamma_c(M) \os{\pi}{\lra}
\Gamma_s(M)$ is a $\R$-principal bundle. By \eqref{basic},
$\theta$ is a connection 1-form of this bundle and
$\pi^*\Omega=\di \theta$ shows that $\Omega$ is the curvature
2-form. Since $\R$ is contractible, $\pi^*$ induces an isomorphism
from $H^2(\Gamma_s(M))$ to $H^2(\Gamma_c(M))$. So $[\pi^*
\Omega]=[\di \theta] $ shows that $[\Omega]=0$, i.e. $\Omega$ is
exact.

``(2)$\Rightarrow$(4)'': If we view $\Lambda \in \wedge^2
TM=\wedge^2(T^*M)^*$, then right translation can move it along
the $\bs$-fibres and make it into a 2-form $\Omega_{\Lambda}$ on
the $\bs$-fibres. By a theorem in \cite{m-von-est}, $\Lambda$ is
integrable if and only if the period group
$Per(\Omega_{\Lambda})=0$. Here, $\Omega_{\Lambda}=
\Omega|_{\bs}^{-1}(x)$. Therefore (2) is equivalent to (4).
\end{proof}

\begin{cor}
If every symplectic leaf in an integrable Poisson manifold $M$ has
exact symplectic form, then the symplectic form $\Omega$ of
$\Gamma_s(M)$ is also exact.
\end{cor}
\begin{proof}
This is a direct conclusion from the above theorem and Theorem
\ref{relation2}.
\end{proof}

\section{Integrability of Poisson manifolds as Jacobi manifolds}

The integrability of $M$ as a Poisson manifold and as a Jacobi
manifold are closely related, but they are not equivalent. In the next chapter, we will see a Poisson
manifold which can be integrated into a contact groupoid but not into a
symplectic one.  In this subsection, we deal with the other
direction, i.e. we assume that $M$ is integrable as a Poisson
manifold and describe its integrability as a Jacobi manifold in
terms of the group bundles $\G$ and $P:=\sqcup_x Per(\Omega|_{\bs^{-1}
(x)})$.

\begin{prop}\label{p}
Suppose that a Poisson manifold $(M, \Lambda)$ can be integrated
into the symplectic groupoid $(\Gamma_s(M), \Omega)$. Then $M$ is
integrable as a Jacobi manifold if and only if $P$ is uniformly
discrete.
\end{prop}
\begin{proof}
By the main theorem in \cite{cf}, $M$ is integrable as a Jacobi
manifold if and only if the groups $\tilde{N}_x(T^*M\oplus_M\R)$
are uniformly discrete. Recalling Lemma \ref{partial}, $\tilde{N}
(T^*M\oplus_M\R)$ is uniformly discrete if sequences
$[\gamma_i]\in\pi_2(L,x_i)$ and $ x_i=x$  satisfy
\begin{equation} \label{lim} \lim_{n\to + \infty} \text{distance}
((\partial_s \gamma_i, -\int_{\gamma_i } \omega_L), (1_{x_i},
0))=0,
\end{equation}
and
\[ \lim_{i\to + \infty}x_i=x,\]
then $(\partial_s\gamma_i, -\int_{\gamma_i } \omega_L)=(1_{x_i},
0)$ for $i$ large enough. Condition \eqref{lim} is equivalent to
$\lim_{i\to +\infty}
\partial_s\gamma_i=1_{x_i}$ and $\lim_{i\to + \infty}-\int_{\gamma_i}\omega=0$. Since
$T^*M$ is integrable, the groups $\tilde{N}_x(T^*M)$ are uniformly
discrete. Therefore $\lim_{i\to+\infty}
\partial_s\gamma_i =1_{x_i}$ implies $\partial_s\gamma_i=1_{x_i}$
for $i$ large enough.

Rephrasing \eqref{lim}, if $\tilde{N}_x(T^*M\oplus_M\R)$
is uniformly discrete then there exist  sequences
$\partial_s\gamma_i =1_{x_i}$ with $\lim_{i\to \infty}x_i =x$
satisfying
\[\lim_{i\to
+\infty} \int_{\gamma_i} \omega=0.\]This implies $\int_{\gamma_i}
\omega_L\equiv 0$ for $i$ large enough.

By Corollary \ref{per},
$Per(\Omega|{\bs^{-1}(x)})=Per_0(\omega_L)$. So the above
condition is  exactly the condition that $P$ is uniformly
discrete.
\end{proof}

\begin{thm}
If a  Poisson manifold $M$  can be integrated into the symplectic
groupoid $(\Gamma_s(M), \Omega)$, then the following statements
are equivalent:
\begin{enumerate}
\item $M$ is integrable as a Jacobi manifold;
\item $\G$ is a Lie groupoid over $M$;
\item $P$ is an \'etale groupoid over $M$.
\end{enumerate}
\end{thm}
\begin{proof}

(1)$\Rightarrow$(2): Recall the short exact
sequence \eqref{ses2}.  When $\Gamma_c(M)$ is a manifold, the
projection $\pi$ is a submersion by definition.  So $\G
\cong \pi^{-1}(M)$ is a smooth submanifold of $\Gamma_c(M)$.
Moreover, $\G$ also inherits a groupoid structure from
$\Gamma_c(M)$. The source map $\bs_{\G}: \G \to M$, sending
$[(1_x, a_0)] $ to $x$ is obviously a submersion. Similarly, the
same result also holds for the target map.  So $\G$ is a Lie
groupoid.

(2)$\Rightarrow$(3): As in the proof above, in the
short exact sequence
\[ 1\lra P \lra\R\times M \os{\phi}{\lra} \G \lra 1,\]
$\phi$ is a submersion. This tells us that $P$ is a closed
submanifold of $\R \times M$. Moreover, $P$ also inherits a
groupoid structure from the trivial groupoid $\R \times M$. The
source map $\bs_P: P \to M$, sending $(x, \int_{\gamma}
\Omega|{\bs^{-1}(x)})$ to $x$, is obviously a submersion.
Similarly, the target map is also submersion. So $P$ is a Lie
groupoid.

The period group $Per(\Omega|{\bs^{-1}(x)})$ is a closed subgroup
of $\R$ because $M$ is a closed submanifold of $\R \times M$.
However, $Per(\Omega|{\bs^{-1}(x)})$ contains at most countably many
elements since second homotopy groups of manifolds are always
countable. So $Per(\Omega|{\bs^{-1}(x)})$ must be discrete.
Therefore, $P$ is an \'etale groupoid.

``(3)$\Rightarrow$(1)'': $P$ being \'etale implies that $P$ is
uniformly discrete. By Proposition \ref{p}, $M$ is integrable as a
Jacobi manifold.
\end{proof}

\section{Relation to prequantization}
Prequantizations of symplectic groupoids were introduced by
Weinstein and Xu in \cite{xw1}, as the first steps of quantizing
symplectic groupoids, for the purpose of quantizing Poisson
manifolds ``all at once''.

In general, a {\em prequantization} $E$ of an symplectic manifold
$(S, \omega)$ is a $S^1$-principal bundle over $S$ with connection
1-form $\theta$  that  has curvature 2-form $\omega$. It
turns out that $(S, \omega)$ is prequantizable if and only if
$\omega$ represents an integral class in $H^2(S, \Z)$.

Generally, the isomorphism classes of principal $S^1$-bundles over any manifold $X$
form an abelian group $\mathcal{P}(X, S^1)$ with ``tensor product'', which is
isomorphic to $H^2(X, \Z)$. The isomorphism is constructed as
follows: for a principal $S^1$-bundle $E$, the curvature 2-form
$\omega$ on $X$ is an integral class and doesn't depend on the
choice of connection. The class $[\omega]$ is called the {\em
characteristic class} of $E$. On the other hand, for any integral
2-form $\omega$, there exists a principal $S^1$-bundle $E$ with
characteristic class represented by $\omega$ \cite{koba63}.
Therefore, there is a \emph{unique} principal $S^1$-bundle $E$
serving as a prequantization for a symplectic manifold $S$ with
integral class; moreover, when $S$ is simply connected, the
cohomology class of connections on $E$ is also unique \cite{bl}.

In our case, the prequantization of the symplectic groupoid
$(\Gamma_s(M), \Omega)$ is closely related to
$\Gamma_c(M)$. $\Z$ as a subgroup of $\R$ acts naturally on
$\Gamma_c(M)$. With this $\Z$-action, we can prove Theorem
\ref{main3}.

Let us recall the content of the theorem. It says:
If $(\Gamma_s(M), \Omega)$ is a symplectic groupoid with
$\Omega\in H^2(\Gamma_s(M), \Z)$, then $M$ can be integrated into
a contact groupoid $(\Gamma_c(M), \theta, 1)$. Furthermore, if we
quotient out by a $\Z$ action, $\Gamma_c(M)/\Z$ is a
prequantization of $\Gamma_s(M)$ with connection 1-form
$\bar{\theta}$ induced by $\theta$. Moreover, $(\Gamma_c(M)/\Z,
\bar{\theta}, 1)$ is also a contact groupoid of $M$.

\begin{proof}[Proof of Theorem \ref{main3}]
First of all, when the symplectic form $\Omega$ on $\Gamma_s(M)$
is an integral class, $Per(\Omega|_{\bs^{-1}(x)})$ is always a
subset of $Per(\Omega)\subset$ the trivial $\Z$-bundle. So $P$ is
always uniformly discrete. Therefore $\Gamma_c(M)$ is
automatically a Lie groupoid.

By Lemma \ref{basic}, $\theta$ is $\R$-invariant, so it descends
to $\Gamma_c(M)/\Z$, i.e. there is a 1-form $\bar{\theta}\in
\Omega^1(\Gamma_c(M)/\Z)$ such that $\pi_{\Z}^*
\bar{\theta}=\theta$, where $\pi_{\Z}$ is the projection from
$\Gamma_c(M)$ to $\Gamma_c(M)/\Z$. Since $\theta\wedge(\di
\theta)^n \neq 0$, we have
$\bar{\theta}\wedge(\di\bar{\theta})^n\neq 0$ too, where
$n=\frac{1}{2}(\dim \Gamma_c(M) -1)$. So, with this 1-form
$\bar{\theta}$, $\Gamma_c(M)/\Z$ is a contact manifold.

Moreover, the Reeb vector field $\frac{\partial}{\partial s}$ can
also descend to a vector field $E$ on $\Gamma_c(M)/\Z$ and becomes
the Reeb vector field of $\bar{\theta}$, i.e.
\[ \call_{E} \bar{\theta}=0 ,\;\;\;\;\; \iota(E) \bar{\theta}=1. \]

Since the period group $Per(\Omega|{\bs^{-1}(x)})$ is a subgroup
of $\Z$, the $S^1$-action on $\Gamma_c(M)/\Z$ induced by the
$\R$-action on $\Gamma_c(M)$ is free and the projection $\pi:
\Gamma_c(M)\to \Gamma_s(M)$ factors through to $\pi_s :
\Gamma_c(M)/\Z \to \Gamma_s(M)$ . Then $\Gamma_c(M)/\Z
\os{\pi_s}{\to} \Gamma_s(M)$ is a principal $S^1$-bundle. By
reasoning similar to that in  Section \ref{int-case}, $\bar{\theta}$ is the
connection 1-form of the $S^1$-principal bundle and $\Omega$ is
the curvature 2-form. So $\Gamma_c(M)/\Z$ is a prequantization
bundle of $\Gamma_s(M)$.

Moreover, the source and target maps from $\Gamma_c(M)$ to $M$ are
$\R$-equivariant. So we can define the source and target maps
$\bar{\bs}$, $\bar{\bt}$ from $\Gamma_c(M)/\Z$ to $M$ as
$\bar{\bs}(h+\Z)= \bs(h)$, and similarly for $\bar{\bt}$, for all
$h\in\Gamma_c(M)$.

If $\bar{\bs}([(a_1, a_0)] + \Z)=\bar{\bt}([(a_1^*, a_0^*)] +
\Z)$, then $ \bs[(a_1, a_0)]=\bt[(a_1^*, a_0^*)])$. We can define
the multiplication by
\[ ([(a_1, a_0)] + \Z)\cdot ([(a_1^*, a_0^*)] + \Z) =[(a_1,
a_0)]\cdot[(a_1^*, a_0^*)] + \Z. \] Notice that
\[
\begin{split}
[(a_1, a_0+s)]\cdot[(a_1^*, a_0^*+t)]&= [(a_1, a_0)]\cdot[(0_x,
s)] \cdot [(a_1^*, a_0^*+t)]\\
&=[(a_1, a_0)]\cdot[(a_1^*, a_0^*+s+t)]\\
&=[(a_1, a_0)]\cdot[(a_1^*, a_0^*)]\cdot[(0_y, s+t)],
\end{split}
\]
so the multiplication is well defined.

Viewing any $x\in M$ as a constant path $0_x$, we have the identity
section
\[ M \hookrightarrow \Gamma_c(M)/\Z, \;\;\;\;\; x \mapsto [(0_x, 0)] +\Z.\]

Moreover, for any $[(a_1, a_0)] + \Z \in \Gamma_c(M)/\Z$, its
inverse element is just $[(\bar{a}_1, \bar{a}_0)] + \Z$, where
$\bar{c}(t)= c(1-t)$ for any path $c$.

It is routine to check that the above  gives us a Lie groupoid
structure on $\Gamma_c(M)/\Z$. It is also easy to see that the
multiplicativity of $\bar{\theta}$ follows from that of
$\theta$. Moreover, $\bar{\bs}$ is a Jacobi map because $\bs_c$ is
a Jacobi map. Therefore, $(\Gamma_c(M)/\Z, \bar{\theta}, 1)$ is a
contact groupoid of $M$.
\end{proof}

Notice that in the proof we have only used the fact that
$Per(\Omega|_{\bs^{-1}}) \subset \Z$ to construct the principal
bundle structure for $\Gamma_c(M)/\Z$. We have the following
corollary:

\begin{cor}
The symplectic groupoid $(\Gamma_s(M), \Omega)$ is prequantizable
if  $Per(\Omega|_{\bs^{-1}(x)})$ $\subset$ $\Z  $. \end{cor}

With the same hypotheses, combining Theorem \ref{relation2} and
Corollary \ref{per}, we have:
\begin{cor}
The symplectic groupoid $\Gamma_s(M)$ is prequantizable if every
leaf of $M$ has an integral symplectic form.
\end{cor}





\chapter{Examples}
In this chapter, we give examples on Weinstein groups (when the
basis manifold of a Weinstein groupoid is a point, we call it a
Weinstein group) and contact groupoids.

\section{Weinstein groups}

\begin{ep}[ $B\Z_2$] $B\Z_2$ is a Weinstein group (i.e. its base space is a point) integrating the trivial Lie algebra
0. The \'etale differentiable stack $B\Z_2$ is presented by
$\Z_2\rightrightarrows \cdot$ (here $\cdot$ represents a point).
We establish all the structure maps on this presentation.

The source and target maps are just projections from $B\Z_2$ to a
point.

The multiplication $m$ is defined by \[m: (\Z_2 \rightrightarrows
\cdot)\times (\Z_2 \rightrightarrows pt) \to (\Z_2
\rightrightarrows pt), \; \text{by} \; m(a, b)= a\cdot b,\] where
$a, b\in \Z_2$. Since $\Z_2 $ is commutative, the multiplication
is a groupoid homomorphism (hence gives rise to a stack
homomorphism). It is easy to see that $m\circ (m\times id) =m
\circ (id \times m)$, i.e. we can choose the 2-morphism $\alpha$
inside the associativity diagram to be $id$.

The identity section $e$ is defined by
\[ e: (pt \rightrightarrows pt)\to (\Z_2 \rightrightarrows
pt), \; \; e(1) =1, \] where 1 is the identity element in the
trivial group $pt$ and $\Z_2$.

The inverse $i$ is defined by
\[i:(\Z_2 \rightrightarrows pt) \to (\Z_2
\rightrightarrows pt), \;  \; i(a)=a^{-1},\] where $a\in \Z_2$. It
is a groupoid homomorphism because $\Z_2$ is commutative.

It is routine to check that  the above satisfies the axioms of
Weinstein groupoids. The local Lie groupoid associated to $B\Z_2$
is just a point. Therefore the Lie algebra of $B\Z_2$ is 0.
Moreover, notice that we have only used the commutativity of
$\Z_2$, so for any discrete commutative group $G$, $BG$ is a
Weinstein group with Lie algebra $0$. Moreover,
\[ \pi_0(BG)=1, \quad \pi_1(BG)=G. \]
So it still does not contradict with the uniqueness of the simply
connected and connected group integrating a Lie algebra.
\end{ep}

\begin{ep} [ ``$\Z_2 * B \Z_2$''] This is an example in which case Proposition
\ref{3asso} does not hold. Consider the groupoid $\Gamma=(\Z_2\times
\Z_2\rightrightarrows \Z_2)$. It is an action groupoid with
trivial $\Z_2$-action on $\Z_2$. We claim that the presented
\'etale differentiable stack $B\Gamma$ is a Weinstein group. We
establish all the structure maps on the presentation $\Gamma$.

The source and target maps are projections to a point.

The multiplication $m$ is defined by,
\[ m:\Gamma \times \Gamma \to \Gamma, \;\text{by} \; m((g_1, a_1),
(g_2, a_2))=(g_1 g_2, a_1 a_2).  \]It is a groupoid morphism
because $\Z_2$ (the second copy) is commutative. We have $m\circ
(m\times id) = m\circ (id \times m) $. But we can also construct a
non-trivial 2-morphism
\[\alpha: \Gamma_0(=\Z_2) \times \Gamma_0 \times\Gamma_0 \to
\Gamma_1, \;\text{by}\; \alpha ( g_1, g_2, g_3) =  ( g_1\cdot
g_2\cdot g_3, g_1\cdot g_2 \cdot g_3).\] Since the $\Z_2$ action
on $\Z_2$ is trivial, we have $m\circ (m\times id) = m\circ (id
\times m)\cdot \alpha $.

The identity morphism $e$ is defined by
\[ e: pt \rightrightarrows pt \to  \Gamma, \; \;
e(pt)=(1, 1),\] where 1 is the identity element in $\Z_2$.

The inverse $i$ is defined by
\[ i: \Gamma \to \Gamma, \;  \; i(g, a)=(g^{-1},
a^{-1}). \] It is a groupoid morphism because $\Z_2$ (the second
copy) is commutative.

It is not hard to check that $B\Gamma$ with these structure maps
is a Weinstein group. But when we look into the further
obstruction of the associativity described in Proposition
\ref{3asso}, we fail there. Let $F_i$'s be the six different ways
of composing four elements as defined in Proposition \ref{3asso}.
Then the 2-morphisms $\alpha_i$'s (basically coming from $\alpha$)
satisfy,
\[ F_{i+1} = F_i \cdot \alpha_i, \quad i=1, ..., 6 \;(F_7=F_1).\]
But $\alpha_i (1, 1, 1, -1)=(-1, -1)$ for all $i$'s except that
$\alpha_2 =id$. Therefore $\alpha_6 \circ \alpha_5\circ
...\circ\alpha_1(1,1,1,-1)=(-1, -1)$, which is not $id(1, 1, 1,
-1) = (-1, 1)$.
\end{ep}

\section{Contact groupoids}
\begin{ep}[Symplectic manifolds]
When $(M, \omega)$ is a symplectic manifold, the symplectic
groupoid $\Gamma_s(M)$ is the fundamental groupoid of $M$
\cite{cafe}. In this case, $Per_0(\omega)=Per(\omega)$, so
$P=Per(\omega) \times M$ is a trivial group bundle over $M$. $P$
is uniformly discrete if and only if $Per(\omega)$ is a discrete
group. Therefore $(M, \omega)$ is integrable as a Jacobi manifold
if and only if $\omega$ has discrete period group.

Suppose $(M, \omega)$ can be integrated into a contact groupoid.
Then according to the discussion above, the period group
$$Per(\omega)=a\cdot \Z,\;\;\;\;\; a\in \R.$$  To simplify the construction,
let us assume that $M$ is simply connected. Then $$\Gamma_s(M)=(M
\times M, (\omega, -\omega)).$$

When $a=0$, the contact groupoid $\Gamma_c(M)$ is simply
$\Gamma_s(M) \times \R$ and the groupoid structure is given by Theorem \ref{main2'}.

When $a\neq 0$, there is a principal $S^1$-bundle  $(E, \theta')$
over $(M, \omega/a)$. If $M$ is simply connected, $E$ is also
simply connected because $Per(\omega/a)=\Z$ and
\[ ...\lra\pi_2(M) \os{\partial= \int_{\gamma} \omega/a} {\lra}
\pi_1(S^1)\lra \pi_1(E) \lra 1. \]
 From this, we can get a principal
$\R/a\cdot\Z:=S^1_a$-bundle $(E, a\theta')$ over $(M, \omega)$.
$S^1_a$ acts diagonally on $E\times E$ and the 1-form $(a\theta',
-a\theta')$ is basic under this action, i.e. it is invariant under
the action and its contraction with the generator of the action is
0. So, $(a\theta', -a\theta')$ can be reduced to a 1-form $\theta$
on the quotient $E \times E/S^1_a$. Thus the contact groupoid
$\Gamma_c(M)$ is  $(E \times_{\R/a \cdot \Z} E, \theta)$.
\end{ep}

\begin{ep}
When the Jacobi manifold $M_0$ is contact with contact 1-form
$\theta_0$, the contact groupoid $\Gamma$ of $M_0$ is $M_0 \times
M_0 \times \R$, i.e. the direct sum of the pair groupoid and $\R$
with multiplication $(x, y, a ) \cdot (y, z, b) = (x, z, a+b) $.

The contact 1-form is $\theta= -(\exp \circ p_3) p_2^* \theta_0 +
p_1^* \theta_0$, where $p_i$, $1\leq i \leq 3$, is the projection
of $\Gamma$ to its $i$-th component. The function is $f= \exp
\circ p_3$.
\end{ep}

\begin{ep}[2-dimensional case]
Let $(M, \Lambda, E)$ be a 2-dimensional Jacobi manifold. Notice
that there is no multivector field in degree 3, so
\[ [\Lambda, \Lambda] = 2 \Lambda \wedge E =0 , \;\;\;\;\;
[\Lambda, E]=0 ,\] i.e. $M$ is a Poisson manifold equipped with a
vector field $E$ such that the Poisson structure is $E$-invariant.

It is known that every 2-dimensional Poisson manifold $M$ is
integrable \cite{cf2}. Actually, by \eqref{ses3} it is not hard to
see that the monodromy group $N_x(T^*M)=0$ because every
symplectic leaf of $M$ is either a point or 2 dimensional (so that
$\g_x(T^*M)=0$). By Lemma \ref{partial}, at point $x$ on a
symplectic leaf $ L$,
\[
\begin{split}
 N_x(T^*M\oplus_M \R) &= \{ (0, \int_{\gamma} \omega_L): \partial
 \gamma\}\\
 &= Per(\omega_L).
 \end{split}
 \]
Therefore, we have
\begin{cor}
A 2-dimensional Jacobi manifold $(M, \Lambda, E)$ is integrable if
and only if $Per(\omega_L)$ is discrete for all leaves $L$.
\end{cor}
\end{ep}

\begin{ep}[Non-integrable case]
As we can see in the last example,  there
are non-integrable Jacobi manifolds already in dimension 2. But as Poisson manifolds,
they are all integrable.

Let $M_a=\R^3$ be a Poisson manifold equipped with Poisson bracket
on coordinate functions $x^i$ as follows:
\[
\{x^2, x^3 \} = a x^1 ,\;\;\; \{x^3, x^1 \} = a x^2, \;\;\; \{x^1,
x^2 \}= a x^3,
\]
where $a=a(r)$ is a function depending only on the radius $r$.
Away from 0, the Poisson bivector field $\Lambda$ is given by
\[ \Lambda= (a\di x^1 +b x^1 r \bar{n}) \partial_2 \wedge
\partial_3 + c.p. , \]
where $\bar{n}= 1/r \sum_{i} x^i \di x^i$, $b(r) = a'(r)/r$, and
$c.p.$ is  short for cyclic permutation. The anchors $\rho_s:
T^*M_a \to TM_a$ and $\rho_c: T^*M_a \oplus_M\R \to TM_a$ are then
\[ \rho_s(\di x^i) = a v^i ,\;\;\;\;\; \rho_c((\di x^i, f))=av^i ,\;\;\; \forall
f\in \ci (M_a), \] where $v^1=x^3 \partial_2-x^2\partial_3$ etc.

Then the symplectic leaves of $M_a$ are spheres centered at the
origin (including the degenerate sphere: the origin itself).
Suppose $a(r)>0 $ for $r>0$. Choose sections  $\sigma_s: TM_a \to
T^*M$ and $\sigma_c: TM_a \to T^*M_a \oplus_{M_a}\R$ as follows:
\[ \sigma_s(v^i)= 1/a (\di x^i - \frac{x^i}{r} \bar{n}) ,
\;\;\;\;\; \sigma_c(v^i)= (\sigma_s(v^i), 0).\] Their compositions
with $\rho_s$ and $\rho_c$ are both identities. Then their
curvatures are
\[ \Omega_s= \frac{ra'-a}{a^2 r^3} \omega \bar{n}, \;\;\;\;\;
\Omega_c= \left( \Omega_s, \frac{r^2}{a} \omega \right), \] where
$\omega=x^1 \di x^2 \wedge \di x^3 + c.p. $. The symplectic form
on $S_r$ induced from the Poisson structure is $\frac{r^2}{a}
\omega$. Since $\int_{S_r} \omega =4 \pi r^3$, the symplectic area
of the sphere $S_r$, $A_a(r)$ is $4 \pi \frac{r}{a}$. By Lemma 3.6
in \cite{cf}, we have
\[ N_{\vec{x}}(T^*M_a)= \{ \int_{\gamma} \Omega_s, [\gamma] \in \pi_2(S_r)\} =A'_a(r)\Z \bar{n},\]
and \[ N_{\vec{x}}(T^*M_a \oplus_{M_a}\R)=\{ \int_{\gamma}
\Omega_c, [\gamma] \in \pi_2(S_r)\}=(A'_a(r) \Z \bar{n}, A_a(r)
\Z). \]

Generally, to measure the uniform discreteness of monodromy groups
$N_x(A)$ of some Lie algebroid $A$ over $M$, we introduce a
distance function $r_N(A)$ on $M$:
\[ r_N(A)(x) = \min_{0\neq \xi \in N_x(A)} \text{distance}(\xi,
0).\] $N_x(A)$ is uniformly discrete if and only if $r_N(A) >0$,
and $\lim_{y \to x} r_N(A)(y)>0$.

In our case, from the equation above, we have,
\[
r_N(T^*M_a)(x) = \begin{cases}
             +\infty& \text{if $r(x)=0$ or $A'_a(r)=0$},\\
             A'_a(r)& \text{otherwise},
             \end{cases}
\]
and
\[
r(T^*M_a\oplus_{M_a}\R)(x)=\begin{cases}
                        +\infty& \text{if $r(x)=0$}, \\
                        A'_a(r)+ A_a(r)& \text{otherwise}.
                        \end{cases}
\]
Therefore, $M_a$ is integrable as a Poisson manifold exactly when
$A'_a$ is nowhere 0 and $\lim_{r\to 0}A'_a(r) \neq 0$ or
$A'_a\equiv0$; $M_a$ is integrable as a Jacobi manifold exactly
when $\lim_{r\to 0}A'_a(r) +A_a(r)\neq 0$. By choosing a suitable
function $a$, (for example $a(r)=1/(\sin r +2)$),  it is easy to
discover an example $M_a$ which is integrable as  a Jacobi
manifold but not as a Poisson one.
\end{ep}



\lhead[\fancyplain{}{}]{\fancyplain{}{Bibliography}}
\rhead[\fancyplain{}{}]{\fancyplain{}{}}

\bibliographystyle{abbrv}
\bibliography{bibz.bib}

\end{document}